\documentclass[11pt]{article}

\usepackage{amsmath, amsfonts, amssymb,latexsym,wasysym,graphicx,stmaryrd,tikz,txfonts}
\input epsf.sty

\usetikzlibrary{decorations.pathmorphing}
\usetikzlibrary{arrows}

\newtheorem{Theorem}{Theorem}[section]
\newtheorem{Lemma}{Lemma}[section]
\newtheorem{Proposition}[Lemma]{Proposition}
\newtheorem{Corollary}[Lemma]{Corollary}
\newtheorem{Definition}[Lemma]{Definition}

\newtheorem{Assumption}[Lemma]{Assumption}

\newcommand{\BEQ}{\begin{equation}}     
\newcommand{\BEA}{\begin{eqnarray}}
\newcommand{\BD}{\begin{displaymath}}
\newcommand{\EEQ}{\end{equation}}       
\newcommand{\EEA}{\end{eqnarray}}
\newcommand{\ED}{\end{displaymath}}

\newcommand{\del}{\delta}
\newcommand{\Del}{\Delta}
\newcommand{\eps}{\varepsilon}          

\newcommand{\supp}{{\mathrm{supp}}}

\newcommand{\Tr}{{\mathrm{Tr}}}

\newcommand{\Vol}{{\mathrm{Vol}}}

\newcommand{\osc}{{\mathrm{osc}}}
\newcommand{\locsup}{{\mathrm{locsup}}}
\newcommand{\Locsup}{{\mathrm{Locsup}}}
\newcommand{\R}{\mathbb{R}}

\newcommand{\Z}{\mathbb{Z}}
\newcommand{\N}{\mathbb{N}}
\newcommand{\D}{\mathbb{D}}

\def\proba{{\mathbb{P}}}

\def\esper{{\mathbb{E}}}

\def\Var{{\mathrm{Var}}}

\def\Errfc{{\mathrm{Errfc}}}
\newcommand{\lele}{\preceq}

\newcommand{\eop}{\hfill $\Box$}        

\newcommand{\half}{{1\over 2}}          
\renewcommand{\vec}[1]{\boldsymbol{#1}} 

\catcode`\@=11
\def\numberbysection{\@addtoreset{equation}{section}
        \def\theequation{\thesection.\arabic{equation}}}
\numberbysection

\begin{document}

\vspace*{1.5cm}
\begin{center}
{\Large \bf PDE estimates for multi-dimensional KPZ equation.}

\end{center}

\vspace{2mm}
\begin{center}
{\bf  J\'er\'emie Unterberger$^b$}
\end{center}

\vskip 0.5 cm

\centerline {$^b$Institut Elie Cartan,\footnote{Laboratoire 
associ\'e au CNRS UMR 7502} Universit\'e de Lorraine,} 
\centerline{ B.P. 239, 
F -- 54506 Vand{\oe}uvre-l\`es-Nancy Cedex, France}
\centerline{jeremie.unterberger@univ-lorraine.fr}

\vspace{2mm}
\begin{quote}

\renewcommand{\baselinestretch}{1.0}
\footnotesize
{We study in this series of articles the Kardar-Parisi-Zhang  (KPZ) equation
$$ \partial_t h(t,x)=\nu\Del h(t,x)+\lambda  V(|\nabla h(t,x)|) +\sqrt{D}\, \eta(t,x), \qquad x\in\R^d $$
in $d\ge 1$ dimensions. The forcing term $\eta$ in the right-hand side is
 a regularized white noise. The deposition rate $V$ is assumed to be isotropic and convex.
Assuming $V(0)\ge 0$, one finds $V(|\nabla h|)\ltimes |\nabla h|^2$ for small gradients,
 yielding the equation which is  most commonly used  in the literature.

\medskip
The present article is dedicated to existence results and PDE estimates for the solution. Our results extend in a non-trivial way  those previously obtained for the noiseless
equation. We prove in particular a comparison principle for sub- and supersolutions of the KPZ equation in 
new functional spaces containing unbounded functions, implying existence and
uniqueness. These new functional spaces made up of functions with "locally bounded
averages", generically called ${\cal W}$-spaces thereafter, and which may be of  interest for the  study
of  parabolic equations in general, allow {\em local} or {\em pointwise}  estimates.  The comparison to the linear heat equation through a
Cole-Hopf transform is an essential ingredient in the proofs, and our results are accordingly valid only for
a function $V$ with at most quadratic growth at infinity. 
}
\end{quote}

\vspace{4mm}
\noindent

 \medskip
 \noindent {\bf Keywords:}
  KPZ equation, viscous Hamilton-Jacobi equations, maximum principle, renormalization,
  multi-scale analysis.

\smallskip
\noindent
{\bf Mathematics Subject Classification (2010):}  35B50, 35B51, 35D40, 35K55, 35R60,  35Q82, 60H15,  81T08, 81T16, 81T18, 
82C41.

\newpage
\tableofcontents


\section{General introduction}

The KPZ equation \cite{KPZ} is a stochastic partial differential equation describing the growth by normal deposition of an interface  in $(d+1)$ space dimensions, see e.g. \cite{BarSta,Car}. 
By definition the time evolution of the height $h(t,x)$, $x\in\R^d$, is given by
\BEQ \partial_t h(t,x)=\nu\Del h(t,x)+2\lambda  \left( \sqrt{1+|\nabla h(t,x)|^2}-1\right) +\sqrt{D}\, \eta(t,x), \qquad x\in\R^d \label{eq:KPZ0} \EEQ
 where $\nu$ ({\em diffusion constant}), $\lambda$ ({\em coupling constant}) are
 positive constants, and $\eta(t,x)$ is some (possibly regularized) white noise.
 The gradient $|\nabla h|$ (the slope of the interface) is assumed to remain throughout small so that the evolution makes physically sense, precluding
e.g. any overhang, so that the non-linear term $\sqrt{1+|\nabla h(t,x)|^2}-1\simeq \half |\nabla h(t,x)|^2$ is essentially quadratic; using this approximation gives the most common
form of this equation in the literature, thereafter called {\em quadratic KPZ equation},

\BEQ \partial_t h(t,x)=\nu\Del h(t,x)+\lambda |\nabla h(t,x)|^2 +\sqrt{D}\, \eta(t,x), \qquad x\in\R^d \label{eq:KPZ2} . \EEQ

 Following these preliminary remarks, we shall call  KPZ equation any equation of the type
\BEQ \partial_t h(t,x)=\nu\Del h(t,x)+\lambda   V(\nabla h(t,x))) +\sqrt{D}\, \eta(t,x), \qquad x\in\R^d \label{eq:KPZ} \EEQ
where the {\em deposition rate} $V$ is  isotropic and convex (hence $V(\nabla h(t,x))=a+b|\nabla h(t,x)|^2+\ldots$ around $0$, with $b\ge 0$). The interest is generally in the large-scale limit of this equation, for $t$ large. 
A well-known naive rescaling argument gives some ideas about the dependence on the dimension of this limit. 
 Namely, the linearized equation, a stochastic heat equation which is a particular instance of {\em Ornstein-Uhlenbeck
 process}, 
\BEQ \partial_t \phi(t,x)=\nu \Del \phi(t,x)+ \sqrt{D}\,  \eta(t,x), \quad (t,x)\in\R_+\times\R^d \label{eq:O.-U.} \EEQ
 is invariant under the  rescaling $\phi(t,x)\mapsto \phi^{\eps}(t,x):= \eps^{-d_{\phi}}
\phi(\eps^{-1} t,\eps^{-\half} x)$, where
\BEQ d_{\phi}:= \half(\frac{d}{2}-1) \EEQ
is the {\em scaling dimension} of the field $\phi$ (or rather {\em half} the scaling dimension, in the
physicists' convention); we used here the equality in distribution,
$\eta(\eps^{-1} t,\eps^{-\half}x)\stackrel{(d)}{=} \eps^{\half(1+\frac{d}{2})}\eta(t,x)$. 
Assuming that $\phi$ is a solution of the KPZ equation instead (say, with quadratic
deposition rate $|\nabla \phi|^2$) yields after rescaling 
\BEQ \partial_t \phi^{\eps}(t,x)=\nu \Del \phi^{\eps}(t,x)+ \eps^{d_{\phi}} \frac{\lambda}{2}|\nabla \phi^{\eps}(t,x)|^2 + \sqrt{D}\,  \eta^{\eps}(t,x), \label{eq:power-counting-vertex-2} \EEQ
where $\eta^{\eps}\stackrel{(d)}{=}\eta$.
 For $d>2$,  the scaling exponent $d_{\phi}$ is $>0$, and the non-linear term scaling coefficient, $\eps^{d_{\phi}}$, vanishes in the limit $\eps\to 0$; in other terms, 
the KPZ equation is sub-critical at large scales in $\ge 3$ dimensions and believed to behave like the corresponding linearized equation up to a redefinition
(called {\em renormalization}) of the {\em diffusion constant} $\nu$ and of the {\em noise strength} $D$.  More
precisely, according
to the general scheme due to K. Wilson \cite{Wilb,Wilc}, the fluctuations of the solution field at time scale
of order $\eps^{-1}\approx 2^j$ and space scale of order $\eps^{-\half}\approx 2^{j/2}$ should be approximately
governed for $j$ large by a linearized equation with scale-dependent coefficients $\nu^{(j)},D^{(j)}$ $(j\ge 0)$, themselves  solutions of 
a certain complicated but explicit discrete dynamical system. Ultimately our purpose
is to confirm these predictions.

\bigskip

\noindent The present work contains some preliminary steps towards this goal,
using {\em deterministic tools}.
Since we cannot capture the large-scale behaviour of the equation without taking
into account  white noise fluctuations,  we do not address the full equation
(\ref{eq:KPZ}) but either (i) the associated homogeneous equation $(D=0$); or (ii)
a KPZ-type equation with general, deterministic right-hand side $g$ (possibly coming
from a realization of white noise, which allows a connection to our original problem), exhibiting an extra
scale-dependent linear damping, which is supposed to mimick the behaviour of the
KPZ solution {\em at some given scale}. Because of the damping we do not see the dimension dependence, so our results actually hold for any $d\ge 1$. For both equations we provide estimates which
are essentially optimal, reproducing the expected Ornstein-Uhlenbeck scaling in case (ii).
Thus our renormalization scheme will ultimately be able to include some of these a priori estimates. The connection to the multi-scale analysis of KPZ equation is explained in some details in the end of the article (sections 5 and 6). One may however choose to
disregard these matters, and see this article as a purely deterministic PDE paper
concerned about existence results and PDE estimates for inhomogeneous viscous Hamilton-Jacobi equations
of a certain type. Our original contribution in this respect is that we want to  allow right-hand 
sides in functional spaces large enough to contain realizations of regularized white noise. Since the latter has unbounded fluctuations in full space, standard existence
theory, mostly based on the maximum principle, cannot be applied. Thus we are led
to solve KPZ equation in new {\em  functional spaces}, {\em modeled} on the {\em space ${\cal H}^0$ of functions with locally bounded averages} (see below). Instead of 
the a priori bounds in supremum norm obtained from the maximum principle, one gets 
{\em poinwise} or {\em local} a priori bounds in {\em pointwise} or {\em local quasi-norms}, using some stronger and more versatile version of the maximum principle for parabolic equations, based
on the comparison to the heat flow.

\bigskip

\noindent While we provide an outline of the article, we shall try to explain
more concretely the above principles.

\medskip

\noindent  Section 2 is concerned with {\em bounded} solutions of the homogeneous
KPZ equation, relying in particular on the comparison principle for non-linear parabolic PDE's. (Precise assumptions for the deposition rate $V$ are listed in section 1). The titles of subsections 2.1, 2.2, 2.3 reflect the three main arguments from
which estimates can be derived; comparison to the linear heat equation (\S 2.1) is the main argument surviving in later sections when we consider unbounded solutions. Some results are
derived with little effort from those already existing in the literature; on the other hand, the bounds on the
gradient and on the higher derivatives of the solution, see Theorems \ref{th:decay-nabla-L1} and \ref{th:bound-higher-derivatives}, may not be
found elsewhere.

\medskip

\noindent The really original material starts in section 3 with the search for solutions in  new spaces
of possibly unbounded functions with good averaging properties. There is a large
variety of choice for such spaces, for which we therefore coin a generic term, "${\cal
W}$-spaces" for the discussion. Generally speaking all ${\cal W}$-spaces are modeled
after 
\BEQ {\cal H}^0:=\{f\in L^{\infty}(\R^d)\ |\ \forall x\in\R^d, f^*(x)<\infty\}, \qquad f^*(x):=\sup_{\tau>0} e^{\tau\Del}|f|(x).\EEQ
Since $e^{\tau\Del}|f|(x)$ is some weighted average of $|f|$ centered at $x$, it makes
sense to speak of elements of ${\cal H}^0$ as {\em functions with locally bounded
averages}. Clearly, $L^{\infty}\subset {\cal H}^0$, but  unbounded
 functions, with arbitrarily large but rare fluctuations,  also belong to ${\cal H}^0$, notably our regularized white noise, $\eta$, for
 which (as is well-known for the supremum of $n$ essentially independent identicallly
 distributed Gaussian
 variables) $\sup_{|x|\le n} |\eta(x)|=O(\sqrt{\log(n)})$. By construction, the solution $f$ of the heat equation $(\partial_t-\Del)f(t,x)=0$  satisfies $|f(t,x)|\le (f_0)^*(x)$, a {\em pointwise} version of the maximum
principle which states that $||f||_{\infty}\le ||f_0||_{\infty}$. Clearly one also has
$(f_t)^*(x)\le (f_0)^*(x)$. Let now $h=h(t,x)$ be a solution of the {\em homogeneous} quadratic KPZ equation (\ref{eq:KPZ2}).  Since Cole-Hopf transformation $h\mapsto e^{\lambda h}$ 
maps solutions  of (\ref{eq:KPZ2})  into
solutions of the heat equation, we get  $(e^{\lambda h_t})^*(x)\le 
(e^{\lambda h_0})^*(x)$. With some extra work (see Lemma \ref{lem:easy}), letting
 \BEQ |||f|||_{{\cal H}^{\lambda}}(x):=\frac{1}{\lambda} \ln \left( (e^{\lambda |f|})^*(x) \right)  \label{intro:f-lambda} \EEQ
one proves:
\BEQ |||h_t|||_{{\cal H}^{\lambda}}(x)\le |||h_0|||_{{\cal H}^{\lambda}}(x).  \label{intro:hth0} \EEQ
 We
have thus defined a new space,
 ${\cal H}^{\lambda}:=\{f\in L^1(\R^d) \ |\ \forall x\in\R^d, |||f|||_{{\cal H}^{\lambda}}(x)<\infty\}$, together with what plays the r\^ole of a family of {\em "pointwise quasi-norms"},
 $|||\ \cdot\ |||_{{\cal H}^{\lambda}}(x)$. The interplay between these spaces is investigated in \S 3.2, where we show in particular satisfactory collective properties of the family (\ref{intro:f-lambda}) with respect to convex operations, justifying the term
 of "quasi-norms" for lack of a better term. 
 
\medskip 
 \noindent
Then the rest of section 3, resp. section 4, are dedicated to existence theorem and
 bounds  of the homogeneous, resp. inhomogeneous KPZ equation in ${\cal W}$-spaces.
Let us first discuss the {\em homogeneous} case. We say that $h=h(t,x)$ solves the  {\em homogeneous KPZ equation} if
\BEQ (\partial_t-\nu\Del)h(t,x)=\lambda V(\nabla h). \label{intro:homogeneous} \EEQ
 We prove  a {\em comparison principle}
  for sub- and supersolutions of the homogeneous KPZ equation in these spaces,
implying existence and unicity for viscosity solutions, which are proved to be classical. The
statement is as follows (see Theorem \ref{th:comparison}):

\bigskip

\noindent {\bf Theorem 1 (comparison principle).} {\em 
Let $\underline{U}\in USC([0,T]\times\R^d)\cap{\cal H}^{2\lambda}([0,T])$ (resp. $\bar{U}\in LSC([0,T]\times\R^d)\cap {\cal H}^{2\lambda}([0,T])$) be a viscosity sub-solution (resp. super-solution)
of  the homogeneous KPZ equation (\ref{eq:1.1}). 
Then $\underline{U}\le \bar{U}$ in $[0,T]\times\R^d$.}

\bigskip

\noindent Bounding the gradient of the solution turns out to be more challenging than getting the almost
trivial bound (\ref{intro:hth0}). One possibility (see \S 3.3 for a discussion) is to introduce {\em local
${\cal W}$-spaces}, ${\cal W}^{1,\infty;\lambda}_j$, for which one replaces the various ${\cal W}$-{\em "pointwise quasi-norms"} by stronger ${\cal W}$-{\em "local quasi-norms"}, $|||\ \cdot\ |||_{{\cal W}^{1,\infty;\lambda}_j}(x)$,  obtained by substituting to $|f(x)|$ its 
local supremum $\sup_{B(x,1)} |f|$ or more generally (in consistence with
section 4, see below) $\sup_{B(x,2^j)} |f|$, where $j$ is some scale. As shown in
\S 3.3, the finiteness of the {\em local} quasi-norm implies a polynomial bound
at infinity of a precise order (which holds for $\eta$ and even for $e^{\lambda|\eta|}$ !) We emphasize that we do also get bounds in {\em local} quasi-norms for the solution
 in terms of the  {\em local} quasi-norm of the initial condition (see discussion in \S 3.3 and after the proof of Lemma \ref{lem:classical} in \S 3.4), so using {\em local} quasi-norms (here and also in the non-homogeneous case treated below) only shortens
 the statements, to the great happiness of the reader, while restraining the generality.  Our main
result is :

\bigskip

\noindent {\bf Lemma \ref{lem:classical} (bound for the homogeneous KPZ equation)} {\em
Let $h$ be the solution of the homogeneous KPZ equation (\ref{intro:homogeneous}) with
$h_0\in {\cal W}^{1,\infty;2\lambda}\cap C^2$. Then $h_t\in {\cal W}^{1,\infty;2\lambda/5}$ and 
\BEQ  |||\, 2^{j/2}\locsup^j |\nabla h_t|\, |||_{{\cal H}^{2\lambda/5}}\le  5 |||h_0|||_{{\cal W}^{1,\infty;2\lambda}}(x)
.\EEQ
}

\bigskip

 \noindent An intelligent study of
the  {\em  full inhomogeneous} equation (\ref{eq:KPZ}) for large time is a much more difficult problem, since it relies in an essential way on the averaging properties of the noise. However,  essentially optimal bounds can be obtained for the {\em scale $j$ infra-red
cut-off equation},
\BEQ \partial_t\psi=\nu\Del \psi-2^{-j} \psi + \lambda V(\nabla \psi)+g, \label{intro:j} \EEQ
where $g$ is some adequate, regular right-hand side. 
  This equation in meant to select the fluctuations of the solution on time ranges of order $2^j$ and space ranges
of order $2^{j/2}$. Thus $g$ should enjoy the same scaling properties as 
the "$j$-th scale projected" regularized white noise $\eta$. Scalings are
discussed in details in section 5; let us just mention at this point (see Remark at the very end of section 5) that only 
smaller scale components $j'=0,\ldots,j-1$ of the right-hand side need to be discarded to get a correct scaling.  Then we show in
section 4 how to solve  and bound  (\ref{intro:j}) along essentially the same lines as in section 3. In the course of the computations we are led to introduce
new ${\cal W}$-spaces, which take into account both the scaling, and the time-dependence
(for the right-hand side). The main result (see Lemma \ref{lem:bound-psin} and
ensuing discussion) is (see \S 4.2 for the definition of the family of  time-dependent 
${\cal W}$-spaces ${\cal W}_j^{1,\infty;\lambda}([0,t])$ and the associated
"local quasi-norms" $||| \ \cdot\ |||_{{\cal W}_j^{1,\infty;\lambda}([0,t])}(x)=
|||\ \cdot\ |||_{\lambda,j}([0,T],x)$).

\bigskip

\noindent{\bf Theorem 2 (a priori bounds for the KPZ equation)}    {\em 
Let 
 $\psi$  be the (viscosity) solution of (\ref{intro:j}) with
initial condition $\psi_0\in{\cal W}^{1,\infty;2\lambda'}_{j}\cap  C^2$ with
$\lambda'>\lambda$, and forcing term $g\in {\cal W}_j^{1,\infty;2\lambda}([0,T])\cap C([0,T],
C^3(\R^d))$. Then 
\BEQ  |||\, \locsup^j \psi_t \, |||_{{\cal H}^{\lambda}}(x) \le 
e^{-2^{-j}t} \,  |||\locsup^j \psi_0|||_{{\cal H}^{\lambda'}}(x) + 
 |||\, \locsup^j g|||_{\lambda,j}([0,t],x) \label{intro:thm2-1} \EEQ
and 
\BEQ |||\, 2^{j/2} \locsup^j |\nabla\psi_t|\, |||_{{\cal H}^{2\lambda/5}}(x) \le 
5 \left( |||g|||_{{\cal W}_j^{1,\infty;2\lambda}([0,t])}(x) + e^{-2^{-j}t} |||\psi_0|||_{{\cal W}_j^{1,\infty;2\lambda'}}(x) \right). \label{intro:thm2-2} \EEQ
}

The reader may easily check that, {\em without the damping term}, (\ref{intro:thm2-1},
\ref{intro:thm2-2}) hold with some modifications on the time interval $[0,2^j]$,
\BEQ   |||\, \locsup^j \psi_t \, |||_{{\cal H}^{\lambda}}(x) \le C
\left( \,  |||\locsup^j \psi_0|||_{{\cal H}^{C\lambda}}(x) + 
 |||\, \locsup^j g|||_{C\lambda,j}([0,t],x)  \right) \EEQ
and 
\BEQ |||\, 2^{j/2} \locsup^j |\nabla\psi_t|\, |||_{{\cal H}^{\lambda}}(x) \le C
 \left( |||g|||_{{\cal W}_j^{1,\infty;C\lambda}([0,t])}(x) +  |||\psi_0|||_{{\cal W}_j^{1,\infty;C\lambda}}(x) \right)  \EEQ
for $C$ large enough. Composing these estimates on the successive time intervals
$[n2^j,(n+1)2^j]$, $n=0,1,\ldots$, one may easily prove a bound for the solution and its
gradient in $|||\ \cdot\ |||_{{\cal H}^{a\lambda}}$ "quasi-norms" for any $a\ge 1$, provided $\psi_0\in \cup_{a\ge 1}
{\cal H}^{a\lambda}$ and $g\in \cup_{a\ge 1} {\cal W}_j^{1,\infty;a\lambda}$,
implying in particular {\em global existence of the solutions of the full inhomogeneous KPZ equation (\ref{eq:KPZ}) in $\cal W$-spaces}.  However, because of  the "loss of regularity"  in $\lambda$, bounds increase exponentially in time and become extremely
bad for $t$ large.

\bigskip

 \noindent Finally sections 5 and 6 are appendices containing
multi-scale decompositions of the propagator $(\partial_t-\nu\Del)^{-1}$  and the
white noise $\eta$, and large-deviation estimates for $\eta$, implying the applicability
of the general arguments developed in  section 4 to the case of the noisy KPZ equation (\ref{eq:KPZ}). The proof of large-deviation estimates in itself is far from trivial
because we need to bound the ${\cal H}^{\lambda}$ "quasi-norm" of $\eta$ (see
(\ref{intro:f-lambda})), which involves its {\em exponential}. Since exponentiated
 Gaussian variables do not admit any exponential moment, we must turn to non-standard
(and not that well-known) deviation estimates found in the Soviet literature of the 60es
(see section 6). Though section 5 is really helpful to motivate the scaling issues
related to (\ref{intro:j}, the reader who is not particularly interested in stochastic
PDEs may safely skip section 6, which is quite involved and of a very different nature with respect to the
previous ones.

\bigskip

\noindent In a companion article \cite{Unt-KPZ2}, we tackle the problem of getting existence/unicity and estimates in ${\cal W}$-quasi norms for solutions of the scale $j$ infra-red cut-off equation (\ref{intro:j}), but with totally different techniques, using the
{\em Hamilton-Jacobi-Bellman formalism}; recall this formalism allows to
represent the solution $\psi$ as the maximum over an admissible class of random paths $X$ driven by Brownian motion of a functional $\int_0^t F(s,X_s)\, ds$. Controlling the
random characteristics allows  {\em less precise} but {\em much more flexible} estimates, extending 
in particular to the case of  deposition rates $V$ growing faster than quadratically at infinity.
It is interesting to compare the results obtained by the two methods. The two
articles are largely dependent one from the other, though the present article may
be read independently from the other {\em except} for a technical point in the proofs:
the {\em unicity statement} in Theorem 2, which we could {\em not} prove by
the techniques developed here.


\section{Model and notations}


We consider throughout the present article either the {\em homogeneous} (or {\em noiseless})  equation
\BEQ \partial_t h=\nu\Del h+\lambda V(\nabla h) \label{eq:1.1}\EEQ
where $\lambda>0$ is a fixed, arbitrary constant,
or the {\em infra-red cut-off, inhomogeneous} equation,
\BEQ \partial_t \psi=\nu\Del \psi-\eps \psi+ \lambda V(\nabla \psi)+g, \label{eq:1.2}\EEQ
where the constant $\eps=2^{-j}$ ($j\ge 0$) is an infra-red cut-off of scale $j$. 
Bounds for the homogeneous equation (\ref{eq:1.1}), resp. inhomogeneous equation
(\ref{eq:1.2}) turn out to be quite different in the end, though they are based of
course on the same principles, so -- in order to avoid any confusion --  we keep throughout the article to the following convention: {\bf solutions of the {\em homogeneous} equation are denoted by $h$, solutions of the {\em inhomogeneous} equation by $\psi$. }

\medskip
 \noindent The assumptions on $V$ are the following:

\begin{Assumption} \label{assumptions}
{\bf The deposition rate $V$ satisfies the following assumptions,
\begin{itemize}
\item[(1)] $V$ is $C^2$;
\item[(2)] $V$ is isotropic, i.e. $V(\nabla h)$ is a function of $y=|\nabla h|$; by abuse of notation we shall
consider $V$ either as a function of $\nabla h$ or of $y$;
\item[(3)] $V$ is convex; 
\item[(4)] $V(0)=0$ and  $V(y)\ge 0$ for all $y\ge 0$;
\item[(5)] (quadratic growth at infinity) $V(y)\le y^2$ for all $y\ge 0$.
\end{itemize} }
\end{Assumption}

It follows immediately from Assumptions (1), (2) and (4) that $V(y)=O(y^2)$ near $y=0$.  Assumption (5) is  thus equivalent (up to a redefinition of the constant $\lambda$) to requiring that $V$ has at most quadratic growth at infinity.

As for the force term $g$, it is assumed to be regular enough and have good averaging properties, depending on
the cut-off scale $j\approx -\log\eps$; the regularized white noise  $\eta$  (as shown in  section 6) satisfies
these properties.

\medskip

Assumption (3)  is a key assumption to get a time decay of the
gradient of the solution, and is also used in the proof of the comparison theorem for unbounded solutions;  Proposition \ref{prop:decay-gradient} (ii), (iii) hold under a stronger assumption.
Assumption (5) allows a comparison of the solutions to those of  the usual KPZ equation corresponding to $V(y)=y^2$, which is linearizable.

\medskip
{\bf Notations.}  The notation:
$f(u)\lesssim g(u)$, resp. $f(u)\gtrsim g(u)$ means: $|f(u)|\le C|g(u)|$, resp. $|f(u)|\ge C|g(u)|$, where $C>0$ is an unessential constant (depending only on the dimension $d$
and on the coefficients of the linearized equation, $\nu$ and $D$). 
Similarly, $f(u)\approx g(u)$ means: $f(u)\lesssim g(u)$ and $g(u)\lesssim f(u)$.
We denote by $L^p$, $p\in[1,\infty]$  the usual Lebesgue spaces with associated norm $||\ ||_{p}$,  by ${\cal W}^{1,\infty}$ the Sobolev space of bounded functions with bounded generalized derivative, and by $C^{1,2}$ the space of functions which are $C^1$ in time and $C^2$ in space. The positive, resp. negative part of a function $f$ is denoted
by $f_+$, resp. $f_-$; by definition, $f^+,f^-\ge 0$, $f=f^+-f^-$ and $f^+ f^-=0$. The {\em oscillation} $\osc_{\Omega} f$
of a continuous function $f$ on a domain $\Omega$ is defined as  $\sup_{\Omega}f-\inf_{\Omega}f$; the {\em average} 
$\frac{1}{|\Omega|} \int_{\Omega} f$ of $f$
on a bounded domain $\Omega$ is denoted by $\fint_{\Omega}f$. The space of
lower, resp. upper semicontinuous functions on a domain $\Omega$ is denoted by $LSC(\Omega)$, resp. $USC(\Omega)$.


\section{Bounds for the homogeneous equation: the case of a bounded initial condition}


We consider in this section the homogeneous equation,
\BEQ \partial_t h=\nu\Del h+\lambda V(|\Del h|)  \label{eq:KPZdet} \EEQ
with initial condition $h_0(x)=h(0,x)$ in ${\cal W}^{1,\infty}$.  One finds in the literature a detailed
study of the particular case $V(y)=y^q$, $q>1$. Most basic results (including existence), based on the
principle of maximum or on a short-time series expansion of the mild solution, depend very little
on the precise form of $V$, provided it is regular enough and, say, polynomially bounded. We quickly review them
now and leave it to the reader to check that they extend to a rate $V$ satisfying Assumptions \ref{assumptions} (1),
(2), (4). 

By \cite{AmoBen} and \cite{BenLau}, 
the Cauchy problem has a unique, global solution $u$ which is classical for positive times, that is,
 $u\in C([0,+\infty)\times\R^d)\cap C^{1,2}((0,\infty)\times\R^d)$.
The comparison principle, in the form proved by Kaplan \cite{Kap} for classical, bounded solutions of
non-linear parabolic equations on unbounded spatial domains, implies that $h_t\ge 0$ for all
$t\ge 0$ (resp. $h_t\le 0$ for all $t\ge 0$) if
$h_0\ge 0$  (resp. $h_0\le 0$) and  yields the  a priori estimates
\BEQ ||h_t||_{\infty}\le ||h_0||_{\infty},\qquad ||\nabla h_t||_{\infty}\le ||\nabla h_0||_{\infty} \qquad 
(t\ge 0). \label{eq:a-priori-estimates}\EEQ

\medskip

We now prove time-decay estimates of the solution for various norms, emphasizing those which are not a
straightforward extension of previously known results for $V(y)=y^q$. Such estimates
 come roughly from three different sources, and are correspondingly split
 into 3 paragraphs (\S 2.1, \S 2.2 and \S 2.3). Generally speaking,  constants appearing in the inequalities deteriorate
 when $\nu\to 0$ whenever parabolic estimates are involved  (see below); Proposition \ref{prop:decay-gradient}
 (ii), (iii) is an outstanding exception.

\medskip

We recall here briefly for non-specialists the {\em maximum principle } and the {\em comparison principle} 
for parabolic PDE's, in a weak
form which is sufficient for section 2. Standard references on the subject are e.g. \cite{Eva}, \cite{CraIshLio},
\cite{Bar}.

\begin{Proposition}[maximum and comparison principle]
Let $u(t,x)$, $(t,x)\in[0,T]\times\R^d$ be a classical solution of the parabolic PDE $\partial_t u(t,x)=\nu\Del u(t,x)+
W(t,x,\nabla u(t,x))$, where $W$ is a smooth function, bounded in any subset of the form $\R\times\R^d\times K$,
$K\subset\R^d$ compact. Assume that $\sup_{[0,T]\times\R^d} |u|<\infty$ and $\sup_{[0,T]\times\R^d}
|\nabla u|<\infty$. Then:
\begin{itemize}
\item[(i)] (weak maximum principle) $\forall t\in[0,T],\ ||u_t||_{\infty}\le ||u_0||_{\infty}$.
\item[(ii)] (weak comparison principle) let $\bar{U}$, resp. $\underline{U}$ be a super-, resp. sub-solution of the above PDE, namely, $\bar{U},\underline{U}\in C^{1,2}([0,T]\times\R^d)$ and 
\BEQ \partial_t \bar{U}(t,x)\ge \nu\Del \bar{U}(t,x)+
W(t,x,\nabla \bar{U}(t,x)),\ \partial_t \underline{U}(t,x)\le\nu\Del \underline{U}(t,x)+
W(t,x,\nabla \underline{U}(t,x)).
\EEQ
 Assume $\bar{U}_0\ge\underline{U}_0$. Then $\bar{U}_t\ge\underline{U}_t$ for all $t\ge 0$.
\end{itemize}
\end{Proposition}

Note that the above proposition extends under appropriate monotonicity hypotheses to parabolic PDE's of the form
$\partial_t u(t,x)=\nu\Del u(t,x)+
W(t,x,u(t,x),\nabla u(t,x))$. However, it is precisely the absence of dependence of $W$ on $u(t,x)$ that makes
two-sided a priori estimates like (\ref{eq:a-priori-estimates}) so easy.


\subsection{Comparison to the linear heat equation}


Assumptions \ref{assumptions} (4)-(5), $0\le V(y)\le y^2$, allows (as we shall presently see) a direct 
{\em comparison with the linear heat equation} if either  $h_0\ge 0$
or  $h_0\le 0$. Bounds for signed initial conditions  follow then from the comparison principle: namely, letting
$\bar{h}$, resp. $\underline{h}$ be the solution of (\ref{eq:KPZdet}) with initial condition $h^+_0$, resp. $-h^-_0$,
one has
\BEQ \underline{h}\le 0, \bar{h}\ge 0; \qquad \qquad \underline{h}\le h\le \bar{h}.\EEQ 
Also, $t\mapsto ||\bar{h}_t||_1$ is increasing, while $t\mapsto ||\underline{h}_t||_1$ is decreasing.

Considering first $\underline{h}$,  the comparison
principle allows one to bound the solution of (\ref{eq:KPZdet}) by the solution of the linear heat equation
with same initial condition, namely, 
\BEQ |\underline{h}(t)|\le e^{t\nu\Del}h_0^-. \label{eq:h-}\EEQ
We now turn to $\bar{h}$ and bound similarly the solution of (\ref{eq:KPZdet}) with positive initial condition by
the solution $u$ of the standard KPZ equation, $\partial_t u=\nu\Del u+\lambda |\nabla u|^2$ with the same initial
condition. The exponential 
transformation $w:=e^{\frac{\lambda}{\nu}u}-1$ turns it into the linear equation $\partial_t w=\nu\Del w$, with
positive initial condition $w_0=e^{\frac{\lambda}{\nu}h^+_0}-1$. The inequality $x\le \frac{\nu}{\lambda}(e^{\frac{\lambda}{\nu}x}-1),\ x\ge 0$ yields
\BEQ ||\bar{h}_t||_{\infty}\le||u_t||_{\infty}\le \frac{\nu}{\lambda} ||w_t||_{\infty}.\EEQ

\medskip
To go further, we assume $w_0\in L^1$ and use the following standard parabolic estimates \cite{Sch} for $q=1$.

\begin{Proposition}[parabolic estimates] There exist constants $C_k$, $k=0,1,\ldots$ depending only on $d$ such that,
for every regular enough function $f_0:\R^d\to\R$ and $p\ge q\ge 1$,
\BEQ ||\nabla^k e^{t\nu\Del}f_0||_{p}\le C_k (\nu t)^{-\frac{d}{2}(\frac{1}{q}-\frac{1}{p})-\frac{k}{2}}||f_0||_q,\qquad k\ge 0. 
\label{eq:parabolic-estimates} \EEQ
\end{Proposition}

Let $\mu$ be the Lebesgue measure on $\R^d$. The well-known identity $\int f(u(x))dx=\int_0^{+\infty} \mu(u>a)f'(a)da$, valid for  $u:\R^d\to\R_+$ measurable and $f:\R_+\to\R$ smooth such that $f(0)=0$, 
 yields for $f(u)=\frac{\nu}{\lambda}(e^{\frac{\lambda}{\nu}u}-1)$
\BEA  ||\bar{h}_t||_{1}&\le & \frac{\nu}{\lambda} ||w_t||_1\le \frac{\nu}{\lambda} ||w_0||_1 \nonumber\\
 &\lesssim& \int h^+_0(x) {\bf 1}_{h^+_0(x)\le \nu/\lambda} dx + 
 \int_{\nu/\lambda}^{+\infty}
 \mu(h^+_0>a) e^{\frac{\lambda}{\nu} a} da \nonumber\\
  &\lesssim & ||h^+_0||_1 \left(1+\int_{\nu/\lambda}^{||h_0^+||_{\infty}} \frac{e^{\frac{\lambda}{\nu}a}}{a}
 da \right)\nonumber\\
  &\lesssim& ||h^+_0||_1 e^{\frac{\lambda}{\nu}||h^+_0||_{\infty}}\EEA
 so 
 \BEQ \bar{I}_{\infty}:=||h^+_0||_1 e^{\frac{\lambda}{\nu}||h^+_0||_{\infty}}  \label{eq:barIinfty}\EEQ
 is an upper bound for $\sup_{t\ge 0} ||\bar{h}_t||_1$ (see \cite{LauSou}, Proposition 2 (iii)); 
 at the same time, one gets
\BEQ ||\bar{h}_t||_{\infty} \le \frac{\nu}{\lambda} ||w_t||_{\infty}\lesssim \frac{\nu}{\lambda}
||w_0||_1 (\nu t)^{-d/2} \lesssim \bar{I}_{\infty} t^{-d/2}.  \EEQ

On the other hand, (\ref{eq:h-}) gives immediately if $h_0^-\in L^1$
\BEQ ||\underline{h}_t||_{\infty}\lesssim ||h_0^-||_1 t^{-d/2}.\EEQ

\medskip

Thus one has shown a global bound for the $L^1$-norm, and a time-decay in $O(t^{-d/2})$ for the sup-norm of solutions of (\ref{eq:KPZdet}) with arbitrary integrable
initial condition $h_0\in {\cal W}^{1,\infty}\cap L^1$.


\subsection{Time-decay of solutions of viscous Hamilton-Jacobi equations}


A second series of results is a particular case of the more general  {\em time-decay of the  gradients of solutions of viscous Hamilton-Jacobi
equations}, which can itself be seen as (1) an extension to non-linear equations of the standard parabolic
estimates; (2) or a multi-dimensional extension of the decay of  solutions of scalar 
conservation laws, see \cite{GilGueKer}, \cite{BenLau}, \cite{BenBenLau}, or \cite{BenKarLau}, section 3 for further
results concerning in particular  single-sided bounds on the Hessian.  Generally speaking, such results rely on
 convexity assumptions on $V$. Here we shall only state the following estimate, 
which is an extension of  \cite{BenLau}, Theorem 1. In Proposition \ref{prop:decay-gradient}, by exception, constants, explicit and
implicit (i.e. hidden by the sign $\lesssim$) are $\nu$-independent.

\begin{Proposition}[time-decay of the gradient] \label{prop:decay-gradient}
\begin{itemize}
\item[(i)] If $yV'(y)-V(y)\ge 0$ , then \BEQ ||\nabla h_t||_{\infty}\lesssim ||h_0||_{\infty} (\nu t)^{-\half}. \label{eq:decay-gradient1} \EEQ    
\item[(ii)] Under the stronger assumption 
\BEQ yV'(y)-V(y)\ge C\min(y^2,y^q),\qquad y\ge 0 \label{eq:stronger-assumption}\EEQ
for some constant $C>0$ and some exponent $q\in(1,2]$, one has 
\BEQ ||\nabla h_t||_{\infty}\lesssim  \left(\frac{||h_0||_{\infty}/\lambda}{t}\right)^{1/q}, \qquad
t\le \frac{||h_0||_{\infty}}{\lambda}  \label{eq:1.14}\EEQ
\BEQ  ||\nabla h_t||_{\infty}\lesssim  \left(\frac{||h_0||_{\infty}/\lambda}{t}\right)^{1/2}, \qquad
t\ge \frac{||h_0||_{\infty}}{\lambda}  \label{eq:1.15}\EEQ
\item[(iii)] Under the even stronger assumption 
\BEQ yV'(y)-V(y)\ge Cy^2 \label{eq:even-stronger-assumption}, \EEQ
 one has for all $t\ge 0$
 \BEQ |\nabla h_t(x)|\lesssim \left(\frac{|h_t(x)|/\lambda}{t}\right)^{1/2}, \qquad x\in\R^d.
 \label{eq:pointwise-decay-gradient} \EEQ
 Hence in particular
\BEQ  ||\nabla h_t||_{\infty}\lesssim  \left(\frac{||h_0||_{\infty}/\lambda}{t}\right)^{1/2}. \label{eq:decay-gradient2}\EEQ 
\end{itemize}
\end{Proposition}

Note that condition (i), $yV'(y)-V(y)\ge 0$ is a consequence of the convexity of $V$ (see  Assumption \ref{assumptions} (3)). On the other hand,  the hypothesis (\ref{eq:stronger-assumption}) in (ii) holds true for functions $V(y)$ that behave like $y^2$ for $y$ small, and like
$y^q$, $1<q\le 2$ for $y$ large; the stronger hypothesis (\ref{eq:even-stronger-assumption}) in (iii) for
functions that behave like $y^2$ both for $y$ small and $y$ large. Note that the decay in (\ref{eq:decay-gradient1}) is produced by the diffusion term $\nu\Del$ in the equation, so it might be called a generalized parabolic estimate; while (\ref{eq:1.14},\ref{eq:1.15})
or (\ref{eq:decay-gradient2}) are diffusion-independent effects of the non-linear term in the equation, and would
also hold true for viscosity solutions of the first-order Hamilton-Jacobi equation obtained by letting $\nu\to 0^+$.
\medskip

{\bf Proof.} We first rescale $h$ and $x$ by letting $x\to x'=\nu^{-\half} x$, $h\to u=\frac{\lambda}{\nu} h$ 
so that $\nabla_x h=\left(\frac{\lambda^2}{\nu}\right)^{-\half} \nabla_{x'}u$, and
$W(y)=\frac{\lambda^2}{\nu} V(\left(\frac{\lambda^2}{\nu}\right)^{-\half}y)$, so that the equation for $u$
\BEQ \partial_t u=\Del u +W(|\nabla u|) \EEQ
is independent of the parameters $\nu,\lambda$. 

\begin{itemize}
\item[(i)] 
Referring to the proof of Lemma 3 in \cite{GilGueKer}, from which \cite{GilGueKer},
Theorem 2 follows immediately, and letting directly $\eps=0$,  we see that $\frac{|\nabla u|^2}{\theta^2(u)}$ is a 
super-solution for the parabolic
operator
\BEQ \widetilde{\cal N}(w):=\Del w+b\cdot \nabla w+cw^2+ew-\partial_t w,\EEQ
where  $c:=2\theta(u)\theta''(u)$, $e=-2(W(\nabla u)-\nabla u\cdot W'(\nabla u))$
(note that for $V$ homogeneous, $V(|\nabla u|)=|\nabla u|^q$,  $\widetilde{\cal N}(\frac{|\nabla u|^2}{\theta^2(u)})={\cal N}(\frac{|\nabla u|^2}{\theta^2(u)})$ where $\cal N$ \cite{GilGueKer} has
instead of $ew$  a sum of two terms, $2(q-1)\theta^{q-1}(u)\theta'(u)w^{1+q/2}-2\frac{\theta'(u)}{\theta(u)}
H(\nabla u)w$, with $H(\xi)=W(\xi)-\xi\cdot\nabla W(\xi)+(q-1)|\xi|^q$, the first term $2(q-1)\theta^{q-1}(u)\theta'(u)
w^{1+q/2}$ compensating the last term $(q-1)|\xi|^q$ in $H$). Choose $\theta=\theta_1$ as in \cite{GilGueKer},
eq. (24), so that $c=-1$. Now, for $V$ isotropic, $e=-2\frac{\theta'(u)}{\theta(u)}(W(y)-yW'(y))$, $y=|\nabla u|$;
this is $\le 0$ under the assumptions (\ref{assumptions}). Hence $t^{-1}$ is a sub-solution of $\widetilde{\cal N}$, and the
comparison principle yields $|\nabla u_t(x)|^2\le \frac{\theta_1^2(u_0(x))}{t}$. Now $||\theta_1||_{\infty}\le 
2||u_0||_{\infty}$. Scaling back to the original variables $h,x$ yields  the first bound, $||\nabla h_t||_{\infty}\lesssim ||h_0||_{\infty} (\nu t)^{-\half}.$

\item[(ii)] The second bound is an 
extension of \cite{BenLau}. Exactly as in (i), one may assume that $u_0\le 0$. Up to an overall change of sign, $u\mapsto -u$, we are in the
conditions of \cite{BenLau}, Theorem 1, with $p=2$, except that $V$ is not necessarily a power function. Letting again $\eps=0$, the function $\Theta$ in eq. (20) p. 2005 is here $\Theta(r)=\Theta(y^2)=2y^2\frac{d}{dy^2}W(y)-W(y)=yW'(y)-W(y)$;
by assumption, 
\BEQ \Theta(r)\gtrsim \min(r,\left(\frac{\lambda^2}{\nu}\right)^{1-\frac{q}{2}} r^{q/2})=
r {\bf 1}_{r\le \frac{\lambda^2}{\nu}}+\left( \frac{\lambda^2}{\nu}\right)^{1-\frac{q}{2}} r^{\frac{q}{2}}
{\bf 1}_{r>\frac{\lambda^2}{\nu}}.\EEQ
 Eq. (14) p. 2003 implies
\BEQ {\cal L}w+Cv^{-2}\Theta(v^2 w)w\le 0\EEQ
with ${\cal L}=\partial_t-\Del$ up to some gradient term vanishing on functions $h$ which are independent of $x$ (see eq. (10) p. 2002 for a precise definition), $v:=\sqrt{u}$, $w=|\nabla v|^2$. Now  $v\le ||u_0||_{\infty}^{\half}$, hence
\BEA &&  v^{-2}\Theta(v^2 w)w\ge w^2 \qquad (w\lesssim \left( \frac{\nu||u_0||_{\infty}}{\lambda^2}\right)^{-1}),
\nonumber\\ && \qquad \qquad 
v^{-2}\Theta(v^2 w)w\ge  \left( \frac{\nu||u_0||_{\infty}}{\lambda^2}\right)^{\frac{q}{2}-1} w^{1+\frac{q}{2}} \qquad (w\gtrsim \left(\frac{\nu ||u_0||_{\infty}}{\lambda^2}\right)^{-1}) \nonumber\EEA
so ${\cal N}w\le 0$, where ${\cal N}$ is the parabolic differential operator
\BEQ {\cal N}:h\mapsto {\cal L}h+C N(h),\qquad N(h)=h^2 {\bf 1}_{h\le \left(\frac{\nu||u_0||_{\infty}}{\lambda^2}\right)^{-1}}+
\left(\frac{\nu||u_0||_{\infty}}{\lambda^2}\right)^{\frac{q}{2}-1} h^{1+\frac{q}{2}} {\bf 1}_{h> \left(
\frac{\nu||u_0||_{\infty}}{\lambda^2}\right)^{-1}}.\EEQ
Note that $N$ is an increasing function. The comparison principle thus implies that $w\le h$ if ${\cal N}h\ge 0$.
Such a function $h=h(t)$ is easily constructed by solving the ordinary differential equations $\partial_t h=-
\left(\frac{\nu||u_0||_{\infty}}{\lambda^2}\right)^{\frac{q}{2}-1} h^{1+\frac{q}{2}}$ for $h\le 
\left(\frac{\nu||u_0||_{\infty}}{\lambda^2}\right)^{-1}$, $\partial_t h=-h^2$ for $h\ge 
\left(\frac{\nu||u_0||_{\infty}}{\lambda^2}\right)^{-1}$, yielding up to unimportant constants $h(t)=
\left(\frac{\nu||u_0||_{\infty}}{\lambda^2}\right)^{\frac{2}{q}-1} t^{-\frac{2}{q}}$ $(t\le
\frac{\nu||u_0||_{\infty}}{\lambda^2})$,
$h(t)=\frac{1}{t}$ $(t\ge \frac{\nu||u_0||_{\infty}}{\lambda^2})$. This gives bounds for $||\nabla v_t||_{\infty}$ by taking the square-root,
and then bounds for $||\nabla u_t||_{\infty}$ by noting that $\nabla u=2u^{1/2}\nabla v$ and $u^{\half}_t(x)\le ||u_0||_{\infty}^{\half}$ (see \cite{BenKarLau}, proof of Proposition 3.1), namely,
\BEQ ||\nabla u_t||_{\infty}\lesssim \left(\frac{\nu}{\lambda^2}\right)^{-\half} \left( \frac{\nu||u_0||_{\infty}}{\lambda^2 t}\right)^{1/q}, \qquad t\le \frac{\nu ||u_0||_{\infty}}{\lambda^2};\EEQ
\BEQ ||\nabla u_t||_{\infty}\lesssim \left(\frac{||u_0||_{\infty}}{t}\right)^{\half},\qquad
 t\ge \frac{\nu ||u_0||_{\infty}}{\lambda^2}.\EEQ
 Hence (\ref{eq:1.14},\ref{eq:1.15}) by rescaling. 
 
 \item[(iii)] Under the stronger assumption (\ref{eq:even-stronger-assumption}), the previous computations
 yields $h(t)=\frac{1}{t}$ for all $t>0$. Hence $|\nabla u_t(x)|\lesssim \left( \frac{|u_t(x)|}{t}\right)^{\half}$. Eq. (\ref{eq:pointwise-decay-gradient}), and then (\ref{eq:decay-gradient2}),  follow by rescaling and using
 the a priori bound $||h_t||_{\infty}\le || h_0||_{\infty}$.

\end{itemize}
\hfill\eop


\subsection{Bounds through integral representation of mild solutions}


The third source of results is the integral form of the equation,
\BEQ h_t=e^{t\nu\Del}h_0+\lambda\int_0^t e^{(t-s)\nu\Del} V(\nabla h_s)ds, \label{eq:integral-form} \EEQ
the solutions of which, traditionally called {\em mild solutions}, are not necessary twice differentiable in space.
(\ref{eq:integral-form}) is used to prove local-in-time well-posedness of the equation, while the  a priori estimates
(\ref{eq:a-priori-estimates}) imply  global existence \cite{AmoBen}. We shall not come back to this; instead,  we give an application to the
proof of various bounds for the gradient and for higher derivatives of the solution. Generally speaking,
$\nu$-dependent constants
(throughout denoted by $C$ and possibly varying from line to line) come out of the computations everywhere. From \cite{AmoBen}, the solution
obtained by iterating (\ref{eq:integral-form}),
\BEQ h^{(0)}_t=h_0,\qquad h^{(k+1)}_t=e^{t\nu\Del}h^{(k)}_0+\lambda\int_0^t e^{(t-s)\nu\Del} V(\nabla h^{(k)}_s)ds
\qquad (k\ge 0) \label{eq:iterative} \EEQ
 in search for a fixed point is obtained as a converging series for
$t<T_1^*$, \BEQ T_1^*=C(\lambda||\nabla h_0||_{\infty})^{-2} \EEQ
 for some constant $C$, and shown to be uniformly smooth: namely, for every $k\ge 0$,
$||\nabla^k h_t||_{\infty}\le C_k||\nabla^k h_0||_{\infty},\ 0\le t\le T_1^*$ provided the initial solution has
bounded derivatives of  order $\le k$. By an appropriate choice of $C$ one may assume that $C_2=C_3=2$. This, in turn, shows, using the a priori bound, $||\nabla h_t||_{\infty}\le ||\nabla h_0||_{\infty}$,  that the solution at any later time also has bounded
derivatives of arbitrary order. We are interested here in quantitative bounds that can be shown to be close to
optimal in some case where explicit computations are possible (see next paragraph).

\bigskip

We shall give two different results. Recall $\bar{I}_{\infty}=||h^+_0||_1 e^{\frac{\lambda}{\nu} ||h^+_0||_{\infty}}$
(see (\ref{eq:barIinfty}). Our first result uses hypothesis (\ref{eq:even-stronger-assumption}), $yV'(y)-V(y)\ge Cy^2$
(see Proposition \ref{prop:decay-gradient} (iii)).

\begin{Theorem}[decay in $L^1$-norm of the gradient] \label{th:decay-nabla-L1}
Assume $h_0\in W^{1,\infty}\cap L^1$ and let $h_t$ be the solution of the KPZ equation (\ref{eq:KPZdet}), with
$V$ satisfying the hypothesis (\ref{eq:even-stronger-assumption}). Then 
\BEQ ||\nabla h_t||_1\lesssim \max(J_{\infty},J_{\infty}^{1+\frac{1}{d}}) (1+t)^{-\half},\qquad J_{\infty}=
\sup\left(1,\frac{||\nabla h_0||_1}{||h_0||_1},\lambda ||\nabla h_0||_{\infty} (1+O(\lambda ||h_0||_{\infty}))
\right) ||h_0||_1 e^{C\lambda ||h_0||_{\infty}}. \label{eq:decay-nabla-L1}\EEQ
\end{Theorem}

The time-decay in $O(t^{-\half})$ of $||\nabla h_t||_1$ is shown in \cite{BenKarLau} for $V(y)=y^2$, but with a constant $J_{\infty}$
which is roughly $e^{\bar{I}_{\infty}}$ and thus far from optimal (see p. 1290 and 1291). The emphasis there was on the
asymptotic convergence for $t\to \infty$ of the solution to a multiple of the heat-kernel (see Theorem 2.3 (a)),
 an interesting result in itself to which we do not come back here. 

\medskip

{\bf Proof.}

Our proof, based on intuition derived from the explicit computations of the next paragraph, shows that there are
different time regimes for $||\nabla h_t||_1$. Initially (i) the $L^1$-norm of the gradient may increase (as is the
case for the $L^1$-norm of the solution when the initial condition is positive); for later times (iii) it decreases like
the square-root of time. There
also appears a regime (ii) for intermediate times, during which the $L^1$-norm of the gradient is shown to be essentially
constant.

These three regimes come from the three essentially different bounds one has on $||\nabla h_t||_{\infty}$; namely,
(i) $||\nabla h_t||_{\infty}\le ||\nabla h_0||_{\infty}$ by the comparison principle; (ii) $||\nabla h_t||_{\infty}\lesssim \sqrt{\frac{||h_0||_{\infty}}{\lambda}} t^{-\half}$ by Proposition \ref{prop:decay-gradient} (iii); (iii) 
\BEQ ||\nabla h_t||_{\infty}\lesssim ||h_{t/2}||_{\infty} t^{-\half}\lesssim I_{\infty} t^{-(d+1)/2},
\label{eq:d+1/2} \EEQ
 as follows from a combination of Proposition \ref{prop:decay-gradient} (i) and of the parabolic estimates
developed in the lines following eq. (\ref{eq:barIinfty}), where 
\BEQ I_{\infty}=
||h_0||_1 e^{C\lambda||h_0||_{\infty}} \label{eq:Iinfty} \EEQ
is an upper bound for $\sup_{t\ge 0} ||h_t||_1$ (in order to get not too complicated formulas, we avoid the unpleasant task of optimizing the constants, and choose $C$ large enough).

\begin{itemize}
\item[(i)] For $t$ small one uses the trivial bound (i), $||\nabla h_t||_{\infty}\le ||\nabla h_0||_{\infty}$,
and applies the iterative scheme (\ref{eq:iterative}) in uniform time slices $[T_0^*,T_1^*]=[0,T_1^*],[T_1^*,T_2^*]=
[T_1^*,2T_1^*],\ldots,[T_{n_0-1}^*,T_{n_0}^*]=[(n_0-1)T_1^*,n_0 T_1^*]$ where $n_0\approx \lambda ||h_0||_{\infty}$,
so that $||\nabla h_0||_{\infty}\approx \sqrt{\frac{||u_0||_{\infty}}{\lambda}} (T_{n_0}^*)^{-\half}$. At some time
comparable with $T_{n_0}^*$,  the bound (ii) on $||\nabla h_t||_{\infty}$ becomes better.  We let
$M_n^{(k)}:=\sup_{[T_n^*,T_{n+1}^*]} ||\nabla u_t^{(k)}||_1$ and $M_n:=\sup_{[T_n^*,T_{n+1}^*]} ||\nabla u_t||_1=
\lim_{k\to\infty} M_n^{(k)}$. By (\ref{eq:iterative}) and the parabolic estimates recalled in (\ref{eq:parabolic-estimates}),
\BEQ M_n^{(k+1)}\lesssim ||\nabla h_{T_n^*}||_1+\lambda \sup_{t\in[T_n^*,T_{n+1}^*]} \int_{T_n^*}^t
(t-s)^{-\half} ||(\nabla h^{(k-1)}(s))^2||_1 ds, \label{eq:M1}\EEQ
together with the interpolation inequality, 
\BEQ ||(\nabla h^{(k-1)}(s))^2||_1=||\nabla h^{(k-1)}(s)||_2^2\le
||\nabla h^{(k-1)}(s)||_1 ||\nabla h^{(k-1)}(s)||_{\infty}, \EEQ one obtains
\BEQ M_n^{(k+1)}\lesssim ||\nabla h_{T_n^*}||_1+C^{-1} M_n^{(k)} \EEQ
where $C^{-1}$ is proportional to the inverse of the constant $C$ in the definition of $T_1^*$. For $C$ large
enough, this yields $\sup_k M_n^{(k)}\le 2M_{n-1}$ and $M_n\le 2M_{n-1}$. Thus
\BEQ M_{n_0}\le ||\nabla h_0||_1 e^{C\lambda ||h_0||_{\infty}}\lesssim \frac{||\nabla h_0||_1}{||h_0||_1} I_{\infty}
\label{eq:2.32} \EEQ
with an appropriate definition of the constant in (\ref{eq:Iinfty}).

\item[(ii)] For $n\ge n_0$ one defines inductively $T_n^*$ by $T_{n+1}^*-T_n^*=C(\lambda||\nabla h_{T_n^*}||_{\infty})^{-2}$. Note that $T_{n_0}^*\approx \frac{||h_0||_{\infty}}{||\nabla h_0||_{\infty}^2}$; by the second
estimate (ii) on $||\nabla h_t||_{\infty}$, $T_{n+1}^*-T_n^*\gtrsim \frac{T_n^*}{\lambda ||h_0||_{\infty}}$,
so 
\BEQ T_n^*\ge \frac{||h_0||_{\infty}}{||\nabla h_0||_{\infty}^2}\left( 1+\frac{C}{\lambda||h_0||_{\infty}}\right)^{n-n_0},\qquad n\ge n_0.  \label{eq:1.32} \EEQ

Instead of  the bound $||\nabla e^{(t-T_n^*)\nu\Del} h_{T_n^*}||_1\le||\nabla h_{T_n^*}||_1$ used in (\ref{eq:M1}), it is more clever for $n$ and $t-T_n^*$ large enough to use the parabolic
estimate $||\nabla e^{(t-T_n^*)\nu\Del}h_{T_n^*}||_1\lesssim (t-T_n^*)^{-\half} I_{\infty}$ if the latter
expression is $\le ||\nabla h_{T_n^*}||_1$. Thus one gets the improved estimate
\BEQ M_n^{(k+1)}\lesssim\sup_{t\in[T_n^*,T_{n+1}^*]} \left( \inf(||\nabla h_{T_n^*}||_1,I_{\infty}
(t-T_n^*)^{-\half}) +\lambda||\nabla h_{T_n^*}||_{\infty} M_n^{(k)} (t-T_n^*)^{\half}\right). \label{eq:M2}\EEQ

If $||\nabla h_{T_n^*}||_1 \gtrsim I_{\infty} (T_{n+1}^*-T_n^*)^{-\half}\approx \lambda I_{\infty}
||\nabla h_{T_n^*}||_{\infty}$, the improved estimate (\ref{eq:M2})
is better than (\ref{eq:M1}) and yields
\BEA && M_n^{(k+1)}\lesssim \sup\left( ||\nabla h_{T_n^*}||_1+\lambda ||\nabla h_{T_n^*}||_{\infty} 
M_n^{(k)} \frac{||h_{T_n^*}||_1}{||\nabla h_{T_n^*}||_1}, \right.\nonumber\\
&&\qquad \qquad \left. \sup_{T_{n+1}^*-T_n^*\ge t-T_n^*\ge (||h_{T_n^*}||_1/||\nabla h_{T_n^*}||_1)^2}
\left( I_{\infty} (t-T_n^*)^{-\half}+\lambda ||\nabla h_{T_n^*}||_{\infty} M_n^{(k)} (t-T_n^*)^{\half} \right)\right)
\nonumber\\ \EEA
The function $x\mapsto \frac{a}{x}+bx$, here $x=\sqrt{t-T_n^*}$, is bounded on any interval of $\R_+$ by the max of its values at the two ends of the
interval. Hence
\BEQ M_n^{(k+1)}\lesssim \sup\left(  ||\nabla h_{T_n^*}||_1+\lambda ||\nabla h_{T_n^*}||_{\infty} 
M_n^{(k)} \frac{I_{\infty}}{||\nabla h_{T_n^*}||_1}, I_{\infty} (T_{n+1}^*-T_n^*)^{-\half}+C^{-1}M_n^{(k)}
\right)\EEQ
with $C^{-1}<1$. Iterating these affine inequalities yields either 
\BEQ M_n\lesssim I_{\infty}  (T_{n+1}^*-T_n^*)^{-\half} \approx \lambda I_{\infty} ||\nabla h_{T_n^*}||_{\infty}
\le I_{\infty} ||\nabla h_0||_{\infty}; \EEQ
or, assuming on the contrary that $||\nabla h_{T_n^*}||_1\gtrsim I_{\infty}(T_{n+1}^*-T_n^*)^{-\half}$,
\BEQ M_n\lesssim  \frac{||\nabla u_{T_n^*}||_1}{1-\lambda I_{\infty} \frac{||\nabla u_{T_n^*}||_{\infty}}{||\nabla
u_{T_n^*}||_1}}.\EEQ
Since the sequence $n\mapsto \lambda I_{\infty}||\nabla h_{T_n^*}||_{\infty}$ is exponentially decreasing
(as follows from the bound (ii) and the fact that the sequence $(T_n^*)$ is exponentially increasing, see
(\ref{eq:1.32})), the
recursive sequence
\BEQ x_{n+1}=\frac{x_n}{1-\lambda I_{\infty}\frac{ ||\nabla h_{T_n^*}||_{\infty}}{x_n}}\approx x_n+
\lambda I_{\infty} ||\nabla h_{T_n^*}||_{\infty}\EEQ
starting from $x_{n_1}\approx I_{\infty}(T_{n+1}^*-T_n^*)^{-\half}\lesssim \lambda||\nabla h_0||_{\infty} I_{\infty}$,
converges to 
\BEQ x_{\infty} \lesssim \lambda||\nabla h_0||_{\infty}I_{\infty} \left(1+\frac{1}{1-(1+\frac{C}{\lambda||h_0||_{\infty}})^{-1/2}}\right)\lesssim \lambda ||\nabla h_0||_{\infty} I_{\infty} (1+O(\lambda ||h_0||_{\infty})).\EEQ

This gives a global bound for $M_n,n\ge 0$,
\BEQ \sup_{n\ge 0}M_n\le \tilde{I}_{\infty}:=\sup\left(\frac{||\nabla h_0||_1}{||h_0||_1},
\lambda ||\nabla h_0||_{\infty} (1+O(\lambda ||h_0||_{\infty})) \right) 
I_{\infty},\EEQ
but no time decay yet in general.

\item[(iii)] For $t\gtrsim I_{\infty}^{2/d}$ we use the estimate (iii), $||\nabla h_t||_{\infty}\lesssim I_{\infty}
t^{-(d+1)/2}$ and prove the time decay in $O(t^{-\half})$. Let $\tilde{M}_n:=\sup_{t\in[2^{n-1},2^n]}
||\nabla h_t||_1$ for $n\ge n_2:=1+\log_2 I_{\infty}^{2/d}$. For all $t\in[2^n,2^{n+1})$,
\BEA ||\nabla h_t||_1 &\lesssim & 2^{-n/2} ||h_{t-2^{n-1}}||_1+\lambda \int_{t-2^{n-1}}^{t} ds (t-s)^{-\half}
||\nabla h_s||_1 ||\nabla h_s||_{\infty} \nonumber\\
&\lesssim& 2^{-n/2}I_{\infty}+\lambda 2^{n/2} (\tilde{M}_n+\tilde{M}_{n+1})I_{\infty} (2^{-n})^{(d+1)/2},\EEA
hence
\BEQ \tilde{M}_{n+1}\lesssim 2^{-n/2}I_{\infty}+\lambda I_{\infty} 2^{-nd/2} \tilde{M}_n+\lambda I_{\infty}
2^{-nd/2}\tilde{M}_{n+1}.\EEQ
For $n\ge n_2$ one has by definition $I_{\infty}2^{-nd/2}\le 1$, so
\BEQ \tilde{M}_{n+1}\lesssim (1+O(\lambda))(2^{-n/2}I_{\infty}+\lambda\tilde{M}_n),\EEQ
while $\tilde{M}_{n_2}\lesssim \tilde{I}_{\infty}$ by (ii),
implying by a straightforward induction 
\BEQ \tilde{M}_n\lesssim 2^{-n/2} I_{\infty}+\tilde{I}_{\infty} \lesssim  2^{-n/2}(I_{\infty}+ I_{\infty}^{1/d}\tilde{I}_{\infty}) \EEQ
and finally
\BEQ \sup_{t\ge I_{\infty}^{2/d}} \sqrt{t} ||\nabla h_t||_1\lesssim I_{\infty}+ I_{\infty}^{1/d} \tilde{I}_{\infty}.\EEQ
Finally,
\BEQ \sup_{t\le I_{\infty}^{2/d}} \sqrt{t} ||\nabla h_t||_1\lesssim I_{\infty}^{1/d} \tilde{I}_{\infty}.\EEQ
Hence the result.\hfill\eop
\end{itemize}

\smallskip

{\bf Remark.} If $V$ does not satisfy (\ref{eq:even-stronger-assumption}), then the beginning of the proof is modified as follows:
substituting to (ii) the bound $||\nabla h_t||_{\infty}\lesssim ||h_0||_{\infty} t^{-\half}$, see Proposition
\ref{prop:decay-gradient} (i), leads to 
$n_0$ defined such as to satisfy $||\nabla h_0||_{\infty}\approx ||h_0||_{\infty} (T_{n_0}^*)^{-\half}$, namely,
$n_0\approx (\lambda ||h_0||_{\infty})^2$, and (compare with (\ref{eq:2.32})) $M_{n_0}\lesssim ||\nabla h_0||_1
e^{C(\lambda ||h_0||_{\infty})^2}$. A bound comparable to (\ref{eq:decay-nabla-L1}) probably holds with
the quadratic exponential $e^{C(\lambda ||h_0||_{\infty})^2}$ substituting $e^{C\lambda ||h_0||_{\infty}}$, which
is clearly not optimal for the quadratic KPZ equation (see next paragraph). 


\medskip
Our second result is valid under our general assumptions on $V$ stated in section 1.

\begin{Theorem}[bounds on higher derivatives] \label{th:bound-higher-derivatives}

Let $h$ be the solution of eq. (\ref{eq:KPZdet}) with initial condition $h_0\in {\cal W}^{3,\infty}$. Then:

\BEQ ||\nabla^2 h_t||_{\infty}\lesssim  P_1(||h_0||_{\infty}, ||\nabla h_0||_{\infty}, ||\nabla^2 h_0||_{\infty})
\frac{\ln (1+t)}{t}   \label{eq:nabla2}\EEQ
\BEQ ||\nabla^3 h_t||_{\infty}\lesssim P_2(||h_0||_{\infty}, ||\nabla h_0||_{\infty}, ||\nabla^2 h_0||_{\infty},
||\nabla^3 h_0||_{\infty}) \left(\frac{\ln^2 (1+t)}{t^{3/2}}\right) \label{eq:nabla3} \EEQ
where $P_1,P_2$ are polynomials.
\end{Theorem}

{\bf Proof.} 

We already know that $\sup_{[0,T_1^*]} ||\nabla^k h_t||_{\infty}\le 2||\nabla^k h_0||_{\infty}$, $k=2,3$ 
for $T_1^*\approx (\lambda||\nabla h_0||_{\infty})^{-2}.$
For $t\ge T_1^*$, $\nabla^2 h_t$ is the solution of an integral equation,
\BEQ \nabla^2 h_t=e^{t\nu\Del}\nabla^2 h_0+\lambda \left[ \int_{(1-\eps)t}^t ds (\nabla e^{(t-s)\nu\Del}) \nabla
(V(\nabla h_s)) + \int_0^{(1-\eps) t} ds (\nabla^2 e^{(t-s)\nu\Del}) V(\nabla h_s) \right].\EEQ
The idea is to commute  the gradient with the  heat operator $e^{(t-s)\nu\Del}$ in order to make the most
of  parabolic estimates; the implied
decay may be put to good use only for $t-s$ large enough, and we shall choose the parameter $\eps$ accordingly.
First, $||e^{t\nu\Del}\nabla^2 h_0||_{\infty}=||\nabla^2 e^{t\nu\Del}h_0||_{\infty}\lesssim ||h_0||_{\infty} t^{-1}$. Then, using $\nabla (V(\nabla h_s))=V'(\nabla h_s)\cdot \nabla^2 h_s$ and Proposition \ref{prop:decay-gradient} (i),
together with the inequality $V'(y)\lesssim y$, consequence of Assumption (\ref{assumptions}) (3), (4) (namely,
$(2y)^2\ge V(2y)\ge V(y)+yV'(y)$)
\BEQ   \lambda\left|\int_{(1-\eps)t}^t ds (\nabla e^{(t-s)\nu\Del}) \nabla
(V(\nabla h_s)) \right|  
 \lesssim \lambda \int_{(1-\eps)t}^t \frac{ds}{\sqrt{t-s}} ||\nabla h_s||_{\infty} \sup_{[(1-\eps)t,t]} ||\nabla^2
h_s||_{\infty} 
\lesssim\lambda  ||h_0||_{\infty} \sqrt{\eps} \sup_{[(1-\eps)t,t]} ||\nabla^2 h_s||_{\infty} \nonumber\\
\EEQ

provided $(1-\eps)t\gtrsim t$. To get a useful inequality we choose $\eps$ so that 
$\lambda||h_0||_{\infty}\sqrt{\eps}\le \frac{1}{4}$, namely, $\eps\approx \min(\frac{1}{4},\frac{1}{(\lambda^
||h_0||_{\infty})^2}) $. Finally (if $t\ge 1$) we split the second integral, $\int_0^{(1-\eps)t}$, into several pieces:

\BEQ \lambda \int_{t/2}^{(1-\eps)t} ds (\nabla^2 e^{(t-s)\nu\Del}) V(\nabla h_s) \lesssim
\lambda \int_{t/2}^{(1-\eps)t} \frac{ds}{t-s} \left( \frac{||h_0||_{\infty}}{\sqrt{t}} \right)^2 \lesssim \lambda
\ln(\eps^{-1}) ||h_0||_{\infty}^2 t^{-1};\EEQ

\BEQ \lambda \int_{1}^{t/2} ds (\nabla^2 e^{(t-s)\nu\Del}) V(\nabla h_s) \lesssim
\lambda t^{-1} \int_1^{t/2} \frac{ds}{s} ||h_0||_{\infty}^2 \lesssim \lambda\frac{\ln t}{t} ||h_0||^2_{\infty};
\label{eq:1.52}\EEQ

\BEQ \lambda \int_0^1 ds (\nabla^2 e^{(t-s)\nu\Del}) V(\nabla h_s) \lesssim
\lambda t^{-1}  \int_0^1 ds ||\nabla h_0||^2_{\infty}= \lambda ||\nabla h_0||^2_{\infty} t^{-1}. \label{eq:1.53}\EEQ

If $t<1$, we merge (\ref{eq:1.52},\ref{eq:1.53}) into
\BEQ \lambda \int_0^{t/2} ds (\nabla^2 e^{(t-s)\nu\Del}) V(\nabla h_s) \lesssim \lambda t^{-1}
\int_0^{t/2} ds ||\nabla h_0||_{\infty}^2=\frac{\lambda}{2} ||\nabla h_0||_{\infty}^2.\EEQ

We finish as in the proof of Theorem \ref{th:decay-nabla-L1} (iii), namely,  letting 
$M_0:=\sup_{[0,T_1^*]} ||\nabla^2 h_t||_{\infty}$ and $M_n:=\sup_{[2^{n-1} T_1^*,2^{n}T_1^*]} ||\nabla^2 h_t||_{\infty}$ ($n\ge 1$),
one has $M_0\le 2||\nabla^2 h_0||_{\infty}$ and 
\BEQ M_{n+1}\le C^{-1} M_n+Q_1(||h_0||_{\infty},||\nabla h_0||_{\infty})2^{-n} +\lambda ||\nabla h_0||_{\infty}^2 \EEQ
for $t<1$,
\BEQ M_{n+1}\le C^{-1} M_n+Q_2(||h_0||_{\infty}, ||\nabla h_0||_{\infty})2^{-n}+Q_3(||h_0||_{\infty},||\nabla
h_0||_{\infty})n2^{-n} \EEQ
for $t\ge 1$. 
 Hence (\ref{eq:nabla2}).

\medskip

This result is used as input to get similar a bound for $||\nabla^3 h_t||_{\infty}$. This time we must move around
three gradients in the best way; this gives three integrals rewritten as 
$$ \int_{(1-\eps)t}^{t} ds
(\nabla e^{(t-s)\nu\Del})\nabla^2 (V(\nabla h_s)), \quad  \int_{t/2}^{(1-\eps)t} ds \nabla^2 e^{(t-s)\Del} 
\nabla(V(\nabla h(s))), \quad \int_0^{t/2} ds  \nabla^3 e^{(t-s)\nu\Del} V(\nabla h_s).$$  
One has $\nabla^2 (V(\nabla h_s))=V'(\nabla h_s)\cdot\nabla^3 h_s+V''(\nabla h_s)(\nabla^2 h_s)^2$, yielding the same constraints on $\eps$, plus a supplementary
quadratic term in $\nabla^2 h_s$. The other terms are computed as before. Details are left to the reader.

\hfill\eop


\subsection{An explicit example: decay of a 'bump' for the quadratic KPZ equation}


We consider here the time-decay (pointwise and with respect to various norms) of the solution of the
quadratic, homogeneous KPZ equation,
\BEQ \partial_t h=\Del h+|\nabla h|^2 \EEQ
with initial "bump" condition $h_0(x)=A {\bf 1}_{|x|\le L}$, where $A,L>0$. The coefficient $\lambda$ in front
of the nonlinearity has been disposed of by a simple rescaling. Note that, if $A=||h_0||_{\infty}\lesssim 1$ (i.e.
for a small initial condition), then the decay of the solution and of its derivatives in $L^p$-norms, $1\le p\le \infty$
follow the parabolic estimates as for the solutions of the  linear heat equation; thus we may assume that $A\gg 1$.
We want to compare the decays obtained by explicit computation to those obtained in much greater generality
in the previous paragraphs.

\medskip
 Through the exponential transformation, $w=\exp h$,
the equation becomes simply the heat equation,
\BEQ \partial_t w=\Del w,\qquad w_0=e^{h_0} \EEQ
so that
\BEA  h_t(x) & = & \ln\left( (2\pi t)^{-d/2} \int dy\,  e^{-(x-y)^2/2t} w_0(y)\right)  \nonumber\\
&=& \ln\left (1+ (e^A-1)\,  (2\pi t)^{-d/2}\int_{|y|\le L} dy\,  e^{-(x-y)^2/2t} \right) \nonumber\\
&\approx& \ln\left (1+ e^A t^{-d/2}\int_{|y|\le L} dy\,  e^{-(x-y)^2/2t} \right).
\EEA
Though the initial data is not in ${\cal W}^{1,\infty}$, this defines a solution. We are interested in its 
behaviour for $t\gtrsim L^2$, corresponding to the approximate amount of  time necessary for the solution to
smoothen up. Then
\BEQ h_t(x)\approx \ln\left(1+e^A \left(\frac{L^2}{t}\right)^{d/2} e^{-x^2/2t}\right).\EEQ There are two regimes:

\begin{itemize}
\item[(i)] (initial regime) Assume $L^2\lesssim t\lesssim L^2 e^{\frac{2}{d}A}$. Define $x_{max}(t)\in\R_+$ as the
solution of the equation $e^A \left(\frac{L^2}{t}\right)^{d/2} e^{-x^2_{max}(t)/2t}=1$; explicitly,
$x_{max}(t)=\sqrt{2t(A-\frac{d}{2}\log(\frac{t}{L^2}))}$. If $|x|\gtrsim x_{max}(t)$ then $h_t(x)\approx e^A
\left(\frac{L^2}{t}\right)^{d/2} e^{-x^2/2t}\lesssim 1$. On the other hand, if $|x|\lesssim x_{max}(t)$, then
$h_t(x)$ is still large, $h_t(x)\approx A-\frac{d}{2}\log\left(\frac{t}{L^2}\right)-\frac{x^2}{2t}$. In
particular, 
\BEQ ||h_t||_{\infty}\approx A-\frac{d}{2} \log\frac{t}{L^2} \label{eq:htinfty}\EEQ
 and 
\BEA  ||h_t||_{1} &\approx&  ||h_t||_{\infty}
\cdot \Vol(B(0,x_{max}(t)))+e^A L^d \int_{|x|>x_{max}(t)} t^{-d/2} e^{-|x|^2/2t}\, dx
\nonumber\\
&\lesssim & t^{d/2} \left(A-\frac{d}{2} \log\frac{t}{L^2}\right)^{\frac{d+1}{2}} +
e^A L^d. \label{eq:htL1} \EEA
Both quantities $t\mapsto t^{d/2}\left(A-\frac{d}{2}\log \frac{t}{L^2} \right)$ and
$t\mapsto t^{d/2}\left(A-\frac{d}{2}\log \frac{t}{L^2} \right)^{\frac{d+1}{2}}$ are
easily checked to be maximal for $A-\frac{d}{2}\log \frac{t}{L^2}\approx 1$, yielding
\BEQ ||h_t||_{\infty}\lesssim e^A L^d t^{-d/2}\approx  ||h_0||_1 \frac{e^{||h_0||_{\infty}}}{||h_0||_{\infty}} t^{-d/2} \EEQ
and
\BEQ ||h_t||_1\lesssim e^A L^d\approx  ||h_0||_1 \frac{e^{||h_0||_{\infty}}}{||h_0||_{\infty}}.\EEQ

\item[(ii)]  (final regime) Assume $t\gtrsim L^2 e^{\frac{2}{d} A}$. Then $h_t(x)\approx e^A \left(\frac{L^2}{t}
\right)^{d/2} e^{-x^2/2t}\lesssim 1$; in other words, the bump has essentially disappeared. Furthermore,
\BEQ ||h_t||_{\infty}\approx e^A L^d t^{-d/2}, \qquad ||h_t||_1\approx e^A L^d\EEQ
saturating the bounds found in case (i).

\end{itemize}

In both cases, bounds for $||\nabla^q h_t||_{\infty}$ or $||\nabla^q h_t||_1$ are obtained by dividing by $t^{q/2}$. Details are left to the reader.

The above computations make it clear that the general bounds for $||\nabla h_t||_{\infty}$, see eq. (\ref{eq:d+1/2}), 
$||\nabla h_t||_1$ and $||\nabla^q h_t||_{\infty}$, $q=2,3$ obtained in Theorems \ref{th:decay-nabla-L1} and \ref{th:bound-higher-derivatives}, are  essentially optimal.


\section{Bound for the homogeneous equation: the case of unbounded initial conditions}


We now want to prove existence of, and  bound,  solutions of the homogeneous KPZ equation,
\BEQ \partial_t h=\Del h+ \lambda V(\nabla h) \label{eq:2.1}\EEQ
with {\em unbounded} initial condition $h_0$ (for simplicity we fix $\nu=1$ from
now on). We would typically like to consider a random initial condition
which is a smoothened white noise (see Appendix A). This raises various problems. First (1), one would like to identify
a functional space preserved by the linear heat equation, for which generalized parabolic estimates hold. Second (2),
one would like to extend the comparison principle to such a functional space, in such a way as to prove existence
of and bound the solution. Finally (3), one would like to identify the solution as the limit of solutions of 
(\ref{eq:2.1}) associated to a sequence of compactly supported (hence bounded) initial conditions converging to the original initial condition, so as to extend to the limit regularity results and estimates obtained in the previous
section. 

\noindent We provide in this section  answers to   questions (1), (2), (3).
We first refer the reader to the Introduction for a short review. \S 3.1 is
devoted to a detailed study of the space ${\cal H}^0$ of functions with locally
bounded averages, and more generally of a family  of spaces ${\cal H}^0_{\alpha}$, 
including ${\cal H}^0\equiv {\cal H}^0_0$, which  enjoy the same type of properties.
The comparison principle  (Theorem 1 in the Introduction) is proved in \S 3.2.
Estimates for the solutions of (\ref{eq:2.1}) and their gradients are proved
respectively in \S 3.3 and \S 3.4.


\subsection{The functional spaces ${\cal H}^0_{\alpha}$}


For $f\in L^1_{loc}(\R^d)$ and $\alpha\ge 0$, $x\in\R^d$, one may define
\BEQ f^*_{\alpha}(x)=\sup_{\tau> 0}\,  (1+\tau)^{\alpha} e^{\tau\Del}|f|(x)\in[0,+\infty]\EEQ
and in particular
\BEQ f^*(x):=f^*_0(x)=\sup_{\tau>0} e^{\tau\Del}|f|(x)\in[0,+\infty]. \EEQ
Note that $f^*\le f_{\alpha}^*\le f_{\beta}^*$ if $\alpha\le\beta$.
If $f$ is bounded, then $f^*(x)\le ||f||_{\infty}$. On the other hand, the kind of random initial conditions
we are interested in (see Appendix A) are a.s. unbounded, but satisfy a.s. $f^*(x)<\infty$ for every $x$
(see Lemma \ref{lem:LD-main} and discussion thereafter); compare with the standard parabolic estimates, $f^*_{d/2}(x)\lesssim ||f||_1$ for $f\in L^1(\R^d)$. 
Note that, if $\alpha=0$, the obvious {\em pointwise} estimates, making part
of what we call {\em pointwise parabolic estimates} in Lemma \ref{lem:pointwise-parabolic-estimates}, 
\BEQ |e^{t\Del}f(x)|\le f^*(x),\qquad t\ge 0\EEQ
and, better still,
\BEQ (e^{t\Del}f)^*(x)\le f^*(x),\qquad t\ge 0\EEQ 
generalizing to 
\BEQ (e^{t\Del}f)^*_{\alpha}(x)\le f^*_{\alpha}(x),\qquad t\ge 0,\EEQ
are improvements on the {\em global} estimate $||e^{t\Del} f||_{\infty}\le ||f||_{\infty}$ which is useless
for unbounded functions. If $f^*_{\alpha}(x)<\infty$ for some $\alpha>0$, then 
\BEQ |e^{t\Del}f(x)|\le (1+t)^{-\alpha} f^*_{\alpha}(x)\EEQ
 decays polynomially in time.

\begin{Definition}[${\cal H}^0$-spaces]
Let, for $0\le\alpha\le d/2$,
\BEQ {\cal H}^0_{\alpha}:=\{ f\in L^{\infty}_{loc}(\R^d) \ |\ \forall x\in\R^d, f^*_{\alpha}(x)<\infty\}\EEQ
and
\BEQ {\cal H}^0:={\cal H}^0_0=\{ f\in L^{\infty}_{loc}(\R^d) \ |\ \forall x\in\R^d, f^*(x)<\infty\}.\EEQ
\end{Definition}

For every $\alpha\le d/2$, ${\cal H}^0_{\alpha}\supset L^1(\R^d)$ (actually, it is easy to prove that ${\cal H}^0_{\alpha}=
\{0\}$ for $\alpha>d/2$). 

It is easy to see that $f\in {\cal H}^0_{\alpha}$ provided there {\em exists}
some $x\in\R^d$ such that $f^*_{\alpha}(x)<\infty$. However the various "norms"
$f\mapsto f^*_{\alpha}(x)$, $x\in\R^d$, are not comparable. In this sense $f^*_{\alpha}(x)$
should be understood as a {\em local, $x$-centered measure} of the size of $f$.

\medskip

Another closely related definition is by averaging: if $f\in L^{\infty}_{loc}(\R^d)$, and $r\ge 0$, $x\in\R^d$, one may define
\BEQ f^{\sharp}_{\alpha}(x)=\sup_{\rho> 0}\,  (1+\rho^2)^{\alpha} 
\frac{\int_{B(x,\rho)} |f(y)|dy}{\Vol(B(x,\rho))}\EEQ 
(called  {\em maximal function} in real analysis for $\alpha=0$, see the classical book by  Stein \cite{Ste})
where $B(x,\rho)=\{y\in\R^d; |y-x|<\rho\}$ is the Euclidean ball and $\Vol(B(x,\rho))$ its volume.  Here also,
$f^{\#}_0\le f^{\#}_{\alpha}\le f^{\#}_{\beta}$ if $0\le\alpha\le\beta$. It is convenient
to denote averages by barred integrals, so that, by definition, $\fint_{\Omega} f=\frac{\int_{\Omega} f}{\Vol(\Omega)}$.
A simple result is the following:

\begin{Lemma} \label{eq:f*ftilde*}
There exists constants $c,C>0$ such that, for every  $f\in C(\R^d)$, $cf^*_{\alpha}(x)\le f^{\sharp}_{\alpha}(x)\le Cf^*_{\alpha}(x).$
\end{Lemma}

{\bf Proof.} We must prove two inequalities. First,
\BEA && \int dy \frac{e^{-(x-y)^2/2t}}{(2\pi t)^{d/2}} |f(y)| = \int dr \frac{e^{-r^2/2t}}{(2\pi t)^{d/2}}
\int_{\partial B(x,r)} dy |f(y)|  \nonumber\\
&& =\int dr \frac{r}{t} \frac{e^{-r^2/2t}}{(2\pi t)^{d/2}} \int_{B(x,r)} dy |f(y)| \nonumber\\
&&\le f^{\sharp}_0(x) \int dr \frac{r}{t} \frac{e^{-r^2/2t}}{(2\pi t)^{d/2}} \int_{B(x,r)} dy \nonumber\\
&&=f^{\sharp}_0(x) \int dy \frac{e^{-|y|^2/2t}}{(2\pi t)^{d/2}}=f^{\sharp}_0(x). \label{eq:H0H0}\EEA
Thus $f^*_0\le f^{\sharp}_0$. If $\alpha>0$ and $t\gtrsim 1$ then
\BEQ (1+t)^{\alpha} \int dy \frac{e^{-(x-y)^2/2t}}{(2\pi t)^{d/2}} |f(y)| \lesssim
 f^{\sharp}_{\alpha} (1+t)^{-1-\frac{d}{2}+\alpha} \int dr (1+r)^{1+d-2\alpha} e^{-r^2/2t} \lesssim f^{\sharp}_{\alpha}(x),\EEQ
 so $f^*_{\alpha}\lesssim f^{\#}_{\alpha}$.
 
Conversely, 
\BEQ e^{r^2\Del}|f|(x)=\int \frac{e^{-\half (|x-y|/r)^2}}{(2\pi)^{d/2}r^d} |f(y)|dy\ge C\fint_{B(x,r)} |f|.\EEQ
\hfill\eop

In particular, an equivalent definition for ${\cal H}^0_{\alpha}$ is:
\BEQ {\cal H}^0_{\alpha}=\{f\in L^{\infty}_{loc}(\R^d)\ |\ \forall x\in\R^d, f^{\#}_{\alpha}(x)<\infty\}.\EEQ

Note also that, if $f$ is lower semicontinuous (in particular, if $f$ is continuous),  $f^{\sharp}_{\alpha}(x)\ge \lim_{r\to 0} \fint_{B(x,r)} |f|\ge
|f(x)|$, and similarly
$f^*_{\alpha}(x)\ge |f(x)|$.
\medskip

{\bf Example.} The equivalence of the "pointwise quasi-norms" $f^*(x)$, $f^{\sharp}(x)$ makes it easy
to construct unbounded functions $f\in {\cal H}^0$. The idea is to modify a bounded function
on  regions with small relative volume.  Define for instance $f$ to be identically equal to $c_1>0$ outside the union of annuli
$\cup_{k\ge 0} B_k$, where $B_k:=B(0,2^k+2^{k\gamma})\setminus B(0,2^{k})$ for some $\gamma\in(-\infty,1)$, and
$f\big|_{B_k}:=c_2(2^k)^{d(1-\gamma')}$ with $c_2>0$ and $\gamma'\ge\gamma$. Then for $k$ large, $\fint_{B(x,2^k)} |f| \approx c_1$ (corresponding to "rare enough" fluctuations) if $\gamma'>\gamma$, and $\fint_{B(x,2^k)} |f| \approx c_1+c_2$
if $\gamma'=\gamma$ (corresponding to a border case where fluctuations are as important as the
bulk behaviour). Hence $f^*(x)<\infty$. On the other hand, $f^*(x)=\infty$ if $\gamma'<\gamma$. One may also allow arbitrarily large fluctuations by letting $\gamma'_k=\gamma_k$ be a sequence which is unbounded below.
Typical realizations of regularized white noise  are more complicated, but large
fluctuations do not contribute to the average on large balls (see section 6 for a more precise picture).
\medskip

\begin{Lemma} Let $\alpha\in[0,d/2]$ and $f\in {\cal H}^0_{\alpha}\cap C(\R^d)$. 
\begin{enumerate}
\item
The functions $f^{\sharp}_{\alpha}$ and $f^*_{\alpha}$ are continuous.
\item  Let furthermore 
\BEQ \beta:=\sup \{\gamma\in[0,1]\ |\ (x,y)\mapsto \frac{|f(x)-f(y)|}{|x-y|^{\gamma}}\ {\mathrm{in}} \ 
L^{\infty}_{loc} \}\in[0,1] \EEQ
be the maximum local H\"older exponent of $f$, and assume $\beta>0$.  Then $f^{\sharp}_{\alpha}$ and $f^*_{\alpha}$ are 
H\"older continuous, with H\"older exponent $\frac{\beta}{1+\beta}\in[0,\half]$.
\end{enumerate}  \label{lem:3.3}
\end{Lemma}

In particular, $f^{\sharp}_{\alpha},f^*_{\alpha}$ are $\left(\half\right)$-H\"older continuous if $f\in C^1$.

\medskip

{\bf Proof.} For the sake of the proof we choose a bounded function $\phi:B(0,1)\to \R_+$ such that $\phi(u)=\phi(|u|)$ is
strictly increasing ,
$\phi(0)=0$, $\frac{\phi(u)}{u}>2$ and  $\frac{\phi(u)}{u}\to_{u\to 0}\infty$; we assume furthermore that $\phi(u)=
o_{u\to 0}(u^{1/3})$, so that the function $\chi(u)=\chi(|u|)=\phi(u)\sqrt{\frac{\phi(u)}{u}}$
satisfies the same properties but $\frac{\chi(u)}{\phi(u)}\to_{u\to 0}\infty$. The core of the proof
is a bound on the modulus of continuity of $f^*_{\alpha}$,
$f^{\sharp}_{\alpha}$ given in terms of these two functions. We assume in the following lines that 
$|y-x|\le 1$.

\begin{itemize}
\item[(i)] Let us first prove that $f^{\sharp}_{\alpha}$ is locally bounded (from which it follows by Lemma 
\ref{eq:f*ftilde*} that $f^*_{\alpha}$ is also). Since $|y-x|\le 1$, then
\BEQ \fint_{B(y,r)} |f|\le \sup_{B(x,2)} |f|,\qquad r\le 1\EEQ
and
\BEQ (1+r^2)^{\alpha} \fint_{B(y,r)} |f|\le \left(\frac{r+1}{r}\right)^d
(1+r^2)^{\alpha} \fint_{B(x,r+1)} |f|\le 2^d f^{\sharp}(x),\qquad r>1.\EEQ
So 
\BEQ \sup_{B(x,1)}f^{\sharp}_{\alpha}\le 2^{\alpha} \max(\sup_{B(x,2)} |f|, 2^d f^{\sharp}_{\alpha}(x)). \label{eq:fsharp-bounded} \EEQ

\item[(ii)] We now obtain a modulus of continuity for  $f^{\sharp}_{\alpha}$. Fix $x\in\R^d$ and let $y$ vary in $B(x,1)$.
Consider first $r\le \phi(x-y)$. Then, letting $\tau_{x-y}f(z)=f(z-(x-y))$,
\BEA \left|\fint_{B(x,r)} |f|-\fint_{B(y,r)}|f|\right| &\le& \fint_{B(x,r)} |f-\tau_{x-y}f| \nonumber\\
&\le& \sup_{x'\in B(x,\phi(x-y)),y'\in B(y,\phi(x-y))} |f(x')-f(y')| \nonumber\\
&\le& osc_{B(x,2\phi(x-y))} (f). \EEA Since $f$ is continuous in an neighbourhood of $x$ this quantity goes to zero
when $y\to x$.

Consider now $r>\phi(x-y)$. Letting $r'=r+|x-y|$ so that $B(x,r')\supset B(y,r)$,
\BEQ \fint_{B(x,r')}|f|\ge \frac{\Vol(B(y,r))}{\Vol(B(x,r'))} \fint_{B(y,r)} |f|\EEQ
hence
\BEQ (1+r'^2)^{\alpha} \fint_{B(x,r')} |f|- (1+r^2)^{\alpha} \fint_{B(y,r)}|f|\ge \left[ \left(\frac{r}{r+|x-y|}\right)^d-1\right]f^{\sharp}_{\alpha}(y)
\gtrsim -\left(\frac{|x-y|}{\phi(x-y)}\right)f^{\sharp}_{\alpha}(y).\EEQ
Similarly, with $r''=r-|x-y|$ (note that $r''>|x-y|$ by hypothesis),
\BEQ (1+r^2)^{\alpha} \fint_{B(y,r)}|f|-(1+(r'')^2)^{\alpha} \fint_{B(x,r'')}|f|\gtrsim -\left(\frac{|x-y|}{\phi(x-y)}\right)f^{\sharp}_{\alpha}(x).\EEQ
Thus, with $M=\sup_{B(x,1)}f^{\sharp}_{\alpha}$, $M<\infty$ by (i),
\BEA && \sup_r (1+r^2)^{\alpha} \fint_{B(y,r)}|f|-\sup_r (1+r^2)^{\alpha} \fint_{B(x,r)}|f|\le \max\left( \sup_{r\le\phi(x-y)} (1+r^2)^{\alpha}
\left\{  \fint_{B(y,r)}|f|- \fint_{B(x,r)}|f|\right\},  \right.\nonumber\\
&& \qquad\qquad\qquad \left. \sup_{r>\phi(x-y)} \left\{ (1+r^2)^{\alpha} \fint_{B(y,r)}|f|-
(1+r'^2)^{\alpha} \fint_{B(x,r')}|f|
\right\} \right) \nonumber\\
&& \qquad \qquad \qquad \lesssim \max\left( \osc_{B(x,\phi(x-y))\cup B(y,\phi(x-y))}(f),M\frac{|x-y|}{\phi(x-y)}\right).\EEA
Exchanging $x$ and $y$ gives the same inequality. Hence we have shown that $f^{\sharp}_{\alpha}$ is continuous, and obtained
more precisely that, for every function $\phi$ satisfying the above hypotheses,
\BEQ \osc_{B(x,u)}f^{\sharp}_{\alpha} \lesssim \max(\osc_{B(x,2\phi(u))}(f),M\frac{u}{\phi(u)}),\qquad u\in(0,1). \label{eq:modulus}\EEQ
In particular, choosing $\phi(u)=2u^{1/(1+\beta)}$ if $\beta>0$ yields $\osc_{B(x,u)}f^{\sharp}_{\alpha}\lesssim u^{\beta/
(1+\beta)}$, so $f^{\sharp}_{\alpha}$ is $\frac{\beta}{1+\beta}$-H\"older continuous.

\item[(iii)] Let us finally obtain a modulus of continuity for  $f^*_{\alpha}$. The proof is a slightly different from (ii) because
the  support of the heat kernel is the whole space; hence we must deal with the  queue of the exponential 
$e^{-u^2/t}$ for $u\gg \sqrt{t}$. Assume first $\sqrt{t}\le\phi(x-y)$. Then
\BEQ \left| e^{t\Del}|f|(x)-e^{t\Del}|f|(y)\right|\le I_1+I_2+I_3,\EEQ
where
\BEQ I_1=\int_{B(x,\chi(x-y))} \frac{e^{-|x-z|^2/2t}}{(2\pi t)^{d/2}} |f(z)-\tau_{x-y}f(z)| dz \le osc_{B(x,\chi(x-y))
\cup B(y,\chi(x-y))} (f)\EEQ
and
\BEA && I_2=\int_{\R^d\setminus B(x,\chi(x-y))}  \frac{e^{-|x-z|^2/2t}}{(2\pi t)^{d/2}} |f(z)|dz =
\int_{|u|>\chi(x-y)} \frac{e^{-|u|^2/2t}}{(2\pi t)^{d/2}} |f(x+u)|du \nonumber\\
&&\qquad \le 2^d \int \frac{e^{-|u|^2/8t}}{(8\pi t)^{d/2}} |f(x+u)|du \ \cdot\ e^{-3\chi^2(x-y)/8t}\nonumber\\
&&\qquad \lesssim  e^{-\frac{3}{8}\frac{\phi(x-y)}{|x-y|}} f^*(x),\EEA
while $I_3=\int_{\R^d\setminus B(y,\chi(x-y))}  \frac{e^{-|y-z|^2/2t}}{(2\pi t)^{d/2}} |f(z)|dz$ is similar to $I_2$.
The exponential factor in front of $f^*$ decreases to $0$ when $y\to x$.

Assume now $\sqrt{t}>\phi(x-y)$. Then
$|x-z|^2\le |x-y|^2+|y-z|^2+2|x-y||y-z|\le (1+\eps)|y-z|^2+(1+\eps^{-1})|x-y|^2$ for every $z\in\R^d$ and $\eps>0$.
Choose $\eps=\frac{|x-y|}{\phi(x-y)}<\half$ so that $\frac{(1+\eps^{-1})|x-y|^2}{t}\lesssim\eps$. Letting 
$t'=t(1+\eps)$, one obtains
\BEA && (1+t')^{\alpha} \int \frac{e^{-|x-z|^2/2t'}}{(2\pi t')^{d/2}} |f(z)|dz - 
(1+t)^{\alpha} \int \frac{e^{-|y-z|^2/2t}}{(2\pi t)^{d/2}} |f(z)|dz \nonumber\\
&&\gtrsim -\eps (1+t)^{\alpha} \int \frac{e^{-|y-z|^2/2t}}{(2\pi t)^{d/2}} |f(z)|dz\ge -\eps f^*_{\alpha}(y).\EEA
Exchanging $x$ and $y$ gives a similar inequality, and one concludes to (\ref{eq:modulus}) as in (ii) by noting
that $e^{-\frac{3}{8}
\frac{\phi(u)}{u}}\lesssim \frac{u}{\phi(u)}$ for
$u<1$.
\end{itemize}

\hfill \eop

\medskip

A result in the same direction is

\begin{Lemma} \label{lem:tto0}
Let $f\in {\cal H}^0\cap C(\R^d)$. Then, for every $t>0$ and $0<\alpha<\half$,
\BEQ |e^{t\Del}f(x)-f(x)|\lesssim \osc_{B(x,t^{\alpha})} f+e^{-\frac{1}{4} t^{2\alpha-1}} f^*_0(x). \label{eq:tto0}\EEQ
Consequently, $e^{t\Del}f\to_{t\to 0}f$ uniformly on every compact.
\end{Lemma}

{\bf Proof.} (\ref{eq:tto0}) follows directly from the inequality
\BEQ |e^{t\Del}f(x)-f(x)|\le \int_{y\in B(x,\eps)} \frac{e^{-|x-y|^2/2t}}{(2\pi t)^{d/2}} |f(y)-f(x)| dy +
\int_{y\in B(x,\eps)^c} e^{-\eps^2/4t} \frac{e^{-|x-y|^2/4t}}{(2\pi t)^{d/2}} (|f(y)|+|f(x)|) dy.\EEQ
Taking $\eps=t^{\alpha}$ with $\alpha<\half$, using the local boundedness of $f^*_0$  (proved in the previous lemma) and letting $t\to 0$ yields the uniform
convergence on a compact set. \hfill\eop

\medskip

Finally, we shall later on need to approximate functions in ${\cal H}^0_{\alpha}$ by functions with compact support, and
use the following lemma:

\begin{Lemma} \label{lem:bump-alpha}
Let $\chi:\R^d\to\R_+$ be a smooth 'bump' scale 1 function, i.e. $\chi|_{B(0,1)}=1$, $\chi|_{\R^d\setminus B(0,2)}=0$.
Denote by $\chi_n(x)=\chi(\frac{x}{n})$ its dilatations for $n\in\N^*$. Then, if $f\in{\cal H}^0_{\alpha}$,
the functions $f_n:=f\cdot\chi_n$, $n\ge 1$ also belong to ${\cal H}^0_{\alpha}$,
 and $(f_n)^*_{\alpha}\to f^*_{\alpha}$, $(f_n)^{\sharp}_{\alpha}\to
f^{\sharp}_{\alpha}$ uniformly on every compact.
\end{Lemma}

\medskip

{\bf Proof.} Let $K\subset\R^d$ compact containing $0$. We 
  prove that $(f_n)^*_{\alpha}\to f^*_{\alpha}$ uniformly on $K$. Let $B(0,r)$ a ball containing $K$,
and assume $n\gg r$. Then $\frac{|y-x|^2}{2t}\ge \frac{|y|^2}{4t}$ for all $t,x,y$ with $t>0$, $x\in K$, $|y|>n$. 
Hence
\BEQ 0\le e^{t\Del}|f|(x)-e^{t\Del}|f_n|(x)\le 2^{d/2} \int_{|y|>n} \frac{e^{-|y|^2/4t}}{(2\pi (2t))^{d/2}} |f(y)|dy
=2^{d/2} e^{2t\Del}(|f|-|f_n|)(0),\EEQ
from which  uniform convergence follows provided simple convergence holds at $0$. But
\BEQ f^*_{\alpha}(0)=\sup_t (1+t)^{\alpha} e^{t\Del} (\sup_n |f_n|)(0)=\sup_{t,n} (1+t)^{\alpha} (e^{t\Del}|f_n|)(0)=\lim_n (f_n)_{\alpha}^*(0)\EEQ by monotone convergence.

The proof for $(f_n)^{\sharp}_{\alpha}$ is similar: let us just state that
\BEQ \fint_{B(x,R)} (|f|-|f_n|)=0 \EEQ if  $B(x,R)\subset B(0,n)$, and
\BEQ \fint_{B(x,R)} (|f|-|f_n|)\le \left(\frac{R+r}{R}\right)^d \fint_{B(0,R+r)} (|f|-|f_n|)\EEQ
otherwise. Details are left to the reader.
\hfill\eop

\bigskip
We may now finally write down our {\em pointwise parabolic estimates}: 
\begin{Lemma}[pointwise parabolic estimates]  \label{lem:pointwise-parabolic-estimates}
Let $\alpha\in[0,\frac{d}{2}]$ and $f\in{\cal H}^0_{\alpha}$.
For every $k\ge 0$,
\BEQ |\nabla^k e^{t\Del}f(x)|\lesssim t^{-\alpha-k/2} f^*_{\alpha}(x)\EEQ
and
\BEQ (\nabla^k e^{t\Del}f)^*_{\alpha}(x)\lesssim t^{-k/2}f^*_{\alpha}(x).\EEQ
\end{Lemma}

{\bf Proof.} 

By differentiating $k$ times the computations leading to (\ref{eq:H0H0}), one gets 
\BEA |\nabla^k e^{t\Del}f(x)| &\lesssim & 
 \int dr \left(\frac{r}{t}\right)^{k+1} \frac{e^{-r^2/2t}}{(2\pi t)^{d/2}} 
\int_{B(x,r)} dy |f(y)| \nonumber\\
&\lesssim & f^{\sharp}_{\alpha}(x) \int dr \left(\frac{r}{t}\right)^{k+1} \frac{e^{-r^2/2t}}{(2\pi t)^{d/2}} (1+r^2)^{-\alpha}
 \Vol(B(x,r)) \nonumber\\
&\lesssim& t^{-\alpha-k/2} f^{\sharp}_{\alpha}(x) \int dr \frac{r}{t} \frac{e^{-r^2/2t}}{(2\pi t)^{d/2}}
\Vol(B(x,r))=t^{-\alpha-k/2} f^{\sharp}_{\alpha}(x)  \label{eq:ppe1} \EEA

and 

\BEA \left(\nabla^k e^{t\Del}f\right)^{\sharp}_{\alpha}(x) &\lesssim & 
\sup_{\rho}\, (1+\rho^2)^{\alpha} \int dr \left(\frac{r}{t}\right)^{k+1} \frac{e^{-r^2/2t}}{(2\pi t)^{d/2}} \fint_{B(x,\rho)} dx'\, 
\int_{B(x',r)} dy |f(y)| \nonumber\\
&\lesssim & f^{\sharp}_{\alpha}(x) \int dr \left(\frac{r}{t}\right)^{k+1} \frac{e^{-r^2/2t}}{(2\pi t)^{d/2}}
 \Vol(B(x,r)) \nonumber\\
&\lesssim& t^{-k/2} f^{\sharp}_{\alpha}(x) \int dr \frac{r}{t} \frac{e^{-r^2/2t}}{(2\pi t)^{d/2}}
\Vol(B(x,r))=t^{-k/2} f^{\sharp}_{\alpha}(x) . \label{eq:ppe2} \EEA
To go from the first to the second inequality in (\ref{eq:ppe2}) we have made use of
the following facts which are easy to prove,
\BEQ (1+\rho^2)^{\alpha}\fint_{B(x,\rho)} dx'\, \int_{B(x',r)} dy |f(y)| \lesssim
\Vol(B(x,r))  (1+\rho^2)^{\alpha}\fint_{B(x,\rho)} dy [f(y)| \lesssim \Vol(B(x,r))
f^{\sharp}_{\alpha}(x), \qquad (\rho\gtrsim r) \EEQ
\BEQ  (1+\rho^2)^{\alpha}\fint_{B(x,\rho)} dx'\, \int_{B(x',r)} dy |f(y)| \lesssim
\Vol(B(x,r))  (1+\rho^2)^{\alpha}\fint_{B(x,r)} dy |f(y)| \lesssim \Vol(B(x,r))
f^{\sharp}_{\alpha}(x), \qquad (\rho\lesssim r). \EEQ
 \hfill\eop

\medskip
In the sequel we restrict for simplicity to the case $\alpha=0$. All results below are easily adapted to the
case $\alpha>0$ or to similar functional spaces with  pointwise bounds  of the form $|||f|||(x)=\sup_{\tau>0} 
F(\tau,e^{\tau \Del}|f|(x))$.  


\subsection{The comparison principle}

We now want to use as initial condition of (\ref{eq:2.1}) functions $h_0$ such that $h_0\in {\cal H}^{\lambda}\cap
C(\R^d)$, where
\BEQ {\cal H}^{\lambda}:=\{h_0\in L^{\infty}_{loc}(\R^d) \ |\ e^{\lambda |h_0|}\in{\cal H}^0\}. \EEQ
The  comparison to the linear heat equation (see subsection 2.1) actually suggests to consider initial
conditions in the unpleasant-looking space,
\BEQ \widetilde{\cal H}^{\lambda}:=\{h_0\in L^{\infty}_{loc}(\R^d) \ |\ e^{\lambda h_0^+},h_0^-\in{\cal H}^0\} \EEQ
However, by Jensen's inequality, $e^{\lambda (h_0^-)^*(x)}\le (e^{\lambda h_0^-})^*(x)$, so ${\cal H}^{\lambda}
\subset \widetilde{\cal H}^{\lambda}$.  
 Note that the definition is compatible with that of ${\cal H}^0$
in the previous paragraph, in the sense that $\frac{1}{\lambda}(e^{\lambda |h_0(x)|}-1)\to |h_0(x)|$ when
$\lambda\to 0$.  Also, by Jensen's inequality,  ${\cal H}^{\lambda}$ is  for $\lambda>0$ a convex subset (but not
a vector subspace) of ${\cal H}^0$, and 
\BEQ |||h_0|||_{{\cal H}^{\lambda}}(x):=\frac{1}{\lambda}\sup_{\tau>0} 
\ln\left( (e^{\tau\Del} e^{\lambda |h_0|})(x) \right)=\frac{1}{\lambda} \ln\left( (e^{\lambda|h_0|})^*(x)\right) \EEQ
defines a {\em family of pointwise "quasi-norms"}, in the sense that
\BEQ |||f|||_{{\cal H}^{\lambda}}(x)\le |||f|||_{{\cal H}^{\lambda'}}(x) \qquad (\lambda\le\lambda');\EEQ
\BEQ |||\mu f|||_{{\cal H}^{\lambda}}(x)\le |\mu| \  |||f|||_{{\cal H}^{|\mu|\lambda}} (x)\qquad (\mu\in\R)\EEQ
(the last inequality is actually an equality);
\BEQ |||f_1+f_2|||_{{\cal H}^{\lambda}}(x)\le \frac{1}{p_1}||p_1 f_1||_{{\cal H}^{\lambda}}(x)+
\frac{1}{p_2}||p_2 f_2||_{{\cal H}^{\lambda}}(x) \qquad (p_1,p_2\ge 1,\ \frac{1}{p_1}+\frac{1}{p_2}=1). \label{eq:3:Holder} \EEQ

  We then expect the solution of (\ref{eq:2.1}) to be "uniformly bounded
in ${\cal H}^{\lambda}$", at least locally in time (thus allowing for further generalizations to 
equations with time-dependent coefficients), and thus  to lie for all $T>0$ in the functional space
\BEQ {\cal H}^{\lambda}([0,T]):=\{h\in L^{\infty}_{loc}(\R\times\R^d)\ |\ \forall x\in\R^d, 
\sup_{t\in[0,T]} (e^{\lambda |h_t|})^*(x)<\infty\}.\EEQ

As mentioned previously, the comparison principle in its different forms usually requires as a cornerstone
assumption the boundedness of the solutions. However, various authors have proved ad hoc comparison principles for
PDE's with unbounded coefficients; the solution lies in functional spaces including functions growing at infinity.
The KPZ equation is a very particular class of Hamilton-Jacobi-Bellman equations for which a comparison principle
holds under {\em quadratic growth conditions}, see Ito \cite{Ito}, Da Lio-Ley \cite{DaLioLey,DaLioLey2}.  Lemma 2.1 and Theorem
2.1 in \cite{DaLioLey} state the following in our case:

\begin{Proposition}  \cite{DaLioLey}

Let $\underline{U}\in USC([0,T]\times\R^d)$ (resp. $\bar{U}\in LSC([0,T]\times\R^d)$) be a viscosity sub-solution (resp. super-solution)
of (\ref{eq:2.1}). Assume there exists $C>0$ such that $|\underline{U}(t,x)|,|\bar{U}(t,x)|\le C(1+|x|^2)$ for all $x\in\R^d,t\le T$.
Then $\underline{U}\le \bar{U}$ in $[0,T]\times\R^d$.

\end{Proposition}

A continuous function $h_0$ with quadratic growth at infinity, $|h_0(x)|\lesssim 1+|x|^2$, is
in general not in  ${\cal H}^{\lambda}$ for any $\lambda\ge 0$. Conversely, a function in ${\cal H}^{\lambda}$, $\lambda\ge 0$ may grow arbitrarily fast   in small
domains $\Omega_n$, $n\to\infty$ with $d(0,\Omega_n)\to_{n\to\infty} \infty$ provided the Lebesgue measure of $\Omega_n$
decreases to zero fast enough.  On the other hand, since the supremum of $n$ i.i.d. random variables grows like
$O(\sqrt{\log n})$,  one does expect random initial data $h_0$ to have a.s. quadratic growth at infinity.
Actually, if $\eps>0$, then a.s.  a random initial condition grows more slowly at infinity than $|x|^{\eps}$.
Thus the above comparison principle holds for such data, and the existence of a sub-solution and a super-solution in 
this class of functions entails by Perron's method the existence and unicity of a viscosity solution of (\ref{eq:2.1}).

It seems however much more natural  in our setting to  prove a comparison principle for functions in 
${\cal H}^{\lambda}$ since
the bounds one expects for the solution will depend on the pointwise maximal estimates $(h_0^-)^*$ and 
$(e^{\lambda h_0^+})^*$ (on the contrary, solutions are expected to have a finite explosion time for initial conditions
with quadratic growth, showing that this is in some sense too large a functional space). As it happens, 
we get such a comparison principle, but only
for solutions in spaces ${\cal H}^{\lambda'}$ with parameter  $\lambda'\ge 2\lambda$ (our proof does not hold for $\lambda'=\lambda$). In some sense ${\cal H}^0$ is the largest natural
functional space for globally defined solutions of parabolic PDE's. We conjecture that this extension of the viscosity solution theory to spaces modelled after ${\cal H}$ (like  ${\cal H}^{\lambda}$ in
the present case) is valid and of interest not only for the KPZ equation, but  probably much beyond for
many nonlinear parabolic PDE's.

\medskip

\noindent Let us state our first main theorem, Theorem 1 in the Introduction, following closely the strategy of Da Lio and Ley:

\begin{Theorem}[comparison principle] \label{th:comparison}

Let $\underline{U}\in USC([0,T]\times\R^d)\cap{\cal H}^{2\lambda}([0,T])$ (resp. $\bar{U}\in LSC([0,T]\times\R^d)\cap {\cal H}^{2\lambda}([0,T])$) be a viscosity sub-solution (resp. super-solution)
of (\ref{eq:2.1}). 
Then $\underline{U}\le \bar{U}$ in $[0,T]\times\R^d$.

\end{Theorem}

The proof is very similar to  \cite{DaLioLey}, section 2. The essential element is the following lemma.

\begin{Lemma} \label{lem:U-muV}
Let $\underline{U}\in USC([0,T]\times\R^d)\cap{\cal H}^{2\lambda}([0,T])$ be a sub-solution, and $\bar{U}\in
LSC([0,T]\times\R^d)\cap{\cal H}^{2\lambda}([0,T])$ be a super-solution of (\ref{eq:2.1}). Then
$\Psi_{\mu}:=\underline{U}-\mu \bar{U}$, $\mu\in(0,1)$ is  a sub-solution of the quadratic KPZ equation,
\BEQ \partial_t \psi=\Del\psi+\frac{\lambda}{1-\mu}|\nabla \psi|^2.\EEQ
\end{Lemma}

Note that $\left(\underline{U},\bar{U}\in {\cal H}^{2\lambda}([0,T]) \right)\Longrightarrow \Psi_{\mu}
\in{\cal H}^{\lambda}([0,T])$ by (\ref{eq:3:Holder}) (hence our choice of parameter, $\lambda'\ge 2\lambda$, see discussion above).

\medskip

{\bf Proof.}

If $U,V\in C^{1,2}$ then the proof is elementary. First,
\BEQ \partial_t \Psi_{\mu}\le \Del \Psi_{\mu}+\lambda(V(\nabla\underline{U})-\mu V(\nabla\bar{U})).\EEQ
Then, since $V$ is convex, 
\BEQ V(a)\le \mu V(b)+(1-\mu)V(\frac{a-\mu b}{1-\mu}),\qquad a,b\in\R^d. \label{eq:V-convex-mu} \EEQ
Finally, applying this inequality to $a=\nabla\underline{U}$, $b=\nabla\bar{U}$, and using Assumption \ref{assumptions}
(4), $V(y)\le y^2$, yields the result.

\medskip

Otherwise the proof is essentially a very particular case of \cite{DaLioLey}, Lemma 2.2. Let us reproduce the main
arguments for the sake of the reader. Let $\psi\in C^2([0,T]\times\R^d)$ and $(\bar{t},\bar{x})$ a strict local maximum
of $\Psi-\psi$; we must prove that $\partial_t \psi(\bar{t},\bar{x})\le \nu\Del\psi(\bar{t},\bar{x})+
\frac{\lambda}{1-\mu}|\nabla\psi(\bar{t},\bar{x})|^2$. This is done by the standard doubling of variables argument, namely, we let
$\Theta(t,x,y):=\psi(t,x)+\frac{|x-y|^2}{\eps^2}$, and $M_{\eps}=(\Psi-\Theta)(t_{\eps},x_{\eps})$ be the
maximum of $\Psi-\Theta$ in a small ball centered at $(\bar{t},\bar{x})$; it is known that $|x_{\eps}-y_{\eps}|=o(\eps)$
and $M_{\eps}\to_{\eps\to 0}\Psi(\bar{t},\bar{x})-\psi(\bar{t},\bar{x})$. By Theorem 8.3 in the User's guide, see
 in \cite{DaLioLey} for the details of computations, one finds, exploiting the hypotheses on $\underline{U}, \bar{U}$,
\BEQ \partial_t \psi(t_{\eps},x_{\eps})+H(t_{\eps},x_{\eps},\nabla\psi(t_{\eps},x_{\eps})+p_{\eps},X)-\mu 
H(t_{\eps},y_{\eps},\frac{p_{\eps}}{\mu},\frac{Y}{\mu})\le 0 \EEQ
where $p_{\eps}=2\frac{x_{\eps}-y_{\eps}}{\eps^2}$, $H(x,t,p,X):=-\lambda V(p)-\nu\Tr(X)$, and $X,Y$ are symmetric $d\times d$ matrices, depending on $\eps$ and on a parameter $\rho>0$, such that
$\Tr(X-Y)\le \Del\psi(t_{\eps},x_{\eps})+O(\rho/\eps^4).$
Hence
\BEQ \partial_t\psi(t_{\eps},x_{\eps})\le \nu\Del\psi(t_{\eps},x_{\eps})+ \lambda\left[
V(\nabla\psi(t_{\eps},x_{\eps})+p_{\eps})-\mu V(\frac{p_{\eps}}{\mu})\right]+O(\rho/\eps^4).\EEQ
Letting $\rho\to 0$ and using (\ref{eq:V-convex-mu}) as above yields
\BEQ \partial_t\psi(t_{\eps},x_{\eps})\le\nu\Del\psi(t_{\eps},x_{\eps})+\frac{\lambda}{1-\mu}|\nabla\psi(t_{\eps},
x_{\eps})|^2.\EEQ
Finally, letting $\eps\to 0$ gives the result. \hfill\eop

\medskip

We shall also need a non-standard comparison lemma for the linear heat equation:

\begin{Lemma} \label{lem:comparison}
Let $\underline{U}\in USC([0,T]\times\R^d)\cap{\cal H}^0([0,T])$ (resp. $\bar{U}\in LSC([0,T]\times\R^d)\cap
 {\cal H}^0([0,T])$) be a viscosity sub-solution (resp. super-solution)
of the linear heat equation. 
Then $\underline{U}\le \bar{U}$ in $[0,T]\times\R^d$.

\end{Lemma}

\noindent In other words, Theorem \ref{th:comparison} holds for $\lambda=0$.
\smallskip

{\bf Proof.} Since the equation is linear, we may (by replacing $\underline{U}$ with 
$(\underline{U}-\bar{U})-e^{t\Del} (\underline{U}_0-\bar{U}_0)$) assume that $\bar{U}=0$ and $\underline{U}_0=0$.
Now, we have no bound at infinity available for $\underline{U}$, and the classical maximum principle does not hold.
Instead we choose a smooth function $\chi\ge 0$ with $\chi\big|_{(-\infty,0]}\equiv 1$, $\supp(\chi)\subset (-\infty,1]$, define
$\chi_n(x):=\chi(|x|-n)$ 
 and obtain the following inequality for $\underline{U}_n:=\underline{U}\chi_n$,
\BEQ (\partial_t-\Del)\underline{U}_n+\Del\chi_n \underline{U}+2\nabla\chi_n\cdot \nabla\underline{U}\le 0.
\label{eq:Un} \EEQ

Assume that $\underline{U}\in C^{1,2}$ is a classical sub-solution to begin with. Then, since $\underline{U}_n$,
$\Del\chi_n \underline{U}$ and $\nabla\chi_n\cdot \nabla\underline{U}$ are bounded, the classical comparison
principle entails
\BEA \underline{U}_n(t,x) & \lesssim & -\int_0^t ds e^{(t-s)\Del}(\Del\chi_n \underline{U}(s))-2
\int_0^t ds e^{(t-s)\Del}(\nabla\chi_n\cdot \nabla\underline{U}(s)) \nonumber\\
&=& \int_0^t ds e^{(t-s)\Del}(\Del\chi_n \underline{U}(s)) -2\int_0^t ds \nabla e^{(t-s)\Del}\cdot
(\nabla\chi_n \underline{U}(s))\nonumber\\
&\lesssim&
\int_0^t ds e^{(t-s)\Del} (|\tilde{\chi}_n\underline{U}(s)|) + 
\int_0^t ds \left| \nabla e^{(t-s)\Del} (\nabla\chi_n\underline{U}(s)) \right|
, \label{eq:3.41} \EEA

where $\tilde{\chi}_n=\max(|\nabla\chi_n|,|\Del\chi_n|)$. Now $\sum_n \tilde{\chi}_n\lesssim 1$, so (by the
pointwise parabolic estimates)
\BEQ \sum_n \int_0^t ds e^{(t-s)\Del} |\tilde{\chi}_n \underline{U}(s)|(x)\lesssim
\int_0^t ds e^{(t-s)\Del}|\underline{U}(s)|(x)\lesssim T \sup_{s\in[0,T]} (\underline{U}(s))^*(x). \label{eq:UUU}\EEQ
Hence (for $x$ fixed) $\int_0^t ds e^{(t-s)\Del}|\tilde{\chi}_n\underline{U}(s)|\to_{n\to\infty} 0.$
Lemma  \ref{lem:pointwise-parabolic-estimates} yields the same bound as (\ref{eq:UUU}), with $T$ replaced by $\sqrt{T}$,
for the term with the gradient.

\medskip
The above proof does not seem to extend to fonctions in $USC([0,T]\times\R^d)\cap{\cal H}^0([0,T])$ by a
density argument (in particular, if $\chi$ is a smooth, positive 'bump' function, then $\chi\ast \underline{U}$ is a
smooth subsolution in the classical sense if $\underline{U}$ is since $(\partial_t-\Del)(\chi\ast \underline{U})=\chi\ast
(\partial-\Del)\underline{U}\le 0$,  but not in the viscosity sense in general if $\underline{U}$ is only upper-semicontinuous). Instead we use another truncation argument, which could also have been used in the classical case.
We fix $x\in\R^d$ and let $n\to\infty$ as above.
Since $\underline{U}\in{\cal H}^0$, it is locally bounded, so the function $\underline{U}$ is a bounded sub-solution of the heat equation on $[0,T]\times B(0,n+1)$. Thus the
classical maximum principle and Green's formula imply that
\BEQ \underline{U}(t,x)\le \int_0^t ds \int_n^{n+1}dr\, \int_{\partial B(0,r)} \nabla_{\vec{n}}G_r(t,x;s,y) \underline{U}(s,y),
\qquad t\le T, x\in B(0,n)\EEQ
where $\nabla_{\vec{n}}$ is the normal derivative and  $G_r$ is the Green function of the heat equation on $\R_+\times B(0,r)$. By standard estimates, $|\nabla G_r(t,x;s,y)|\lesssim (t-s)^{-\half} G(t,x;s,y)$ if $y\in \partial B(0,r)$, where
$G(t,x;s,y)$ is the usual heat kernel on $\R\times\R^d$. One has thus obtained an estimate very similar to
(\ref{eq:3.41}), and the end of  the proof is the same.

 \hfill\eop

\bigskip

{\bf Proof of Theorem \ref{th:comparison}.}

By Lemma \ref{lem:U-muV}, $\partial_t \Psi_{\mu}\le \Del\Psi_{\mu}+\frac{\lambda}{1-\mu}|\nabla\Psi_{\mu}|^2$. Equivalently,
$(\partial_t-\Del)\left(e^{\frac{\lambda}{1-\mu}\Psi_{\mu}}\right)\le 0$. By Lemma \ref{lem:comparison}, 
\BEQ \Psi_{\mu}(t,x)\le \frac{1-\mu}{\lambda} \ln \left[ \int \frac{e^{-|x-y|^2/2t}}{(2\pi t)^{d/2}} e^{\lambda h_0(y)}
dy \right].\EEQ
Letting $\mu\to 1$, one finds $\Psi_{\mu}\le 0$. \hfill\eop


\subsection{Bounds for the solution}


 Let $h_0\in{\cal H}^{2\lambda}\cap C(\R^d)$. Then $\underline{h}_t:=-e^{t\nu\Del}h_0^-$ is a sub-solution, and
$\bar{h}_t:=\frac{1}{\lambda}\ln(e^{t\nu\Del} e^{\lambda h_0^+})$  a super-solution of (\ref{eq:2.1}), and
the pointwise parabolic estimates, together with Jensen's inequality, imply that
$\underline{h},\bar{h}\in C([0,T]\times\R^d)\cap {\cal H}^{2\lambda}([0,T])$.
Perron's method (see User's guide \cite{CraIshLio}, Theorem 4.1), in combination with the comparison principle
of the previous paragraph, shows that 
\BEQ h(x):=\sup\{\tilde{h}(x)\ |\ \underline{h}\le \tilde{h}\le\bar{h}\ {\mathrm{and}}\ \tilde{h}\ 
{\mathrm{is\ a\ subsolution\ of}}\ (\ref{eq:KPZdet})\}  \EEQ
is the unique viscosity solution in $C([0,T]\times\R^d)\cap {\cal H}^{2\lambda}([0,T])$ of (\ref{eq:2.1}) for every $T>0$. We simply call $h$ the solution on $[0,T]$ of (\ref{eq:KPZdet}) with initial condition $h_0$.

Anticipating on section 4, the analogue of the space ${\cal W}^{1,\infty}$ in our setting is

\begin{Definition} 
\BEQ {\cal W}^{1,\infty;2\lambda}_{j}:=\left\{h_0\in {\cal W}^{1,\infty}_{loc}\ |\ \locsup^j h_0\in {\cal H}^{2\lambda}, \ 
 2^{j/2} \locsup^j |\nabla h_0| \in {\cal H}^{2\lambda} \right\} \qquad (j\ge 0), \EEQ
\end{Definition}
 where $\locsup^j$ (the "scale $j$ local supremum") operates on functions in $L^{\infty}_{loc}$ in the following way, 
\BEQ \locsup^{j} f(x):=\sup_{y\in B(x,2^{j/2})} |f(y)|. \label{eq:locsup} \EEQ
Apparently, it is necessary to consider an initial condition $h_0$ such that
$h_0$ and $\nabla h_0$ are locally bounded (see
proof of Lemma \ref{lem:classical}) if one wants the solution to exist; these conditions are of course automatically
verified if $h_0$ is in $C^1$. Assuming the {\em local suprema} $\locsup^j h_0$, 
$\locsup^j |\nabla h_0|$ to be in ${\cal H}^{2\lambda}$ ensures that the spaces
${\cal W}^{1,\infty;\lambda}$ are stable under the flow (well, not quite, see Lemma
\ref{lem:classical} for an exact statement).
 
   The value of $j$ is at this point arbitrary, and it is quite possible to  take $j=0$. 
However the scaling is important starting from section 4, so we chose to let the dependence on $j$ explicit.

Just like  ${\cal H}^{2\lambda}$ before, ${\cal W}^{1,\infty;2\lambda}_j$ is a convex subset
of $C(\R^d)$, and
\BEQ |||h_0|||_{{\cal W}^{1,\infty;2\lambda}_j}(x):=\max\left(   |||\, \locsup^j h_0\, |||_{{\cal H}^{2\lambda}}(x),\,   |||\, 2^{j/2} \locsup^j |\nabla h_0|\,  |||_{{\cal H}^{2\lambda}}(x) \right)
\label{eq:Wj1inftylambda}  \EEQ
defines now a family of  {\em local quasi-norms}, as pointed out in the Introduction.
At this point we must explain clearly why we distinguish {\em pointwise} quasi-norms from
{\em local quasi-norms}. A function $f$ such that $(\locsup^j f)^*(x)<\infty$ cannot have
arbitrarily large fluctuations, contrary to a function $f$ satisfying simply $f^*(x)<\infty$
(recall the family of examples from section 3.1). Namely (choosing $j=0$ for simplicity), if $(\locsup^0 f)^*(x)<\infty$, then
$(1+|y-x|)^{-d} |f(y)| \lesssim e^{(y-x)^2\Del}\locsup^0(f)(x)$ is bounded uniformly in $y$,
hence $|f(y)|=O(|y-x|^d)$ grows at most polynomially. Ultimately this comes from the
fact that the integral $e^{(y-x)^2\Del} \locsup^0(f)(x)$ is equivalent (up to a multiplicative factor $O(1)$)
 to a weighted sum of the values of $\locsup^0(f)$ over cells of the unit lattice: one gets
 {\em local} estimates instead of {\em pointwise} estimates.  In particular, $|||\locsup^0 h_0|||_{{\cal H}^{2\lambda}}(x)<\infty$ (see \ref{eq:Wj1inftylambda}) implies: $|f(y)|=O_{|y|\to\infty}
(\log|y|)$ (which holds if $f=\eta$ is a regularized white noise).

On the other hand, bounds also hold if one replaces the local quasi-norm $|||\locsup^j h_0|||_{{\cal H}^{2\lambda}}(x)$ in (\ref{eq:Wj1inftylambda}) by the smaller local supremum $\locsup^j( |||h_0|||_{{\cal H}^{2\lambda}})(x)$, and thus go back to our previous {\em pointwise} "quasi-norms". Note that this quantity is finite as soon as $h_0\in {\cal H}^{2\lambda}\cap
L^{\infty}_{loc}$ (see proof of Lemma (\ref{lem:3.3}) (i)).  There is also a way to move out the  local supremum of the quasi-norm for $\nabla h_0$ in (\ref{eq:Wj1inftylambda}) -- see next
subsection for details. Then of course one can only prove bounds for the corresponding
(local supremum of) pointwise quasi-norms at time $t$. Despite allowing more general
initial conditions (i.e. with arbitrary large fluctuations), this has the inconvenient 
of complicating the statements. For applications to inhomogeneous KPZ equation with random
forcing in section 4, the "local quasi-norm" version (\ref{eq:Wj1inftylambda}) will suffice.

\medskip

\medskip

A first easy result is:

\begin{Lemma} \label{lem:easy}
Let $h$ be the viscosity solution on [0,T] of (\ref{eq:2.1}) with initial condition $h_0\in{\cal H}^{2\lambda}\cap 
C(\R^d)$. Then,
for every $t\in[0,T]$,
\BEQ (e^{a\lambda |h_t|})^*(x)\le (e^{a\lambda |h_0|})^*(x),\qquad a\ge 1.\label{eq:ut*bound}\EEQ
In particular, 
\BEQ |h_t(x)|\le |||h_0|||_{{\cal H}^{\lambda}}(x).\EEQ
\end{Lemma}

{\bf Proof.} 
By the comparison principle, Theorem \ref{th:comparison}, $|h|\le u$, where $u$ is the solution of the quadratic
KPZ equation $\partial_t u=\Del u+\lambda|\nabla u|^2$ with initial condition $|h_0|$. Then $e^{\lambda u}$
is a solution of the linear heat equation, hence  (by Jensen's inequality and pointwise parabolic estimates)
\BEQ \left(e^{a\lambda |h_t|} \right)^*(x)\le \left(e^{a\lambda u_t}\right)^*(x)=
\left( (e^{t\Del} e^{\lambda |h_0|})^a\right)^*(x)\le \left(e^{t\Del}e^{a\lambda |h_0|}\right)^*(x)
\le (e^{a\lambda |h_0|})^*(x) \label{eq:3.50} \EEQ   for $a\ge 1$. 
\hfill\eop

Thus $h$ extends to $t\in\R_+$ and satisfies (\ref{eq:ut*bound}) for arbitrary $t$.

\medskip

Note that (\ref{eq:ut*bound}) still holds true when one inserts local suprema: for
$a\ge 1$, 
\BEA \left(e^{a\lambda\,  \locsup^j h_t}\right)^*(x) &=&\sup_{\tau>0} e^{\tau\Del}
\locsup^j e^{a\lambda|h_t|}(x) \nonumber\\
&\le & \sup_{\tau>0} e^{\tau\Del} \locsup^j e^{t\Del} e^{a\lambda|h_0|}(x) 
\nonumber\\ & \le & \sup_{\tau>0} e^{(\tau+t)\Del} e^{a\lambda \, \locsup^j
h_0}(x)\le \left( e^{a\lambda\, \locsup^j h_0}\right)^*(x). \label{eq:3.66} \EEA

\smallskip


\subsection{Bounds for the gradient}


Let  $h$  be the (viscosity) solution of (\ref{eq:2.1}) with
initial condition $h_0\in{\cal W}^{1,\infty;2\lambda}_{j}\cap  C^2$.

\medskip

 The main
task in this subsection is to prove {\em a priori bounds} for the discrete gradient, \\ 
$\sup_{\vec{\eps},\vec{\eps'}\in B(0,1)}
\frac{|h_t(x+2^{j/2}\vec{\eps})-h_t(x+2^{j/2}\vec{\eps}')|}{|\vec{\eps}-\vec{\eps}'|}$.
If $h_t$ is differentiable, then this quantity is equal to $2^{j/2}\locsup^j |\nabla h_t|(x)$. Proving differentiability can then be done using a cut-off argument as follows.
Let $\chi$ be a 'bump' function as in 
Lemma \ref{lem:bump-alpha},  $\chi^{(L)}(x):=\chi(\frac{x}{L})$ ($L\in\N^*$), and
$h^{(L)}$ be the solution of the homogeneous KPZ equation with initial condition $h^{(L)}_0=h_0(x)\chi^{(L)}(x)$. From Lemma \ref{lem:easy}, one knows that the sequence
$(h^{(L)})_{L\ge 1}$ is locally uniformly bounded, i.e.
for every compact $K\subset\R^d$, and every $T>0$, $\sup_{L\ge 1}\sup_{t\le T}\sup_K |h_t^{(L)}|<C(K,T)$. Since $h^{(L)}$ is compactly supported and $C^2$, $h^{(L)}$ is
known to be classical, hence a priori bounds for the discrete gradient of $h_t^{(L)}$ hold
ipso facto for the local supremum of its gradient, $2^{j/2}\locsup^j |\nabla h^{(L)}_t|(x)$. Now we use the following argument:

\begin{Lemma} \label{lem:h(n)}
Assume $(h^{(L)})_{L\ge 1}$ is locally uniformly differentiable, i.e. for all
compact $K\subset \R^d$, and all $T>0$, $\sup_{L\ge 1}\sup_{t\le T}\sup_K |\nabla h_t^{(L)}|<C(K,T)$.
 Then $h^{(L)}$, resp. $\nabla h^{(L)}$  converges to $h$, resp. $\nabla f$ uniformly on every compact. The function $h$ is a classical solution of the KPZ equation.
\end{Lemma}

{\bf Proof.} 
By Ascoli's theorem and the classical diagonal extraction procedure, one may construct a subsequence $h^{(L_m)}$
converging locally uniformly. By the  stability principle for continuous viscosity solutions (see e.g.
\cite{Bar}, Theorem 3.1), the limit is a solution of the KPZ equation with initial condition $h_0$. Since the
solution is unique, we have shown that $h^{(L_m)}\to_{m\to\infty} h$ in $C(\R_+\times\R^d)$. Since the
sequence $(h^{(L)})_L$ is pre-compact and all subsequences converge to $h$, the sequence $(h^{(L)})$ itself converges
to $h$.

Now Schauder estimates applied first to the equations $(\partial_t-\Del)h^{(L)}=\lambda
V(\nabla h^{(L)})$, and then to the equations $(\partial_t-\Del)(h^{(L)}-h^{(L')})=\lambda (V(\nabla h^{(L)})-V(\nabla h^{(L')}))$, imply that
$\nabla^2 h^{(L)}$ are locally uniformly bounded, and then (using the uniform convergence
of the sequence $(h^{(L)})$  on every compact) that the gradient sequence
$(\nabla h^{(L)})$ also converge uniformly on every compact. Then the standard
local existence theory for the KPZ equation implies that $h$ is a classical solution.
 \hfill\eop

\bigskip

We may now come back to a priori bounds. 
As mentioned above, there is a "local quasi-norm" version, and a "pointwise quasi-norm" version.
We concentrate on the "local " version, and then sketch a derivation of the
"pointwise" version.

\begin{Lemma} \label{lem:classical} (see Introduction)
\begin{itemize}
\item[(i)]
Assume  $h_0\in{\cal W}^{1,\infty;2\lambda}_j$. Then the solution $h$ is classical for $t>0$.
Furthermore,
\BEQ |\nabla h(t,x)|\le 4 |||h_0|||_{{\cal W}_j^{1,\infty;2\lambda}}(x).
\label{eq:3.67} \EEQ
\item[(ii)]
 (same hypothesis)  Then $h_t\in {\cal W}_{j}^{1,\infty;2\lambda/5}$ and
\BEA  |||\, 2^{j/2}\locsup^j |\nabla h_t|\, |||_{{\cal H}^{2\lambda/5}} &\le &  4\, 
 |||\locsup^j h_0|||_{{\cal H}^{2\lambda}}(x)+ \,   |||\, 2^{j/2} \locsup^j |\nabla h_0|\, |||_{{\cal H}^{2\lambda}}(x) \nonumber\\
 &\le & 5 |||h_0|||_{{\cal W}^{1,\infty;2\lambda}}(x)
.\EEA
\end{itemize}
\end{Lemma}

{\bf Proof.} 

\begin{itemize}
\item[(i)] Let $\vec{\eps}\in B(0,1)\setminus\{0\}$. We introduce the following notations,
\BEQ \del_{\eps}^j f(\cdot):=\frac{f(\cdot+2^{j/2}\vec{\eps})-f(\cdot)}{|\vec{\eps}|},\qquad 
\tilde{\del}_{\eps}^jf(\cdot):=\frac{f(\cdot+2^{j/2}\vec{\eps})-(1-|\vec{\eps}|)f(\cdot)}{|\vec{\eps}|}.\EEQ
Note that 
\BEQ \tilde{\del}^j_{\eps}f=\del^j_{\eps} f+f \label{eq:deltilde-del}. \EEQ
  By Lemma \ref{lem:U-muV},  $|\vec{\eps}|
\tilde{\del}^j_{\vec{\eps}}h_t(\cdot)=
h(t,\cdot+2^{j/2}\vec{\eps})-(1-|\vec{\eps}|) h(t,\cdot)$ is a sub-solution of the KPZ equation $\partial_t \psi=\Del\psi+\frac{\lambda}{|\vec{\eps}|}
|\nabla\psi|^2$, hence 

\BEQ \tilde{\del}^j_{\vec{\eps}} h_t(x)\le \frac{1}{\lambda} \log e^{t\Del} (e^{\lambda \tilde{\del}_{\vec{\eps}}h_0})(x).\EEQ
On the other hand, exchanging the r\^oles of $x$ and $x+2^{j/2}\vec{\eps}$, 
\BEA  \tilde{\del}^j_{\vec{\eps}}h_t(x) &=& -\frac{1-|\vec{\eps}|}{|\vec{\eps}|} \left(h_t(x)-
\frac{1}{1-|\vec{\eps}|} h_t(x+2^{j/2}\vec{\eps}) \right) \nonumber\\
&=& - (1-|\vec{\eps}|) \tilde{\del}^j_{-\vec{\eps}} h_t(x+2^{j/2}\vec{\eps}) + (2-|\vec{\eps}|)
h_t(x+2^{j/2}\vec{\eps})  \label{eq:3.72} \EEA
hence the two-sided bound,
\BEA &&  |\tilde{\del}^j_{\vec{\eps}} h_t(x)| \le \frac{1}{\lambda} \log e^{t\Del}
\left(\exp \lambda \sup_{\vec{\eps},\vec{\eps}'\in B(0,1)} \frac{|h_0(2^{j/2}\vec{\eps}+\cdot)-
(1-|\vec{\eps}-\vec{\eps}'|)h_0(2^{j/2} \vec{\eps}'+\cdot)|}{|\vec{\eps}-\vec{\eps}'|} \right)(x) + 2\,\locsup^j h_t(x)
\nonumber\\
&& \ \  \le  \frac{1}{2\lambda} \log e^{t\Del}
\left(\exp 2\lambda \sup_{\vec{\eps},\vec{\eps}'\in B(0,1)} \frac{|h_0(2^{j/2}\vec{\eps}+\cdot)-
h_0(2^{j/2}\vec{\eps}'+\cdot)|}{|\vec{\eps}-\vec{\eps}'|} \right)(x) + \frac{3}{2\lambda} \log e^{t\Del} \exp( 2\lambda\,  \locsup^j h_0)(x)  \nonumber\\
&&  \ \  \le  \frac{1}{2\lambda} \log e^{t\Del}
\left(\exp 2\lambda \, \cdot\, 2^{j/2} \locsup^j |\nabla h_0| \right)(x) + \frac{3}{2\lambda} \log e^{t\Del} \exp( 2\lambda\,  \locsup^j h_0)(x)  \nonumber\\ \label{eq:3.71} \EEA

From this we deduce in particular the pointwise estimate, 
\BEQ |\del^j_{\vec{\eps}} h_t(x)|\le |h_t(x)|+ |\tilde{\del}^j_{\vec{\eps}}h_t(x)|\le  5\, |||h_0|||_{{\cal W}^{1,\infty;2\lambda}_j}(x) \EEQ
which is uniform in $\vec{\eps}$, and also
\BEQ \left(e^{\frac{1}{2}\lambda |\tilde{\del}_{\vec{\eps}^j}h_t|}\right)^*(x)  \le 
e^{\frac{1}{2}\lambda\, |||\, 2^{j/2} \locsup^j |\nabla h_0| \, |||_{{\cal H}^{2\lambda}}(x)} e^{\frac{3}{2} \lambda\, |||\, \locsup^j h_0|||_{{\cal H}^{2\lambda}}(x)}.\EEQ

Applying the above arguments to $h^{(L)}$, $L\ge 1$, and letting $\vec{\eps}\to 0$, one  obtains the first estimate (\ref{eq:3.67}) for $h^{(L)}$, with $|||h_0|||_{{\cal W}_j^{1,\infty;2\lambda}}(x)$ replaced by $|||h_0^{(L)}|||_{{\cal W}_j^{1,\infty;2\lambda}}(x)$
. The latter quantity is bounded  uniformly in $L$, implying that $h$ is classical
by Lemma \ref{lem:h(n)}. In particular, $h$ is differentiable, so we have actually proved
(\ref{eq:3.67}).

\item[(ii)]

Letting $x$ move around in the ball $B(x,2^{j/2})$ we see that the bound (\ref{eq:3.71})
is also valid for\\  $F(t,x):=\sup_{\vec{\eps},\vec{\eps'}\in B(0,1)} \frac{|h_t(x+2^{j/2}\vec{\eps})-
(1-|\vec{\eps}-\vec{\eps}'|)
h_t(x+\vec{\eps'})|}{|\vec{\eps}-\vec{\eps}'|}$. Hence (applying H\"older's inequality
with conjugate exponents $(p,q)=(\frac{5}{4},5)$)

\BEA e^{\frac{2}{5}\lambda |||\, 2^{j/2}\,  \locsup^j |\nabla h_t|\, |||_{{\cal H}^{2\lambda/5}}(x)}  &\le& \left[ \left( e^{\frac{1}{2}\lambda F(t,\cdot)} \right)^*(x)\right]^{4/5} \left[ \left( e^{2\lambda\,  \locsup^j h_t}\right)^*(x) \right]^{1/5} \nonumber\\
&\le&  e^{\frac{2}{5} \lambda ||| 2^{j/2} \locsup^j  |\nabla h_0|\, |||_{{\cal H}^{2\lambda}}(x)}
e^{\frac{8}{5}\lambda  |||\locsup^j h_0|||_{{\cal H}^{2\lambda}}(x)}, \EEA
whence the result.

\end{itemize}
 \hfill \eop

\bigskip

Let us now briefly explain how to derive the weaker "pointwise quasi-norm" version of this bound.
Leaving the supremum over $\vec{\eps},\vec{\eps}'$ outside of the heat kernel, one obtains
instead of (\ref{eq:3.71})
\BEA  && e^{\half \lambda |\tilde{\del}^j_{\vec{\eps}} h_t|}(x)  =  
\max \left( e^{\half \lambda \tilde{\del}^j_{\vec{\eps}} h_t(x)}  \, ,\,  
 e^{-\half \lambda \tilde{\del}^j_{\vec{\eps}} h_t(x)} \right)
\nonumber\\
&&\quad \le  \max \left( \left( e^{t\Del}(e^{\lambda|\tilde{\del}_{\vec{\eps}}^j h_0|})(x) \right)^{1/2}  \, ,\,  
\left( e^{t\Del} ( e^{\lambda|\tilde{\del}_{-\vec{\eps}}^j h_0|})(x+2^{j/2} \vec{\eps})
\right)^{1/2} \left( e^{t\Del} (e^{2\lambda |h_0|})(x+2^{j/2} \vec{\eps}) \right)^{1/2}
\right)
\nonumber\\
&&\quad \le  \max\left( \left( e^{t\Del}(e^{2\lambda|{\del}_{\vec{\eps}}^j h_0|})(x) \right)^{1/4}
\left( e^{t\Del}(e^{2\lambda |h_0|})(x)\right)^{1/4} \, ,\,  
\left( e^{t\Del} ( e^{2\lambda |{\del}_{-\vec{\eps}}^j h_0|})(x+2^{j/2} \vec{\eps})
\right)^{1/4} \left( e^{t\Del} (e^{2\lambda |h_0|})(x+2^{j/2} \vec{\eps}) \right)^{3/4}
\right) \nonumber\\  \label{eq:pointwise-1}
\\
&&\quad \le   \left( e^{t\Del}(e^{2\lambda|{\del}_{\vec{\eps}}^j h_0|})(x) \right)^{1/4}
\left( e^{t\Del}(e^{2\lambda |h_0|})(x)\right)^{1/4} +
\left( e^{t\Del} ( e^{2\lambda |{\del}_{-\vec{\eps}}^j h_0|})(x+2^{j/2} \vec{\eps})
\right)^{1/4} \left( e^{t\Del} (e^{2\lambda |h_0|})(x+2^{j/2} \vec{\eps}) \right)^{3/4}
\label{eq:pointwise-2}
\EEA
hence 
\BEA &&  e^{\tau\Del} e^{\half \lambda|\tilde{\del}_{\vec{\eps}}^j h_t|}(x) \le  \left( e^{(\tau+t)\Del}(e^{2\lambda|{\del}_{\vec{\eps}}^j h_0|})(x) \right)^{1/4}
\left( e^{(\tau+t)\Del}(e^{2\lambda |h_0|})(x)\right)^{1/4} \nonumber\\
&& \qquad\qquad\qquad  +
\left( e^{(\tau+t)\Del} ( e^{2\lambda |{\del}_{-\vec{\eps}}^j h_0|})(x+2^{j/2} \vec{\eps})
\right)^{1/4} \left( e^{(\tau+t)\Del} (e^{2\lambda |h_0|})(x+2^{j/2} \vec{\eps}) \right)^{3/4}
\EEA

and (letting $\del^j_{\vec{\eps},\vec{\eps}'} f(x):=\frac{f(x+2^{j/2}\vec{\eps})-f(x+2^{j/2}
\vec{\eps}')}{|\vec{\eps}-\vec{\eps}'|}$ and similarly
 $\tilde{\del}^j_{\vec{\eps},\vec{\eps}'}f(x):=\frac{f(x+2^{j/2}\vec{\eps})
 -(1-|\vec{\eps}-\vec{\eps}'|)f(x+2^{j/2}
\vec{\eps}')}{|\vec{\eps}-\vec{\eps}'|}$)
\BEA \sup_{\vec{\eps},\vec{\eps'}\in B(0,1)} e^{\tau\Del} e^{\frac{2}{5}\lambda |\del^j_{\vec{\eps},\vec{\eps}'} h_t|}(x) &\le & \sup_{\vec{\eps},\vec{\eps}'\in B(0,1)}
\left( e^{\tau\Del} e^{\frac{1}{2}\lambda |\tilde{\del}^j_{\vec{\eps},\vec{\eps}'} h_t|}(x)
\right)^{4/5}  \left( e^{\tau\Del} (e^{2\lambda |h_t|})(x+2^{j/2} \vec{\eps}') \right)^{1/5}
\nonumber\\
&\le & 2  \left[ \sup_{\vec{\eps},\vec{\eps}'\in B(0,1)} \left( e^{(\tau+t)\Del} ( e^{2\lambda{\del}_{\vec{\eps},\vec{\eps}'}^j h_0})(x+2^{j/2} \vec{\eps})
\right)^{1/5} \right] \left[ \sup_{\vec{\eps}\in B(0,1)}  \left( e^{(\tau+t)\Del} (e^{2\lambda |h_0|})(x+2^{j/2} \vec{\eps}) \right)^{4/5} \right] \nonumber\\
&\le & 2 e^{2\lambda ||| h_0|||_{{\cal W}_j,point}^{1,\infty;2\lambda}(x)},
\label{eq:3.78}  \EEA
where (compare with (\ref{eq:Wj1inftylambda}))
\BEQ  ||| h_0|||_{{\cal W}_{j,point}^{1,\infty;2\lambda}}(x) := \max\left(
\sup_{\vec{\eps}\in B(0,1)} |||h_0|||_{{\cal H}^{2\lambda}}(x+2^{j/2}\vec{\eps}),  \sup_{\vec{\eps},\vec{\eps'}\in B(0,1)} |||\del^j_{\vec{\eps},\vec{\eps'}}  h_0|||_{{\cal H}^{2\lambda}}(x)
 \right) \EEQ
is the aforementioned "pointwise quasi-norm". Combining Lemma \ref{lem:easy} with 
(\ref{eq:3.78}), we get a "pointwise" version of the "local" bounds of Lemma
\ref{lem:classical},
\BEQ |||h_t|||_{{\cal W}_{j,point}^{1,\infty;\frac{2}{5}\lambda}}(x)\le \frac{5}{2\lambda} \ln 2 \, +\, 5\, |||h_0|||_{{\cal W}_{j,point}^{1,\infty;2\lambda}}(x).\EEQ

Note that for $h_0$ small (in the appropriate pointwise ${\cal W}$-quasi norm),
one obviously expects $h_t$ to be small. Letting $a,b\ge 1$ be the two terms
appearing in (\ref{eq:pointwise-1}), one may bound $\max(a,b)$ $(a,b\ge 1)$ by
$ab$ instead of $a+b$. This way, we get rid of the unwanted additive factor $\frac{5}{2\lambda} \ln 2$, at the price of some more loss of regularity in the $\lambda$-exponents.

\bigskip

As a side application, let us  consider a rate $V=V(y)$ satisfying assumption (\ref{eq:even-stronger-assumption}), i.e. behaving
like $y^2$ for $y$ small or large, and show how to generalize
the conclusions of Proposition \ref{prop:decay-gradient} (iii).

\medskip

\begin{Corollary} Let $V$ satisfy assumption (\ref{eq:even-stronger-assumption}), $yV'(y)-V(y)\ge Cy^2$. Then
\BEQ |\nabla h_t(x)|\lesssim \left(\frac{|||h_0|||_{{\cal H}^{\lambda}}(x)/\lambda}{t}\right)^{1/2}. \EEQ
\end{Corollary}

{\bf Proof.} By (\ref{eq:pointwise-decay-gradient}),
\BEQ |\nabla h^{(L)}_t(x)|\lesssim \left(\frac{|h^{(L)}_t(x)|/\lambda}{t}\right)^{1/2}, \qquad x\in\R^d.\EEQ
By Lemma \ref{lem:easy}, $|h^{(L)}_t(x)|\le |||h_0^{(L)}|||_{{\cal H}^{\lambda}}(x).$ Hence, for every 
$\vec{\eps}\in B(0,1)$, 
$\frac{1}{|\vec{\eps}|} |h_t^{(L)}(x+\vec{\eps})-h_t^{(L}(x)| \lesssim \sup_{B(x,1)}
\left(\frac{|||h_0^{(L)}|||_{{\cal H}^{\lambda}}(\cdot)/\lambda}{t}\right)^{1/2}$. 
 The corollary
follows by letting first $L\to\infty$ and then $\vec{\eps}\to 0$. \hfill \eop


\section{Bounds for the infra-red cut-off inhomogeneous equation}


We introduce in this section the {\em scale $j$ infra-red cut-off} KPZ equation
(see eq. (\ref{intro:j}) in the Introduction, or (\ref{eq:KPZdowntoj}) below) and
prove the estimates for the solutions stated in Theorem 2 of the Introduction.
\S 4.1 is a somewhat lengthy motivation for eq. (\ref{intro:j}), in connection
to the general, motivating goal of showing diffusive large scale limit  for $d\ge 3$,
and to the multi-scale analysis of the linearized problem (Ornstein-Uhlenbeck's
equation) in section 5. The reader ill-at-ease with the scaling analysis may take
eq. (\ref{eq:KPZdowntoj}) for granted and jump directly to \S 4.2 and \S 4.3, where
we introduce new ${\cal W}$-spaces adapted to the time-dependent forcing term, $g$,
and then prove Lemma \ref{lem:bound-psin}, from which we deduce Theorem 2. Arguments
are generally strongly based on the computations of \S 3.3 and \S 3.4, together
with a {\em Trotter formula} sorting out the contribution of the
right-hand side.


\subsection{General philosophy of scale decompositions}


In this section, we  start our study of the inhomogeneous KPZ equation,
\BEQ \partial_t \psi(t,x)=\nu\Del\psi(t,x)+\lambda V(\nabla\psi(t,x))+g(t,x) \label{eq:full-KPZ} \EEQ
where $g(t,x)$ is a continuous forcing term. For the time being, we only consider an infra-red cut-off version of this
equation, see (\ref{eq:1.2}) or Definition \ref{def:KPZdowntoj} below. We  only require
here  good scale-dependent averaging properties for $g$ (see precise assumptions below). For the complete study (to be
developed in the further articles) we shall take for  $g$  a regularized white noise, denoted by $\eta$.

\medskip

The general motivation in the subsequent analysis is to exhibit an effective {\em scale separation} mechanism. 
In other words, let $G$ be the Green kernel,
\BEQ G:g\mapsto (Gg)(t):=\int_0^{+\infty} e^{\nu s\Del}g(t-s)ds \EEQ
(called {\em propagator} in the physics
literature). Eq. (\ref{eq:full-KPZ}) is equivalent to the integral equation,
\BEQ \psi=G(\lambda V(\nabla\psi)+g) \label{eq:integral-equation}.\EEQ
Now we want to write $G$ as a sum $G=\sum_{j\ge 0} G^j$ over {\em scales}, in such a way that
\begin{itemize}
\item[(1)] $G^j$ is "negligible" except at time-, resp. space distances of order 
$2^j$, resp. $2^{j/2}$;
\item[(2)] $\psi$ is well approximated by the sum $\sum_j \psi^{(j)}$, where $\psi^{(j)}$ is 
the solution  of the {\em single-scale} integral equation
\BEQ \psi^{(j)}=G^j(\lambda^{(j)} V(\nabla\psi^{(j)})+g^{(j)}), \label{eq:4.3} \EEQ 
where $g^{(j)}$ has typical fluctuations at time-, resp. space distances of 
order $2^j$, resp. $2^{j/2}$; and
\item[(3)] $\psi^{(j)}$ by the solution $\phi^{(j)}$ of the linearized equation, $\phi^{(j)}=G^j g^{(j)}$, at least for
$\lambda$ small enough or $j$ large enough.
\end{itemize}

The approximations in (1), (2), (3) are responsible for the renormalization procedure in which $\lambda$ becomes
the scale-dependent parameter $\lambda^{(j)}$ (actually $\lambda$ is not renormalized
in the case of the KPZ$_3$ model because it is super-renormalizable in the infra-red,
i.e. subcritical at large scales), $g$ becomes $g^{(j)}$, and $G^j$ also receives correction terms (see further article
in our series).

\medskip
At this point we are not interested in the renormalization procedure and would like in principle to consider a single-scale equation such
as (\ref{eq:4.3}),
\BEQ \psi^{(j)}=G^j(\lambda V(\nabla\psi^j)+g^{(j)}).\EEQ

The easiest way to select  fluctuations at time, resp. space distances of order $2^j$, resp. $2^{j/2}$ is to set
 \BEQ (G^j f)(t)=\int  ds \, \bar{\chi}^j(s) e^{s\nu\Del}f(t-s),\EEQ
 where $\bar{\chi}^j$ is a cut-off function s. t. $\bar{\chi}^j(s)=0$ if $s\ll 2^j$ or
 $s\gg 2^j$
(see Definition \ref{def:non-local-cut-off}). Coming back to $g=\eta$ to mimic the behaviour of the noisy KPZ
equation, we are led to set $\phi^j=G^j\eta$, $\eta^j=(\partial_t-\Del)\phi^j$. 
Recall
\BEQ d_{\phi}:=\half(\frac{d}{2}-1)\EEQ
is the  {\em scaling dimension} of the solution of the Ornstein-Uhlenbeck or of the KPZ equation, see Introduction.
 It is proved in Appendix A that 
\BEQ \esper[\phi^j(t,x)\phi^{j'}(t',x')]\lesssim 2^{-|j-j'|}(2^{-\max(j,j')})^{2d_{\phi}} e^{-c2^{-\max(j,j')}|t-t'|-c2^{-\max(j,j')/2}|x-x'|} \label{eq:Cphi} \EEQ
\BEQ \esper[\eta^j(t,x)\eta^{j'}(t',x')]\lesssim 2^{-2|j-j'|}(2^{-\max(j,j')})^{2+2d_{\phi}} e^{-c2^{-\max(j,j')}|t-t'|-c2^{-\max(j,j')/2}|x-x'|}
\label{eq:Ceta}\EEQ
for some constant $c>0$.  Consider first the diagonal covariance ($j=j'$): since $\phi^j$ and $\eta^j$  are 
Gaussian, (\ref{eq:Cphi}), (\ref{eq:Ceta}) essentially mean that the following {\em scalings} hold,
\BEQ \phi^j(t,x)=O(2^{-jd_{\phi}}), \qquad \eta^j(t,x)=O((2^{-j})^{1+d_{\phi}}), \label{eq:scaling}\EEQ
with random prefactors. The bounds in Appendix A also yield an order of magnitude of the gradients, with a
supplementary $2^{-j/2}$ factor,
\BEQ \nabla\phi^j(t,x)=O((2^{-j})^{\half+d_{\phi}}), \qquad \nabla\eta^j(t,x)=O((2^{-j})^{\frac{3}{2}+d_{\phi}})
\label{eq:scaling-gradient}.\EEQ
  For $j\not= j'$ one has
an extra decaying exponential factor in $2^{-|j-j'|}$ which lies at the root of the scale separation mechanism. 

\medskip

We shall not pursue along this road in this article. The reason
is that the integral equation (\ref{eq:4.3}) is a delay, non-local equation which does not satisfy at all the
maximum principle, and we have no a priori bounds for its solutions, save in the 
{\em perturbative} regime where $g$ or $\eta$ is small. So we  introduce instead in the sequel
 a very simple {\em infra-red cut-off of scale $j$} for the propagator, namely, we replace $\nu\Del$ by $\nu\Del-2^{-j}$.  Denoting by $G^{j\to}$ the Green kernel of the operator
 $\nu\Del-2^{-j}$, one has the explicit formula
\BEQ G^{j\to}(t,x;t',x')={\bf 1}_{t>t'} e^{-2^{-j}(t-t')} p_{\nu(t-t')}(x-x'),
\label{eq:Gjfleche} \EEQ
which makes apparent an exponential decay in time and space: since $\inf_{s>0} (\frac{|x-x'|^2}{2\nu s}+s2^{-j})\approx
2^{-j/2} |x-x'|$, 
\BEQ G^{j\to}(t,x;t',x')\lesssim (t-t')^{-d/2} e^{-c2^{-j}(t-t')-c2^{-j/2}|x-x'|}. \label{eq:exp-decay} \EEQ
for some constant $c>0$.  The idea is that $G^{j\to}$ is a good substitute for the sum $\sum_{k\le j} G^k$. We also
replace the force term $g$ by $g^j$ such that $g^j(t,x)=O((2^{-j})^{1+d_{\phi}})$ as for $\eta^j$, see (\ref{eq:scaling}). 
Thus the new equation is the following.

\begin{Definition}[inhomogeneous KPZ equation with scale $j$ infra-red cut-off] \label{def:KPZdowntoj}
The inhomogeneous KPZ equation with scale $j$ infra-red cut-off is
\BEQ \partial_t\psi=(\Del-2^{-j})\psi+\lambda V(\nabla\psi)+g. \label{eq:KPZdowntoj}\EEQ 
\end{Definition}

As in \S 3.4, we have chosen $\nu=1$ for simplicity. The integral form of this equation is
\BEQ \psi=G^{j\to}(V(\nabla\psi)+g).\EEQ
Note that the kernel $G^{j\to}$ has no {\em ultra-violet cut-off}, in the sense that it behaves like
the full Green kernel  $G$ for time separations $|t-t'|\ll 2^j$. Because $g$
{\em has} an ultra-violet cut-off, 
 it actually turns out that 
the solution $\psi$ of (\ref{eq:KPZdowntoj})
has the correct scaling, $\psi(t,x)=O(2^{-jd_{\phi}})$, see (\ref{eq:scaling}), under appropriate assumptions
on $g$ that 
we now proceed to write down. Note that, conversely, since $G^{j\to}$ has an {\em infra-red}
cut-off, it is not really necessary to put an infra-red cut-off on $g$ too
(see remark at the very end of section 5).

\subsection{Functional spaces of scale $j$}


As in the case of the homogeneous equation, we need a "local supremum" operation
adapted to space-time functions $g$. Generalizing (\ref{eq:locsup}) in a straightforward
way, taking into account the parabolic scaling, we let
\BEQ \Locsup^j g(t,x):=\sup_{s\in (t-2^j,t+2^j)} \sup_{y\in B(x,2^{j/2})} |g(s,y)|.\EEQ

We shall assume that the right-hand side, $g$, sits in a new convex subspace ${\cal W}_j^{1,\infty;\lambda}([0,T])\subset
C([0,T], W^{1,\infty}_{loc}(\R^d))$ that we now proceed to define, in its stronger "local quasi-norm" version
(a weaker, somewhat ugly "pointwise quasi-norm" version also exists),

\newpage

\begin{Definition} \label{assumption:etaj}
For $g\in C([0,T], W^{1,\infty}_{loc}(\R^d))$, let 
\begin{itemize}
\item[(i)] for all $x\in\R^d$,
\BEQ |||g|||_{\lambda,j}([0,T],x):= 2^{-j} \int_0^T e^{-2^{-j}s}
|||\, 2^j g(T-s,\cdot)\, |||_{{\cal H}^{\lambda}}(x)\, ds; \label{eq:etaj1}\EEQ
\item[(ii)] for all $x\in\R^d$,
\BEQ |||g|||_{{\cal W}_j^{1,\infty;\lambda}([0,T])}(x):= \max\left(\, |||\, \Locsup^j g |||_{\lambda,j}([0,T],x), \,  |||\, 2^{j/2} \Locsup^j |\nabla g| \, |||_{\lambda,j}([0,T],x) \right).  \label{eq:etaj2}\EEQ
\end{itemize}
If $|||g|||_{{\cal W}_j^{1,\infty;\lambda}([0,T])}(x)<\infty$ for all $x\in \R^d$, 
then we say that $g\in {\cal W}_j^{1,\infty;\lambda}([0,T])$.
\end{Definition}

\medskip

If $|||g|||_{{\cal W}^{1,\infty;2\lambda}_j([0,t])}(x)=O(2^{-jd_{\phi}})$ then
Theorem 2 in the Introduction (proved in the following subsection) ensures that
$\psi(t,x)=O(2^{-jd_{\phi}})$, $\nabla\psi(t,x)=O((2^{-j})^{\half+d_{\phi}})$ as
expected (see (\ref{eq:scaling-gradient})). It is proved in section 6 that, indeed,
$|||\eta^j||||_{{\cal W}^{1,\infty;2\lambda}_j([0,t])}(x)=O(2^{-jd_{\phi}})$ a.s.


\subsection{Bounds}


 Consider an
initial condition $\psi_0\in{\cal W}^{1,\infty;2\lambda'}_{j}\cap  C^2$ with
$\lambda'>\lambda$, and forcing term $g\in W_j^{1,\infty;2\lambda}([0,T])\cap C([0,T],
C^3(\R^d))$, for
some large but finite time horizon $T$. We prove here our second main theorem, Theorem 2
in the Introduction.

\medskip

 We use the
following notations in this paragraph. The homogeneous nonlinear semi-group generated by the homogeneous
KPZ equation (\ref{eq:2.1}) is denoted by $\Phi^{\lambda}(t)$, i.e.
$\Phi^{\lambda}(t)h_0$ is the solution at time $t$ of the homogeneous KPZ equation with initial condition
$h_0\in {\cal W}^{1,\infty}$. Let also
$\tau_k(s):C(\R^d)\to C(\R^d), f\mapsto \tau_k(s)f$ by $\tau_k(s)f(x):=\int_{kt/n}^{s+kt/n} g(u,x) du+
f(x)$ ($0\le k\le n-1$). Treating each term in (\ref{eq:KPZdowntoj}) separately, 
we get three equations: (i) $\partial_t\psi=-2^{-j}\psi$, with solution $\psi(t)=e^{-2^{-j}t}\psi(0)\simeq (1-2^{-j}t)\psi(0)$ for small $t$; (ii) $(\partial_t-\Del)\psi=\lambda V(\nabla\psi)$, with solution $\psi(t)=\Phi^{\lambda}(t)\psi(0)$; (iii) $\partial_t\psi=g$, with solution $\psi(s+kt/n)=\tau_k(s)\psi_{kt/n}$. Alternating the action of these three
non-linear semi-groups, we obtain

\begin{Definition} \label{def:Trotter}
 
  Let, for $k=0,\ldots,n$,
\BEQ \psi^{(n)}_{kt/n}(x):=\left( (1-2^{-j}t/n) \Phi^{\lambda}(t/n) \tau_{k-1}(t/n) \right)
\left( (1-2^{-j}t/n) \Phi^{\lambda}(t/n) \tau_{k-2}(t/n)\right)\cdots
\left( (1-2^{-j}t/n) \Phi^{\lambda}(t/n) \tau_{0}(t/n) \right)
 \psi_0(x).\EEQ
\end{Definition}

Having a "Trotter formula" in this setting means proving that $\psi^{(n)}$ converges
in some norm to $\psi$, solution of (\ref{eq:KPZdowntoj}).
Trotter formulas have been shown with some generality for non-linear monotonous operators acting on Hilbert spaces
\cite{Bre}.
However here the natural spaces, $L^{\infty}$, ${\cal W}^{1,\infty}$ and their
localized counterparts, ${\cal H}^{\lambda}$, ${\cal W}_j^{1,\infty;\lambda}$, are not  Hilbert spaces. To show this lemma we therefore follow 
instead the proof of convergence of
"viscous splitting" algorithms for the Navier-Stokes equation, as found in \cite{BerMaj}, \S 3.4, resting
on their {\em stability} and {\em consistency}. {\em Stability} means that 
the sequence $(\psi^{(n)})_n$ is bounded in the relevant norms. Once one has
proved stability, one may prove {\em consistency}, i.e. prove that $\psi^{(n)}-\psi$ converges
to $0$ when $n\to\infty$.

\begin{Lemma}[stability] \label{lem:bound-psin}
Let $n>2^{-j}t$ and $p:=(2^{-j}\frac{t}{n})^{-1}-1$.
Assume $\psi_0\in {\cal W}_j^{1,\infty;2\lambda\frac{p+1}{p}}$ (i.e.
$\lambda'\ge \lambda\frac{p+1}{p}$). Then the following bounds hold,
\BEQ  |||\, \locsup^j \psi^{(n)}_t \, |||_{{\cal H}^{a\lambda}}(x) \le (1+O(\frac{1}{p}))
e^{-2^{-j}t} |||\locsup^j \psi_0|||_{{\cal H}^{a\lambda\frac{p+1}{p}}}(x) +  |||\, \locsup^j g|||_{a\lambda,j}([0,t],x), \ \ a\in [1,2]  \label{eq:4.20}\EEQ
and 
\BEQ |||\, 2^{j/2} \locsup^j |\nabla\psi^{(n)}_t|\, |||_{{\cal H}^{2\lambda/5}}(x) \le 
5(1+O(\frac{1}{p})) \left( |||g|||_{{\cal W}_j^{1,\infty;2\lambda}([0,t])}(x) + e^{-2^{-j}t} |||\psi_0|||_{{\cal W}_j^{1,\infty;2\lambda\frac{p+1}{p}}}(x) \right). \label{eq:4.21} \EEQ
\end{Lemma}

{\bf Proof.} Note first that the condition $\lambda'\ge \lambda\frac{p+1}{p}$ is always
verified for $n$ large enough since by hypothesis $\lambda'>\lambda$. 
We shall rely on the following two elementary bounds,
\BEQ (e^{\lambda a \locsup^j |\Phi^{\lambda}(s) f|})^*(x)\le (e^{\lambda a\locsup^j  |f|})^*(x),\qquad (a\ge 1) \label{eq:lambdaaPhis} \EEQ
(see (\ref{eq:3.66})) and H\"older's inequality
\BEQ (e^{\lambda |f +\tilde{f}|})^*(x)\le \left[ (e^{\lambda \frac{p+1}{p} |\tilde{f}|})^*(x)\right]^{\frac{p}{p+1}}
\left[ (e^{\lambda (p+1)|f|})^*(x)\right]^{\frac{1}{p+1}}.
 \label{eq:Holder}\EEQ

\smallskip

\noindent Choose $p=(2^{-j}\frac{t}{n})^{-1}-1$ in (\ref{eq:Holder}). \\ For $0\le x\le \ln(2)$, 
$e^x-1\le x+x^2\le \frac{x}{1-x}$. Hence (letting $x=2^{-j}\frac{t}{n}$)
$$e^{2^{-j}\frac{t}{n}}\le \frac{p+1}{p}= \frac{1}{1-2^{-j}t/n} \to_{n\to\infty} 1, \qquad   
\left(\frac{p}{p+1}\right)^n \to_{n\to\infty} e^{-2^{-j}t}.$$
  Thus, by (\ref{eq:lambdaaPhis}), for $a\ge 1$, 
\BEA  \left(e^{a\lambda \frac{p+1}{p} \locsup^j \psi^{(n)}_{(k+1)t/n}(\cdot)}\right)^*(x) &=& \left( e^{a\lambda\,  \locsup^j(\Phi^{\lambda}(t/n)\circ
\tau_k(\frac{t}{n}))(\psi^{(n)}_{kt/n}(\cdot))} \right)^*(x) \nonumber\\
&\le&  \left[  (e^{a\lambda \frac{p+1}{p} \locsup^j \psi^{(n)}_{kt/n}(\cdot)})^*(x)|
\right]^{\frac{p}{p+1}} \left[ (e^{a\lambda (p+1)\int_{kt/n}^{(k+1)t/n} \locsup^j g(u,\cdot) \, du })^*(x)\right]^{\frac{1}{p+1}} \nonumber\\
&\le &  \left[  (e^{a\lambda \frac{p+1}{p} \locsup^j \psi^{(n)}_{kt/n}(\cdot)})^*(x)|
\right]^{\frac{p}{p+1}} \left[ (e^{a\lambda 2^j \Locsup^j g(kt/n,\cdot) \, du })^*(x)\right]^{\frac{1}{p+1}}.\EEA

By induction on $k$, this gives
\BEQ  \left(e^{a\lambda\, \locsup^j \psi^{(n)}_t}\right)^*(x) \le  \left(e^{a\lambda \frac{p+1}{p} \locsup^j \psi^{(n)}_{t}}\right)^*(x)\le   \left[\prod_{k=0}^{n-1} (A^{(n)}_k(x))^{\frac{1}{p+1} (\frac{p}{p+1})^k}  \right]
(A^{(n)}_n(x))^{(\frac{p}{p+1})^n},  \label{eq:AAA} \EEQ
where
\BEQ A^{(n)}_k(x)= \left(e^{a\lambda 2^j \Locsup^j g(t-kt/n,\cdot)  }
\right)^*(x)\qquad (k=0,\ldots,n-1),\qquad A^{(n)}_n(x)=(e^{a\lambda\frac{p+1}{p}\locsup^j \psi_0})^*(x) \EEQ

Hence 
\BEA  \frac{1}{a\lambda}\ln \left(e^{a\lambda\,  \locsup^j \psi^{(n)}_t}\right)^*(x)
&\le & (\frac{p}{p+1})^{n-1} |||\psi_0|||_{{\cal H}^{a \frac{p+1}{p} \lambda}}(x) +
\sum_{k=0}^{n-1} \frac{1}{p+1} (\frac{p}{p+1})^k |||2^j \, \Locsup^j g(t-\frac{kt}{n},\cdot)|||_{{\cal H}^{a\lambda}}(x) \nonumber\\
&\le& e^{-2^{-j} \frac{n-1}{n} t} |||\psi_0|||_{{\cal H}^{a\frac{p+1}{p}\lambda}}(x) + 
2^{-j} \frac{t}{n} \sum_{k=0}^{n-1} e^{-2^{-j} \frac{k}{n} t}  |||2^j \, \Locsup^j g(t-\frac{kt}{n},\cdot)|||_{{\cal H}^{a\lambda}}(x) \nonumber\\
&\le &  (1+O(2^{-j}\frac{t}{n})) \left( e^{-2^{-j}t} |||\psi_0|||_{{\cal H}^{a \frac{p+1}{p}\lambda }}(x) +  |||\, \Locsup^j g|||_{a\lambda,j}[0,t],x) \right),
\label{eq:4.28}  \EEA 

as claimed in (\ref{eq:4.20}).

Note for further use that the exponents in (\ref{eq:AAA}) sum up to 1,
\BEQ (\frac{p}{p+1})^n+ \sum_{k=0}^{n-1} \frac{1}{p+1} (\frac{p}{p+1})^k =1.
\label{eq:sum=1} \EEQ

\bigskip
\medskip
The proof of (\ref{eq:4.21}) is similar but requires a further elaboration on the arguments developed in the course
of the proof of Lemma \ref{lem:classical}, to which we refer the reader for the notations. 
Let $\vec{\eps}\in B(0,1)$ and $a\ge 1$.
First

\BEA   e^{\lambda \frac{p+1}{p} \tilde{\del}_{\vec{\eps}}^j \psi^{(n)}_{(k+1)t/n}(x) } &=&
e^{\lambda \tilde{\del}^j_{\vec{\eps}} \Phi^{\lambda}(\frac{t}{n})(\psi^{(n)}_{kt/n}
+ \int_{kt/n}^{(k+1)t/n} g_u\, du)(x)} \nonumber\\
& \le& e^{\frac{t}{n}\Del} \left( e^{\lambda \tilde{\del}^j_{\vec{\eps}}  (\psi^{(n)}_{kt/n}
+ \int_{kt/n}^{(k+1)t/n} g_u\, du)(x)}  \right) \nonumber\\
& \le& \left[ e^{\frac{t}{n}\Del} (e^{\lambda \frac{p+1}{p} \tilde{\del}_{\vec{\eps}}^j
\psi_{kt/n}^{(n)}} )(x) \right]^{\frac{p}{p+1}} \left[ e^{\frac{t}{n}\Del} ( e^{\lambda
(p+1) \int_ {kt/n}^{(k+1)t/n} \tilde{\del}_{\vec{\eps}}^j g_u\, du} )  (x) \right]^{\frac{1}{p+1}}.
\EEA

By induction on $k$, this yields
\BEA e^{\lambda \frac{p+1}{p} \tilde{\del}_{\vec{\eps}}^j \psi^{(n)}_{(k+1)t/n}(x) }
&\le& \left[ e^{t\Del} \left( e^{\lambda\frac{p+1}{p} \tilde{\del}_{\vec{\eps}}^j
\psi_0} \right)(x) \right]^{(\frac{p}{p+1})^n} \prod_{k=0}^{n-1}
\left[ e^{\frac{kt}{n}\Del} \left( e^{\lambda(p+1) \int_{t-(k+1)t/n}^{t-kt/n}
\tilde{\del}_{\vec{\eps}}^j g_u\, du} \right)(x) \right]^{\frac{1}{p+1} (\frac{p}{p+1})^k} \nonumber\\
& \le & \left[\prod_{k=0}^{n-1} (B^{(n)}_k(x))^{\frac{1}{p+1} (\frac{p}{p+1})^k}  \right]
(B^{(n)}_n(x))^{(\frac{p}{p+1})^n}, \EEA
where
\BEQ B^{(n)}_k(x)= e^{\frac{kt}{n}\Del} \left( e^{\lambda 2^j \left[ 2^{j/2} \Locsup^j |\nabla g|(t-kt/n,\cdot) +
\Locsup^j g(t-kt/n,\cdot)\right] }\right)(x)
\qquad (k=0,\ldots,n-1),\EEQ
\BEQ \qquad \qquad  B^{(n)}_n(x)=e^{t\Del} \left(e^{\lambda\frac{p+1}{p} \left[2^{j/2} \locsup^j |\nabla\psi_0| + \locsup^j \psi_0\right] } \right)(x) \EEQ

For the reverse inequality, proceeding as in (\ref{eq:3.72}), we get
\BEQ e^{-\lambda \frac{p+1}{p} \tilde{\del}_{\vec{\eps}}^j \psi_t^{(n)}(x)} =
e^{\lambda \frac{p+1}{p} (1-|\vec{\eps}|) \tilde{\del}_{-\vec{\eps}}^j \psi_t^{(n)}(x+2^{j/2} \vec{\eps})} e^{-\lambda \frac{p+1}{p} (2-|\vec{\eps}|) \psi_t^{(n)}(x+2^{j/2}
\vec{\eps})},\EEQ
whence the two-sided, uniform inequality,
\BEQ e^{\lambda \frac{p+1}{p} \sup_{\vec{\eps},\vec{\eps}'\in B(0,1)}
|\tilde{\del}_{\vec{\eps}-\vec{\eps}'}^j \psi_t^{(n)}(x+2^{j/2}\vec{\eps}')|} \le 
 \left[\prod_{k=0}^{n-1} (B^{(n)}_k(x))^{\frac{1}{p+1} (\frac{p}{p+1})^k}  \right]
(B^{(n)}_n(x))^{(\frac{p}{p+1})^n}  e^{2\lambda \frac{p+1}{p} \locsup^j
\psi_t^{(n)}(x)} ,\EEQ

from which (using H\"older's inequality with exponents (\ref{eq:sum=1}))
\BEA &&
\left(e^{ \half \lambda  \sup_{\vec{\eps},\vec{\eps}'\in B(0,1)}
|\tilde{\del}_{\vec{\eps}-\vec{\eps}'}^j \psi_t^{(n)}(x+2^{j/2}\vec{\eps}')|} \right)^*(x)
\nonumber\\
&& \qquad \le  \left(e^{ \half \lambda \frac{p+1}{p} \sup_{\vec{\eps},\vec{\eps}'\in B(0,1)}
|\tilde{\del}_{\vec{\eps}-\vec{\eps}'}^j \psi_t^{(n)}(x+2^{j/2}\vec{\eps}')|} \right)^*(x)
\nonumber\\
&& \qquad \le \left[ \left( (B_n^{(n)}(x))^{(\frac{p}{p+1})^n} \prod_{k=0}^{n-1}
(B_k^{(n)}(x))^{\frac{1}{p+1} (\frac{p}{p+1})^k}  \right)^*(x) \right]^{1/2} 
\left[ \left(e^{2\lambda \frac{p+1}{p} \locsup^j \psi_t^{(n)}(\cdot)} \right)^*(x)\right]^{1/2} \nonumber\\
&& \qquad \le \left\{ (B^{(n)}_n)^*(x))^{(\frac{p}{p+1})^n} \prod_{k=0}^{n-1} ((B^{(n)}_k)^*(x))^{\frac{1}{p+1} (\frac{p}{p+1})^k} \right\}^{1/2} \left[ \left(e^{2\lambda \frac{p+1}{p} \locsup^j \psi_t^{(n)}(\cdot)} \right)^*(x)\right]^{1/2}
 \nonumber\\
&& \le \left\{  \left[ \left( e^{2\lambda\, \cdot\, 2^{j/2} \locsup^j |\nabla \psi_0(\cdot)|}
\right)^*(x) \right] ^{(\frac{p}{p+1})^n} \left[ \left(e^{2\lambda\, \cdot\,  \locsup^j \psi_0(\cdot)|} \right)^*(x)\right]^{(\frac{p}{p+1})^n}  \right. \nonumber\\
&& \qquad \left. \prod_{k=0}^{n-1} 
 \left( \left( e^{2\lambda\, \cdot\, 2^{3j/2} \Locsup^j |\nabla g(t-kt/n,\cdot)|}
\right)^*(x) \left(e^{2\lambda\, \cdot\, 2^j \Locsup^j  g(t-kt/n,\cdot)|} \right)^*(x)
 \right)^{\frac{1}{p+1} (\frac{p}{p+1})^k}  \right\}^{1/4} \nonumber\\
 && \qquad \qquad \qquad \qquad \left[ \left(e^{2\lambda \frac{p+1}{p} \locsup^j \psi_t^{(n)}(\cdot)} \right)^*(x)\right]^{1/2} \nonumber\\
 && \qquad \le e^{\lambda (1+O(2^{-j}t/n)) |||g|||_{{\cal W}_j^{1,\infty;2\lambda}([0,t])}(x)}  
 e^{\lambda e^{-2^{-j} t} |||\psi_0|||_{{\cal W}^{1,\infty;2\lambda\frac{p+1}{p}}_j}(x)} \left[ \left(e^{2\lambda \frac{p+1}{p} \locsup^j \psi_t^{(n)}(\cdot)} \right)^*(x)\right]^{1/2}
\EEA
(compare with (\ref{eq:4.28})).
Finally, using H\"older's inequality with conjugate exponents $(p,q)=(\frac{5}{4},5)$, 
\BEA && \left(e^{\frac{2}{5} \lambda 2^{j/2} \locsup^j |\nabla \psi_t^{(n)}|} \right)^*(x) \le 
\left( e^{\frac{2}{5} \lambda \left[ \sup_{\vec{\eps},\vec{\eps}'\in B(0,1)}
|\tilde{\del}_{\vec{\eps}-\vec{\eps}'}^j \psi_t^{(n)}(x+2^{j/2}\vec{\eps}')| + 
\locsup^j \psi_t^{(n)}(\cdot) \right]} \right)^*(x) \nonumber\\
&& \le \left[\left(e^{ \half \lambda  \sup_{\vec{\eps},\vec{\eps}'\in B(0,1)}
|\tilde{\del}_{\vec{\eps}-\vec{\eps}'}^j \psi_t^{(n)}(x+2^{j/2}\vec{\eps}')|} \right)^*(x)
\right)^{4/5}
\left[ \left(e^{2\lambda\, \locsup^j \psi_t^{(n)}(\cdot)} \right)^*(x) \right]^{1/5}
 \nonumber\\
&& \le  e^{\frac{4}{5} \lambda (1+O(2^{-j}t/n)) |||g|||_{{\cal W}^{1,\infty;2\lambda}_j([0,t])}(x)} 
 e^{\frac{4}{5} \lambda e^{-2^{-j} t}  |||\psi_0|||_{{\cal W}^{1,\infty;2\lambda\frac{p+1}{p}}_j}(x)}  
\left[ \left(e^{2\lambda \frac{p+1}{p} \locsup^j \psi_t^{(n)}(\cdot)} \right)^*(x)\right]^{3/5} 
 \EEA

Applying now our previous bound (\ref{eq:4.28}) yields (\ref{eq:4.21}).

\hfill \eop

\bigskip

Let us now turn to the proof of {\em consistency}. Since this is an essentially perturbative,
short-time argument, it introduces non-linear terms, typically, $|\nabla\psi^{(n)}|^2$, whose ${\cal H}^{\lambda}$-norm cannot be assumed to be bounded. Hence
we use the same cut-off procedure as in \S 3.4, and introduce instead the doubly-indexed sequence $(\psi^{(L,n)})_{L,n}$, $L,n\in\N$, constructed
as in Definition \ref{def:Trotter} but with cut-off initial data and right-hand side,
$\psi_0 \rightsquigarrow \psi^{(L)}_0(\cdot):=\psi_0(\cdot) \chi^{(L)}(\cdot)$, 
$g(s,\cdot)\rightsquigarrow g^{(L)}(s,\cdot):=g(s,\cdot) \chi^{(L)}(\cdot)$, where
$\chi^{(L)}$ is a cut-off function as in Lemma \ref{lem:h(n)}.  Since $g^{(L)}$ is 
regular and bounded, the standard theory of existence for KPZ equation implies that
$\psi^{(L)}$ is classical. For sake of convenience, we slightly
modify the notation of Definition \ref{def:Trotter} (but not the scheme of approximation) by letting 

\BEQ \psi^{(L,n)}_{(k+1)t/n}(x):=\tau_k(t/n) (1-2^{-j}t/n) \Phi^{\lambda}(t/n) \psi_{kt/n}^{(L,n)}
\qquad (k\ge 0), \psi^{(L,n)}_0(x):=\tau_0(t/n)\psi^{(L)}_0(x).\EEQ

Introduce, for $0\le s\le t/n$,
\BEA  \psi^{(L,n)}_{s+kt/n}(x) &:=& \tau_k(s) (1-2^{-j}s) \Phi^{\lambda}(s) \psi_{kt/n}^{(L,n)}
\nonumber\\
&=& \int_{kt/n}^{s+kt/n} g^{(L)}(u,x)\, du + (1-2^{-j}s)\Phi^{\lambda}(s) \psi^{(L,n)}_{kt/n}(x).\EEA

Then
\BEA \partial_s \psi^{(L,n)}_{s+kt/n}(x)&=& g^{(L)}_{s+kt/n}(x)-2^{-j} \Phi^{\lambda}(s)
\psi^{(L,n)}_{kt/n}(x) + (1-2^{-j}s) (\Del (\Phi^{\lambda}(s) \psi^{(L,n)}_{kt/n}(x)) +
\lambda V(\nabla(\Phi^{\lambda}(s) \psi^{(L,n)}_{kt/n}(x))) \nonumber\\
&=& (\Del-2^{-j})\psi_{s+kt/n}^{(L,n)}(x)+\lambda V(\nabla\psi_{s+kt/n}^{(L,n)}(x))+g^{(L)}_{s+kt/n}(x) + \left( A_1^{(L,n)}+A_2^{(L,n)}+A_3^{(L,n)} \right)(s,x), \nonumber\\ \EEA
where:
\BEQ |A_1^{(L,n)}(s,x)|=\left|-2^{-j} ((1-2^{-j}s)^{-1}-1)\psi_{s+kt/n}^{(L,n)}(x) \right|\lesssim
2^{-2j} \frac{t}{n} |\psi^{(L,n)}_{s+kt/n}(x)|; \EEQ
\BEQ |A_2^{(L,n)}(s,x)|=\left| (1-2^{-j}s) \lambda V((1-2^{-j}s)^{-1}\nabla \psi_{s+kt/n}^{(L,n)}(x))
-\lambda V(\nabla \psi^{(L,n)}_{s+kt/n}(x)) \right|\lesssim 2^{-j} \frac{t}{n} \lambda
|\nabla \psi_{s+kt/n}^{(L,n)}(x)|^2;\EEQ
\BEA && |A_3^{(L,n)}(s,x)| = \left| (1-2^{-j}s) \lambda \left( V((1-2^{-j}s)^{-1} \nabla(\psi_{s+kt/n}^{(L,n)}(x)-\int_{kt/n}^{s+kt/n} g^{(L)}(u,x)du) - V((1-2^{-j}s)^{-1}
\nabla \psi^{(L,n)}_{s+kt/n}(x)) \right) \right.  \nonumber\\
&& \qquad\qquad\qquad \qquad \left. + 2^{-j} (1-2^{-j}s)^{-1} \int_{kt/n}^{s+kt/n} g^{(L)}(u,x)du -
\int_{kt/n}^{s+kt/n} \Del g^{(L)}(u,x)du \right| \nonumber\\
&& \qquad \qquad \lesssim   \frac{t}{n} \left\{ 2^{-j} \Locsup^j g^{(L)}(kt/n,x) + \lambda  \left(|\nabla \psi_{s+kt/n}^{(L,n)}(x)|^2 + (\Locsup^j |\nabla g^{(L)}_{kt/n}|(x))^2 \right) + \Locsup^j \Del g^{(L)}_{kt/n}(x) \right\}.\nonumber\\
\EEA

Let now $\Psi^{(L,n)}_u(x):=\psi^{(L,n)}_u(x)-\psi^{(L)}_u(x)$, $\frac{kt}{n}<u:=s+\frac{kt}{n}<(k+1)\frac{t}{n} $.
Subtracting the evolution equations for $\psi^{(L,n)}$ and $\psi^{(L)}$, one obtains
\BEA  \partial_u \Psi^{(L,n)}_u(x) &=&  (\Del-2^{-j})\Psi^{(L,n)}_u(x)+ \lambda \left(V(\nabla
\psi^{(L,n)}(x))-V(\nabla\psi^{(L)}(x)) \right) + (A_1^{(L,n)}+A_2^{(L,n)}+A_3^{(L,n)})(s,x)
\nonumber\\  &=&
(\Del-2^{-j})\Psi^{(L,n)}_u(x)+  a(t,x)\cdot\nabla\Psi^{(L,n)}_u(x) + (A_1^{(L,n)}+A_2^{(L,n)}+A_3^{(L,n)})(s,x),  \label{eq:subtracting-evolution} \EEA
where (as follows from Lemma \ref{lem:bound-psin} and standard bounds for $\nabla\psi^{(L)}$) $|a(t,x)|\le C$, with
 \\ $C=C(||g^{(L)}||_{\infty}, ||\nabla g^{(L)}||_{\infty},
 ||\psi_0^{(L)}||_{\infty}, ||\nabla\psi_0^{(L)}||_{\infty})$.
By the usual comparison principle,
\BEQ ||\Psi^{(L,n)}_{(k+1)t/n}||_{\infty} \le ||\Psi^{(L,n)}_{kt/n}||_{\infty} + 
O( (\frac{t}{n})^2) \left( 2^{-j} ||g^{(L)}||_{\infty} + ||\Del g^{(L)}||_{\infty} + \lambda  \left(||\nabla g^{(L)}||_{\infty}^2 + ||\nabla \psi^{(L,n)}||_{\infty}^2 \right) \right),\EEQ
from which by induction $||\psi_t^{(L,n)}-\psi^{(L)}_t||_{\infty} \le C \frac{t^2}{n}+
||\psi_0^{(L,n)}-\psi_0^{(L)}||_{\infty} \lesssim C(\frac{t^2}{n}+\frac{t}{n})$, with\\  $C
=C(||g^{(L)}||_{\infty}, ||\nabla g^{(L)}||_{\infty}, ||\nabla^2 g^{(L)}||_{\infty},
 ||\psi_0^{(L)}||_{\infty}, ||\nabla\psi_0^{(L)}||_{\infty})$.
Hence $\psi^{(L,n)}\to\psi^{(L)}$ locally uniformly. Differentiating (\ref{eq:subtracting-evolution})
one  prove similarly that $\nabla\psi^{(L,n)}\to\nabla\psi^{(L)}$ locally uniformly
(at this point we need $g(t,\cdot)$ to be $C^3$). Thus the bounds of Lemma
\ref{lem:bound-psin} hold for the limit $\psi^{(L)}$.

\medskip\noindent Finally Lemma \ref{lem:h(n)} allows
to conclude that $\psi^{(L)}\to \psi$ locally uniformly, with $\psi$ solution of (\ref{eq:KPZdowntoj}), and $\nabla\psi^{(L)}\to\nabla\psi$ locally uniformly, with
the limit, $\psi$, satisfying the same bounds as in Lemma \ref{lem:bound-psin}.

\medskip \noindent On the
other hand, in absence of a comparison principle for the inhomogeneous KPZ equation,
we {\em cannot} conclude to the unicity of the limit.  The difficulty here is
to control the dependence of $\psi(t,x)$ on the (possibly large!) values of the data $(\psi_0,g)$ at
space locations $y$ at distance $|x-y|\to\infty$. We {\em did not manage}, by purely PDE
arguments, to show that the sequence $(\psi^{(L)})$ is Cauchy for the uniform convergence
on compacts. At this point it is more natural to solve the KPZ equation by using 
characteristics. In \cite{Unt-KPZ2} (see \S 2.2) it is shown that characteristics going
far astray from their starting point $x$ hardly contribute to the value of $h_t$ at $x$,
implying, with more generality than required here, that $(\psi^{(n)})$ is a Cauchy sequence whenever $\psi^{(n)}$ are the solutions of the KPZ equations $(\partial_t-\Del+2^{-j})\psi^{(n)}=\lambda V(\nabla \psi^{(n)})+g^{(n)}$ with initial condition $\psi_0^{(n)}$, for
all sequences of bounded data $\psi_0^{(n)}\in {\cal W}^{1,\infty}$, $g^{(n)}\in C([0,T],{\cal W}^{1,\infty})$ such that 
\begin{itemize}
\item[(i)] $|||\psi_0^{(n)}|||_{{\cal W}^{1,\infty}_j}(x)$,
$|||g^{(n)}|||_{{\cal W}^{1,\infty}_j([0,t])}(x)$ are uniformly bounded; 
\item[(ii)] for all
$K\subset\R^d$ compact,
$\psi_0^{(n)}\to_{n\to\infty} \psi_0$ in ${\cal W}^{1,\infty}(K)$ and $g^{(n)}\to_{n\to\infty} g$ in $C([0,t],{\cal W}^{1,\infty})$.
\end{itemize}

This concludes the proof of Theorem 2 in the Introduction.


\section{Scale decompositions}


As a general motivation for this section, consider the Ornstein-Uhlenbeck process (\ref{eq:O.-U.}),
\BEQ \partial_t\phi=\nu\Del\phi+\eta \EEQ
where $\eta$ is a {\em regularized white noise}.
Our precise choice of regularization is the following: we define $\eta_{reg}$ to be a "kick force", namely, we choose
an infinite number of independent copies $(\xi^{reg}_{n+\half})_{n\in\N}$ of regularized space white noises and let
$\eta(t):=\xi^{reg}_{n+\half}$ be constant on $t\in(n,n+1)$. For definiteness we take 
$\xi^{reg}=e^{\nu^{(0)}\Del}\xi$, where $\xi$ is a standard space white noise.
Thus  $\tilde{\eta}$ is is 
the piecewise continuous in time, smooth in space, centered Gaussian process with  covariance 
\BEQ \esper\left[ {\eta}(t,x),\tilde{\eta}(t',x')
\right]=\del_1(t,t') p_{2\nu}(x-x'), \EEQ 
where:  $\del_1(t,t')=1$ if $t,t'$ are in the same unit time interval
$(n,n+1)$ for some $n\in\Z$, $0$ else; and  $p_{\tau}(x-x'):=\frac{1}{(2\pi\tau)^{d/2}} e^{-|x-x'|^2/2\tau}$ is the
standard heat kernel.
Note that the choice of a piecewise continuous "kick force" instead of a time delta-correlated
noise avoids the use of the stochastic calculus toolbox.  

\medskip

 Let $G=(\partial_t-\nu\Del)^{-1}$ be the  Green kernel of the
linear heat equation; formally, $\phi=G\eta$. Thus scale $j$ fluctuation fields  $\phi^j$ and $\eta^j$ should be
in direct link, namely, $\phi^j=G\eta^j$. A natural way to
accomplish this is to cut $G$ itself into scales, $G=\sum_j G^j$, and set $\phi^j=G^j\eta$, $\eta^j=(\partial_t-\nu
\Del)\phi^j$. 

\medskip

The {\em stationary} Ornstein-Uhlenbeck process,
\BEQ \phi(t,x)=\int_{-\infty}^t ds\,  e^{(t-s)\nu\Del}\eta_s(x),  \label{eq:phi} \EEQ
 solution of (\ref{eq:O.-U.}),
has covariance kernel (assuming
 e.g. $t\ge t'$)
\BEA  \esper[\phi(t,x)\phi(t',x')] & =& \int_{-\infty}^t ds  \int_{-\infty}^{t'} ds' \int dy dy'  p_{\nu(t-s)}(x-y)
p_{\nu(t'-s')}(x'-y') \del_1(s,s') p_{2\nu}(y-y') \nonumber\\
&\approx & \int_0^{+\infty} du \left(\int dy  dy' p_{\nu (t-t'+u)}(x-y) p_{\nu u}(x'-y') p_{2\nu}(y-y') \right). \label{eq:stat-cov1}\EEA
The regularization has a measurable effect only around the diagonal $t=t',x=x',u=0$. Away from the diagonal
the last integral (\ref{eq:stat-cov1}) behaves like $\int^{+\infty} du \ p_{\nu (t-t'+2u)}(x-x')=\int^{+\infty} du\,  \frac{e^{-|x-x'|^2/2
\nu(t-t'+2u)}}{(2\pi\nu(t-t'+2u))^{d/2}}$, an integrable function at infinity since $d/2>1$. Thus
\BEQ \left|  \esper[\phi_t(x)\phi_{t'}(x')] \right|\lesssim \int_0^{+\infty} du\,  (t-t'+u)^{-d/2} \left(
1+O\left(\frac{|x-x'|}{\sqrt{t-t'+u}}\right)\right)^{-N},\qquad N\ge 1 \\ \EEQ
is bounded by a constant times $(t-t')^{1-d/2}$ if $|x-x'|\lesssim \sqrt{t-t'}$, and 
by $\int_{|x-x'|^2}^{+\infty} s^{-d/2} ds=\frac{C}{|x-x'|^{d-2}}$ (the Green kernel of the Laplacian on $\R^d$)
in the contrary case.

\bigskip

We now want to cut $\phi$ into scales, i.e. understand how it behaves typically for time separations of order
$2^j$ ($j\ge 0$), or space separations of order $2^{j/2}$.  The main task is to cut $G$ into dyadic scales, $G=\sum_{j\ge 0}G^j$; then (as discussed above) we
 define
\BEQ \phi^j=G^j\eta,\qquad \eta^j=(\partial_t-\nu\Del)\phi^j. \EEQ
With these definitions, $\sum_{j\ge 0}\phi^j=G\eta=\phi$ is the Ornstein-Uhlenbeck field, 
and $\sum_{j\ge 0}\eta^j=(\partial_t-\nu\Del)G\eta=\eta$.

We proceed as follows. Let $\bar{\chi}:\R_+\to\R_+$ be a smooth 'bump' function of scale $1$
supported away from the origin, say, $\bar{\chi}\big|_{[2^{-\half},2^{\half}]}\equiv 1$, $\bar{\chi}\big|_{\R_+\setminus
(2^{-1},2)}\equiv 0$, chosen in such a way that 
\BEQ \bar{\chi}^{0}(\cdot):=\sum_{n\ge 0} \bar{\chi}(2^n \cdot),\qquad \bar{\chi}^j(\cdot):=\bar{\chi}(2^{-j}\cdot)
\qquad (j\ge 1) \EEQ 
form a partition of unity, i.e. $\sum_{j\ge 0} \bar{\chi}^j\equiv 1$ on $\R_+$, with
 $\supp \bar{\chi}^0\subset B(0,2)$, $\supp(\bar{\chi}^j)\subset
 B(0,2^{j+1})\setminus B(0,2^{j-1})$  $(j\ge 1)$.

\begin{Definition}[cut-off] \label{def:non-local-cut-off}
Let $G^j$ be the operator
\BEQ (G^j g)(t):=\int \bar{\chi}^j(s) e^{s\nu\Del} g(t-s) ds,\qquad
j\ge 0\EEQ
and
\BEQ \phi^j=G^j\eta,\qquad
\eta^j=(\partial_t-\nu\Del)\phi^j.\EEQ
\end{Definition}
Clearly, $\sum_{j\ge 0}G^j=G$ and $\sum_{j\ge 0}\phi^j=\phi$ is the solution of the Ornstein-Uhlenbeck
equation (\ref{eq:O.-U.}).
Note that, for $j\ge 1$, $\eta^j(t)=\int (\bar{\chi}^j)'(s) e^{s\nu\Del}\eta(t-s)\, ds
$ is smooth, while $\eta^0(t)=\eta(t)+\int (\bar{\chi}^0)'(s)  e^{s\nu\Del}\eta(t-s)\, ds$ has an extra "kick force" term. 

\medskip

Let $t\ge t'$. Assume $j\ge 1$. The diagonal covariance kernel $C^j_{\phi}(t,x;t',x')=\esper[\phi^j_{t}(x)\phi^j_{t'}(x')]$
 is non-zero only for $t-t'\lesssim 2^j$, in which case (recall $d_{\phi}:=\half(\frac{d}{2}-1)$)
\BEA  C^j_{\phi}(t,x;t',x')&\lesssim &\int_0^{2^j} du\,  \left(e^{\nu(M^{j-1}+u)\Del} e^{\nu(M^{j-1}+u-(t-t'))\Del}\right) p_{2\nu}(x-x') \nonumber\\
&\lesssim & 2^j p_{c\nu 2^j}(x-x') \nonumber\\
&\lesssim & (2^{-j})^{2d_{\phi}} e^{-c' 2^{-j/2} |x-x'|}\EEA
for some constants $c,c'>0$. A similar formula holds for $j=0$: if $t-t'\lesssim 1$,
\BEQ C^0_{\phi}(t,x;t',x') \lesssim \int_0^1 du \left(e^{\nu u\Del} e^{\nu(u-(t-t'))\Del}\right) p_{2\nu}(x-x')\lesssim  p_{c\nu}(x-x') \lesssim e^{-c'|x-x'|}. \EEQ

Then the off-diagonal covariances 
\BEQ C^{j,j'}_{\phi}(t,x;t',x')=\esper[\phi^j_{t}(x)\phi^{j'}_{t'}(x')]\EEQ
are similarly shown to satisfy for $j\ge j'$ the estimate
\BEQ |C^{j,j'}_{\phi}(t,x;t',x')|\lesssim 2^{j'} p_{c\nu 2^j}(x-x') \lesssim 2^{-|j-j'|} (2^{-j})^{2d_{\phi}}
e^{-c' 2^{-j/2} |x-x'|}.\EEQ
Since $C_{\phi}^{j,j'}(t,\cdot;t',\cdot)=0$ for $|t-t'|\gg 2^j$, one may clearly also write
\BEQ |C^{j,j'}_{\phi}(t,x;t',x')|\lesssim 2^{-|j-j'|} (2^{-j})^{2d_{\phi}}
e^{-c2^{-j}|t-t'|-c2^{-j/2} |x-x'|}.\EEQ

\medskip

Finally gradients applied to the heat kernel produce by standard parabolic estimates   small factors
of order $O(2^{-\max(j,j')/2})$. Let us recapitulate.

\begin{Lemma}[covariance kernel estimates] \label{lem:cov-ker-est}
Let
\BEQ C^{j,j'}_{\phi}(t,x;t',x')=\esper[\phi^j_{t}(x)\phi^{j'}_{t'}(x')], C^{j,j'}_{\eta}(t,x;t',x')=\esper[\eta^j(t,x)\eta^{j'}(t',x')]\EEQ 
and
\BEQ C^{j}_{\phi}:=C^{j,j}_{\phi},\qquad C^j_{\eta}:=C^{j,j}_{\eta}.\EEQ

Then, for $j\ge j'$, 
\BEQ \left|\nabla_x^{p} \nabla_{x'}^{p'} C^{j,j'}_{\phi}(t,x;t',x')\right|\lesssim
2^{-|j-j'|} 2^{-\frac{j}{2}(p+p')} (2^{-j})^{2d_{\phi}}
e^{-c2^{-j}|t-t'|-c2^{-j/2} |x-x'|} \\ \EEQ
and
\BEQ \left|\nabla_x^{p}\partial_t^q \nabla_{x'}^{p'}\partial_{t'}^{q'} C^{j,j'}_{\eta}(t,x;t',x')\right|\lesssim
2^{-|j-j'|} 2^{-\frac{j}{2}(p+p')} (2^{-j})^{2+2d_{\phi}}
e^{-c2^{-j}|t-t'|-c2^{-j/2} |x-x'|}. \\ \EEQ
Furthermore, if $j\ge 0$,
\BEQ \esper[(\eta^j_t(x)-\eta^j_t(y))^2]\lesssim (2^{-j})^{3+2d_{\phi}}|x-y|^2 \label{eq:eta-distance} \EEQ
and
\BEQ \esper[(\eta^j_t(x)-\eta^j_s(x))^2]\lesssim (2^{-j})^{4+2d_{\phi}}|t-s|^2 \label{eq:5.20bis}.\EEQ
\end{Lemma}

The  last two estimates (\ref{eq:eta-distance}), (\ref{eq:5.20bis})  follows immediately from Taylor's formula: letting $\vec{v}:=\frac{y-x}{|y-x|}$,
\BEQ \esper[(\eta^j_t(x)-\eta^j_t(y))^2]\le \int_0^{|y-x|} dz \int_0^{|y-x|} dz' 
\left|\nabla_{\vec{v}}\nabla'_{\vec{v}} C^j_{\eta}(t,x+z\vec{v};t,x+z'\vec{v}) \right|\EEQ
and similarly for $\esper[(\eta^j_t(x)-\eta^j_s(x))^2]$.

\medskip

One has thus obtained a very elaborate version of the scalings (\ref{eq:scaling}), $\phi^j(t,x)=O(2^{-jd_{\phi}}),
\eta^j(t,x)=O((2^{-j})^{1+d_{\phi}})$, together with a first indication of the scale-separation mechanism: the prefactors
in powers of $2^{-|j-j'|}$ show clearly that fields of widely separated scales are effectively independent.

\medskip
\noindent {\em Remark.} Note that the {\em low-momentum fields}, $\phi^{\to j}(t,x):=
\sum_{k\ge j} \phi^k(t,x), \eta^{\to j}(t,x):=\sum_{k\ge j}\eta^k(t,x)$ verify
the same scaling as the single-scale fields, namely, $\phi^j(t,x)=O(2^{-jd_{\phi}}),
\eta^j(t,x)=O((2^{-j})^{1+d_{\phi}})$.


\section{Appendix. Large deviations estimates for the single-scale noisy equation}


\subsection{Introduction}


We consider here the noisy KPZ equation with scale $j$ infra-red cut-off,
\BEQ \partial_t \psi=(\nu\Del-2^{-j})\psi+\lambda V(\nabla\psi)+\eta^j\EEQ
with right-hand side $\eta^j=G^j \eta$ defined as in section 5.
Recall the conclusion of the discussion at the end of \S 4.2: by
Theorem 2 (see Introduction), if $|||\eta^j|||_{{\cal W}^{1,\infty;2\lambda}}([0,t],x)=O(2^{-jd_{\phi}})$, then $|\psi(t,x)|=O(2^{-jd_{\phi}}), |\nabla\psi(t,x)|=O((2^{-j})^{\half+d_{\phi}})$. 

\medskip
\noindent
We  show in this section that $|||\eta^j|||_{{\cal W}^{1,\infty;2\lambda}_j}([0,\infty),x)$ is a.s. bounded, and
prove large deviation estimates for this quantity when it is much larger than $O(2^{-jd_{\phi}})$.  Contrary to the previous sections, this one is of essentially probabilistic nature. Non-specialists who are not particularly interested in stochastic PDEs may
safely skip it.

\medskip
\noindent
The random variables appearing in the definition of the pointwise "quasi-norms" associated with  ${\cal W}^{1,\infty;2\lambda}$ are essentially time- and space-averages of a large number of independent {\em log-normal variables}, such as
$e^{4\lambda 2^j|\eta^j(t,x)|}$. Log-normal variables  have large tails in $e^{-a(\ln z)^2}$ and thus no exponential moment,
hence standard large-deviation theory (notably Cram\'er's theorem) does not give any valuable information on the
probability that such averages become large. Some authors have been considering this problem, notably Russians, starting from the 60es; one may cite Linnik \cite{Lin},
Nagaev \cite{Nag1,Nag2}, Rozovski \cite{Roz}, see also e.g. Kl\"uppelberg and Mikosch \cite{KluMik} for a renewal of
the theory with a view to applications in insurance. The theory is not easily accessible, partly because written
originally in Russian journals in the 60es and 70es (in particular in {\em Teoriya Veroyatnostei i ee Primeneniya}, later translated to English as
{\em Theory of Probability and its  Applications}), partly for the lack of a theory as general and satisfactory as the standard large-deviation
theory. 

\medskip
Let us just point out the difficulties (this very short abstract is taken from an inspiring review in \cite{MikNag}).
Choose a random variable $X$ with finite first and second moments; by translation and rescaling we may assume
that $\esper[X]=0,\esper[X^2]=1$. Let $S_n:=X_1+\ldots+X_n$, $M_n:=
\max(X_1,\ldots,X_n)$, where $X_1,\ldots,X_n$ are independent copies of $X$. Let finally 
$\bar{F}^X(x):=\proba[X>x],\ \bar{F}_n^X(x):=\proba[S_n>x]$ and $\Errfc(x):=\int_x^{+\infty}
 \frac{e^{-y^2/2}}{\sqrt{2\pi}} \, dy$ be resp. the queues of $X$, of $S_n$ and of a standard Gaussian variable. By the central limit theorem,
one expects 
\BEQ \bar{F}_n^X(x)\approx \Errfc(x/\sqrt{n}), \label{eq:LD1}\EEQ
at least if $x\approx \sqrt{n}$. On the other hand, one clarly has if $X\ge 0$
\BEQ \bar{F}_n^X(x)\ge \proba[M_n>x]\sim_{x\to\infty} n\bar{F}^X(x). \label{eq:LD2}\EEQ
{\em Subexponential distributions} (including log-normal distributions) are precisely defined by the asymptotic relation
$\bar{F}_n^X(x)\sim_{x\to\infty}  n\bar{F}^X(x)$, implying a heavy queue. For distribution with lighter queues
 (such as e.g. Gaussian distributions), the inequality in (\ref{eq:LD2}) is very rough, in the sense that typically $n\bar{F}^X(x)\ll \bar{F}_n^X(x)$ for every $x\ge x_0$,
with $x_0$ independent from $n$. 

Thus, one expects, specifically for subexponential distributions, a {\em central limit theorem behaviour} as in (\ref{eq:LD1}) for
$x\ll c_n$, with $c_n$ defined by $\Errfc(c_n/\sqrt{n})\approx \bar{F}_n^X(c_n)$, and an {\em extreme-value regime},
\BEQ \bar{F}^X_n(x)\sim n\bar{F}^X(x),\qquad x\gg d_n\EEQ
with $d_n\ge c_n$, in which $n\bar{F}^X(x)\gg \Errfc(x/\sqrt{n})$.
Optimal sequences $c_n,d_n$ have been identified for various types of subexponential distributions; for a
standard log-normal variable $X=e^Z$, $Z\sim {\cal N}(0,1)$, one finds $c_n,d_n\approx n^{\half}\ln n$. One major
drawback of this picture is that it doesn't say anything about the behaviour of $\bar{F}_n^X(x)$ in the window
$c_n\lesssim x\lesssim d_n$ (in our case, for $x\approx c_n$ since $c_n=d_n$), which is expected to be a mixture
of (\ref{eq:LD1}) and (\ref{eq:LD2}). The complicated asymptotics, valid on the whole real line, proved by 
Rozovski \cite{Roz} -- a veritable tour de force -- give a more complete answer. 

\medskip
This being said, our problem does not fit exactly into this frame, since (1) we are only interested in {\em 
upper bounds} for $\bar{F}_n^X$, moreover in the {\em extreme-value regime}, with $x\gtrsim n$; on the other
hand (2) the variables $X_1,\ldots,X_i,\ldots,X_j,\ldots,X_n$  (chosen as local space or space-time averages of
the noise) are not independent, but have {\em correlations which decrease exponentially}
with the scaled distance $d^j$ (see below) or equivalently with $|j-i|$; (3) we need {\em scale-dependent estimates} for $\bar{F}^X_n$ since $X\approx e^{2^{-jd_{\phi}}|Z|}$,
$Z\sim {\cal N}(0,1)$ is strongly $j$-dependent. However all the previous results are strongly dependent on the
particular form of the distribution, in particular on the first and second moments, and it is often difficult
to retrace the $j$-dependence of the constants in the bounds.

\bigskip

Our main  result is the following. 

\begin{Theorem} \label{th:LD1}
Let $j\in\N$ and $\lambda>0$. Then the function $x\mapsto |||\eta^j|||_{{\cal W}^{1,\infty;\lambda}_j(\R_+)}(x)$
is a.s. everywhere defined (i.e. finite). Furthermore, the following large deviation estimates holds for every
$x\in\R^d$,
\BEQ \proba[\sup_{B(x,2^{j/2})} |||\eta^j|||_{{\cal W}^{1,\infty;\lambda}_j(\R_+)}(x)>A 2^{-jd_{\phi}}]\lesssim 
A^{-c\ln(A)}, \qquad A\ge 1\EEQ
where $c>0$ is some constant.
\end{Theorem}

As follows from Theorem 2, this implies (up to the replacement of $\lambda$ by $2\lambda$) that the solution $\psi$ of the full KPZ equation with scale $j$ infra-red cut-off (\ref{eq:KPZdowntoj}) is defined a.s. for all positive times $t\ge 0$ and sits in 
the space ${\cal W}^{1,\infty;2\lambda/5}_j$, with $|||\psi|||_{{\cal W}^{1,\infty;2\lambda/5}_j}(x)=O(2^{-jd_{\phi}})$ for every $x\in\R^d$, with a random multiplicative prefactor $A(x)$ whose queue is bounded {\em locally
in $x$} by that of
a log-normal distribution. (Note that the prefactor $A(x)$ is {\em not globally} bounded!)

\bigskip

The proof  includes both Gaussian inequalities taken from the monograph \cite{Adl}, and
an adaptation to weakly correlated variables of a result about large deviations for subexponential distributions
\cite{Nag2}. We shall need quite a few preliminary results before the proof, given at the very end of the present
section. 

We finish this introductory paragraph with the tiny bit of stochastic domination and Gaussian inequalities used
in the sequel, and a little bit of geometry.

\medskip

\begin{Definition}
Let $X:\Omega\to\R,Y:\Omega'\to\R$ be two real-valued random variables, defined a priori on two different
probability spaces. Then $X$ is {\em stochastically dominated} by $Y$ if
\BEQ \forall x\in\R, \proba[X>x]\le \proba[Y>x]. \EEQ
We then write $X\lele Y$.
\end{Definition}
By Strassen's theorem \cite{Lind}, if $X\lele Y$, there exists a coupling between $X$ and $Y$, i.e. random variables $X',Y':\Omega''\to\R$ defined on the same probability space, with $X'\stackrel{(d)}{=}X,\ Y'\stackrel{(d)}{=} Y$, 
and such that $X'\le Y'$.

\medskip

\begin{Proposition} (see \cite{Adl}) \label{prop:Slepian}
Let $(Z_1,\ldots,Z_n)$ be a centered Gaussian vector, and $\phi:\R^n\to\R$ be a convex function with
polynomial growth at infinity. Then $\esper[\phi(Z_1,\ldots,Z_n)]$ is an increasing function of the coefficients
$c_{ij}=\esper[Z_i Z_j]$, $i,j=1,\ldots,n$.
\end{Proposition}
This technical lemma, due to Slepian (whose short proof relies on a Gaussian integration by parts) is one of the main tools for
Gaussian inequalities. It extends to Gaussian fields $(Z_x)_{x\in\R^d}$ and convex functionals $\phi$ under adequate regularity assumptions.

\medskip

\begin{Proposition}[Borell-Tsirelson-Ibragimov-Sudakov or BTIS inequality] (see \cite{Adl}) \label{prop:BTIS}\\
Let $(Y_x)_{x\in {\cal D}}$, ${\cal D}\subset B(0,1)$ be a centered Gaussian process, such that 
$\sigma^2_{{\cal D}}=\sup_{x\in{\cal D}} \esper[Y^2_x]<\infty$, and
$\del(x,y):=\sqrt{\esper[ (Y_x-Y_y)^2]}\lesssim |x-y|.$ Let $||Y||_{\infty}:=\max_{x\in{\cal D}} |Y_{x}|$. Then a.s.
$||Y||_{\infty}<\infty$, $\esper[||Y||_{\infty}]\lesssim 1$ and
\BEQ \proba[||Y||_{\infty}-\esper[||Y||_{\infty}]>u]\le e^{-u^2/2\sigma_{{\cal D}}^2}.\EEQ
\end{Proposition}
This is actually a particular case of the BTIS inequality. For a Gaussian process $Y$ indexed by an abstract set $\cal D$, $\esper[||Y||_{\infty}]$ is bounded by the integral of the square-root of the entropy $\log N(\eps)$, $\esper[||Y||_{\infty}]
\lesssim \int_0^{+\infty} \sqrt{\ln N(\eps)}d\eps$, where $N(\eps)$ is the minimum number of balls of diameter
$\le\eps$ (with respect to the metric induced by $\del(\cdot,\cdot)$) covering $\cal D$. In our proposition, 
$\ln N(\eps)=0$ for $\eps\gg 1$ since $\sup_{x,y\in{\cal D}} \del(x,y)\lesssim 1$, and $N(\eps)=O(\eps^d)$
otherwise by hypothesis, hence the result.

\medskip
The above proposition applies for fixed $t_0,x_0$ to $Y_x:=2^{j(1+d_{\phi})}\eta^j(t_0,x_0+2^{j/2}x)$, 
with ${\cal D}=B(0,1)$. It follows from  Lemma \ref{lem:cov-ker-est} in Appendix A that $\sigma^2_{{\cal D}}\approx 1$ and $d(x,y)\lesssim |x-y|$. Thus 
\BEQ \esper \left[ \sup_{B(x_0,2^{j/2})} |\eta^j_{t_0}| \right] \lesssim 2^{-j(1+d_{\phi})} \EEQ
and there exists a constant $C\lesssim 1$ such that
\BEQ \proba[2^{j(1+d_{\phi})} \sup_{B(x_0,2^{j/2})} |\eta^j_{t_0}|>u+C]\le e^{-u^2/2C}.\EEQ
One easily deduces that 
\BEQ 2^{j(1+d_{\phi})} \sup_{B(x_0,2^{j/2})} |\eta^j_{t_0}|\lele C' (|Z|+1) \label{eq:suplele}\EEQ
if $Z\sim {\cal N}(0,1)$.

\bigskip

Recall from section 3.1 that $f^*\le f^{\sharp}$ $(f\in C(\R^d,\R))$  -- note, and this is very important, that  the inequality is {\em exact}, with 
a coefficient one --, where $f^{\#}(x)=\sup_{\rho>0} \fint_{B(x,\rho)}
|f|.$ We cannot bound directly a supremum over a continuous parameter (here $\rho$), so it is natural to start
by rewriting $(\eta^j_{t_0})^*$ in terms of its local averages or suprema on balls of radius $2^{j/2}$, over which
we have a good control. However, we cannot obviously cover $\R^d$ (nor $B(x,\rho)$) by disjoint balls of fixed
radius, and taking into account error terms due to overlaps or boundary effects would cost a multiplicative coefficient,
which we cannot afford to do. Hence we first transform balls centered at $x$ into cubes by letting
\BEQ \Phi:\R^d\to\R^d,\quad y\mapsto \Phi(y)=x+\frac{|y-x|}{||y-x||_{\infty}}(y-x)\EEQ
where $||y-x||_{\infty}:=\sup(|y_1-x_1|,\ldots,|y_d-x_d|)$ is the supremum norm. The Euclidean norm $|\ \cdot\ |$ and
the supremum norm $||\ \cdot\ ||_{\infty}$ are equivalent, hence (the easy proof is left to the reader)
$\Phi$ and $\Phi^{-1}$ are uniformly Lipschitz.  Thus $\sup_{B(x,\rho)} |\eta^j_{t_0}|=\sup_{\tilde{B}(x,\rho)}
|\eta^j_{t_0}\circ\Phi^{-1}|$, where $\tilde{B}(x,\rho)=\{y\in\R^d\ |\ ||y-x||_{\infty}=\rho\}$ is a cube.
The field $\eta^j_{t_0}\circ\Phi^{-1}$ has the same general properties as $\eta^j_{t_0}$ (scaling, exponentially
decreasing covariance) as stated in Lemma \ref{lem:cov-ker-est}, so (by abuse of notation) we simply denote
$\eta^j_{t_0}\circ\phi^{-1}$ by $\eta^j_{t_0}$ in the sequel.

\begin{Definition}[scale $j$ cubes]
Let
$\D^j$ be the set of all {\em scale $j$ cubes}, i.e. of all primitive cells $[k_1 2^{j/2},(k_1+1)2^{j/2}]\times\ldots\times[k_d 2^{j/2},(k_d+1)2^{j/2}]$,
$k_1,\ldots,k_d\in\Z$  
of the square lattice $2^{j/2}\Z^d$.

We denote by $x_{\Del}=(x_{\Del,1},\ldots,x_{\Del,d})$ the center of a cube $\Del\in\D^j$.
\end{Definition}

We now show how to bound an average  $\fint_{\tilde{B}(x,\rho)} |f|$, $f\in C(\R^d,\R)$ over a cube of arbitrary radius in terms of the 
local suprema $f_{\Del}:=\sup_{\Del} |f|,\ \Del\in \D^j$. We give the proof in dimension 2 to simplify notations
(in general, we would need the whole cellular decomposition of a cube). Let, for $\rho>0$,
\BEQ \tilde{B}^j(x,\rho):=\cup\{\Del\in\D^j\ |\ \Del\subset \tilde{B}(x,2^{j/2}\rho)\},\qquad
\partial \tilde{B}^j(x,\rho):=\cup\{\Del\in\D^j\ |\ \Del\cap \tilde{B}(x,2^{j/2}\rho)\not=\emptyset\}\setminus
\tilde{B}^j(x,\rho),  \label{eq:Bjxrho1} \EEQ
and $n:=\tilde{B}^j(x,\rho)$.
The boundary $\partial\tilde{B}^j(x,\rho)$ decomposes into 8 disjoint subsets,
\BEQ \partial\tilde{B}^j_{{\mathrm{right}}}(x,\rho):=\cup\{\Del=[x_{\Del,\min},x_{\Del,\max}]\times
[y_{\Del,min},y_{\Del,max}]\ | \ a<x_{\Del,min}<b<x_{\Del,max},\ c\le y_{\Del,min}<y_{\Del,max}\le d\} \EEQ
and similary $ \partial\tilde{B}^j_{{\mathrm{left}}}(x,\rho), \partial\tilde{B}^j_{{\mathrm{up}}}(x,\rho), \partial\tilde{B}^j_{{\mathrm{down}}}(x,\rho)$ for the sides of the square;
\BEQ  \partial\tilde{B}^j_{{\mathrm{up,right}}}(x,\rho)=\cup\{\Del=[x_{\Del,\min},x_{\Del,\max}]\times
[y_{\Del,min},y_{\Del,max}]\ |\ a<x_{\Del,min}<b<x_{\Del,max},\ c< y_{\Del,min}<d<y_{\Del,max}\} \EEQ
and similary for the three other corners. We let $c_{{\mathrm{right}}}:=\frac{\Vol(\partial\tilde{B}^j_{{\mathrm{right}}}
(x,\rho)\cap \tilde{B}(x,\rho))}{\Vol((\partial\tilde{B}^j_{{\mathrm{right}}}
(x,\rho))}$, and similarly $c_{{\mathrm{left}}},\ldots$ be the corresponding volume ratios. Let
\BEQ F(c_{{\mathrm{right}}},c_{{\mathrm{left}}},\ldots):=\frac{\sum_{\Del\in \tilde{B}^j(x,\rho)} f_{\Del}+
c_{{\mathrm{right}}} \sum_{\Del'\in\partial\tilde{B}^j_{{\mathrm{right}}}(x,\rho)} f_{\Del'}+\ldots}{n+
c_{{\mathrm{right}}} \sharp \partial\tilde{B}^j_{{\mathrm{right}}}(x,\rho)+\ldots};\EEQ
note that $\fint_{\tilde{B}(x,\rho)} f\le F(c_{{\mathrm{right}}},\ldots)$ since $c_{{\mathrm{right}}}$ is the
uniform volume ratio $\frac{\Vol(\Del)\cap \tilde{B}^j(x,\rho)}{\Vol(\Del)}$ of all scale $j$ cubes at the right
border, as follows from the fact that the border is straight. Then trivially $F(c_{{\mathrm{right}}},c_{{\mathrm{left}}},\ldots)\le \max\left(  F(0,c_{{\mathrm{left}}},\ldots), F(0,c_{{\mathrm{right}}},\ldots) \right)$; this same elementary
remark may be repeated for the eight $c$ coefficients. Thus we have proved that
\BEQ \fint_{\tilde{B}(x,\rho)} f\le \max_{B^j} \frac{\sum_{\Del\in B^j} f_{\Del}}{\sharp B^j},\EEQ
where the $B^j$ range among
$2^8$ subsets of squares, and by definition  $\tilde{B}^j(x,\rho)\subset B^j\subset \tilde{B}^j(x,\rho)\cup \partial\tilde{B}^j(x,\rho)$.


\subsection{A first preliminary result: large deviations for the noise}


We prove in this paragraph the following result.

\begin{Lemma} \label{lem:LD-main}
Let $j\in\N$ and $t_0\in\R_+$. Then the function $x\mapsto  (\eta^j)^*(t_0,x)$ is a.s. everywhere defined
(i.e. finite). Furthermore, the following large deviation estimates holds,
\BEQ \proba[\sup_{B(x,2^{j/2})} (\eta^j)^*(t_0)>A 2^{-j(1+d_{\phi})}]\le e^{-c(A-C)_+^2} \label{eq:6.14} \EEQ 
for some constants $c,C>0$, where $(A-C)_+^2=(A-C)^2 {\bf 1}_{A>C}.$
\end{Lemma}

It is actually reasonable to expect, on account of the central limit theorem, that $|\eta^j(t_0,x)|-\esper[|\eta^j(t_0,x)|]\in {\cal H}^0_{\alpha}$
for every $\alpha<d/4$,  and that the norm in ${\cal H}^0_{\alpha}$ satisfies large deviation estimates
as in (\ref{eq:6.14}), but we do not prove this. 
The above result, however natural it may be, is not really needed anywhere in the article, but the proof of Theorem
\ref{th:LD} is based on the arguments developed for the proof of the lemma. \medskip

{\bf Proof.} In the sequel $c,c',C>0$ are constants possibly varying from line to line (contrary to $c_0,m_0$,
see below, which are fixed once and for all). As already recalled,  $(\eta^j)^*(t_0,x)\le(\eta^j)^{\sharp}(t_0,x)=
\sup_{\rho>0}  \fint_{B(x,2^{j/2} \rho)} dy |\eta^j(t_0,y)|.$  Also, from the results of Appendix A,
the correlations of the field $(\eta^j_{t_0}(x))_{x\in\R^d}$ decay exponentially with the
{\em scaled distance} $d^j(x,x'):=2^{j/2}|x-x'|$, in the sense that, for a certain constant $c_0$, 
\BEQ |\esper[\eta^j_{t_0}(x)\eta^j_{t_0}(x')]|\lesssim 2^{-2j(1+d_{\phi})} e^{-c_0 d^j(x,x')}.  \label{eq:c0} \EEQ

\medskip We split the proof into several points.

\begin{itemize}
\item[(i)] In order to use the exponential decay, we first choose $m_0\ge 2$ large enough (depending on further
considerations), and partition $\D^j$ into $m_0^d$ disjoint susets $\D^j_{\mu}$, $\mu\in\{1,\ldots,m_0\}^d$, 
with $\Del=[k_1 2^{j/2},(k_1+1)2^{j/2}]\times\ldots\times[k_d 2^{j/2},(k_d+1)2^{j/2}]\in\D^j_{\mu}\Leftrightarrow
k_i\equiv \mu_i\mod m_0$. Two points $x,x'$ located in disjoint cubes $\Del\not=\Del'$ in the same sublattice
$\D^j_{\mu}$ are thus at distance $d^j(x,x')\gtrsim m_0$, which amounts (up to rescaling) to replacing $c_0$ by $m_0 c_0$ in (\ref{eq:c0}); in the sequel, we may thus assume that $c_0$ is large enough. By abuse of notation, we also denote by $\D^j_{\mu}$ the subset $\cup \{\Del;\Del\in\D^j_{\mu}\}
\subset\R^d$.
Clearly, 
\BEQ \fint_{B(x,2^{j/2}\rho)} dy |\eta^j(t_0,y)|\le \sup_{\mu} \fint_{B(x,2^{j/2}\rho)\cap\D^j_{\mu}} dy |\eta^j(t_0,y)|.
\label{eq:6.16} \EEQ
If $\psi:\R\to\R$ is increasing, then
\BEQ \esper\left[\psi\left(\fint_{B(x,2^{j/2}\rho)} dy |\eta^j(t_0,y)|\right)\right]\le \sum_{\mu}
\esper\left[\psi\left(\fint_{B(x,2^{j/2}\rho)\cap\D^j_{\mu}} dy |\eta^j(t_0,y)|\right)\right] \label{eq:mu}.\EEQ

\item[(ii)] Next, we want to bound the average $\fint_{B(x,2^{j/2}\rho)\cap\D^j_{\mu}} dy |\eta^j(t_0,y)|$
over some fixed sublattice by the average of a finite number of variables representing the supremum of $|\eta^j|$ on
each cube. For that (note that the following construction is $\mu$-dependent, which we do not always specify)
we introduce i.i.d. copies $(\eta^j_{\Del})_{\Del\in\D'_{\mu}}$ of the field $\eta^j_{t_0}\big|_{[-\half 2^{j/2},
\half 2^{j/2}]^d}$ restricted to some reference cube, and define a new random field $\tilde{\eta}^j$ on $\D^j_{\mu}$,
 \BEQ \tilde{\eta}^j(x):=\sum_{\Del'\in\tilde{B}^j(x,\rho)} e^{-c_0 d^j(\Del,\Del')} \eta^j_{\Del'}(x-x_{\Del}),\qquad x\in\Del \label{eq:new-random-field} \EEQ
 separately on each cube $\Del\in\D^j_{\mu}$, where
\BEQ d^j(\Del,\Del'):=2^{-j/2}\sup_{x\in\Del}\inf_{y\in\Del'} |x-y| \EEQ
is the set distance measured in scaled units, and 
\BEQ \tilde{B}^j(x,\rho):=\{\Del\in\D^j_{\mu} \ |\ 
\Del\subset \tilde{B}(x,2^{j/2}\rho)\} \label{eq:Bjxrho2}\EEQ
(compare with the previous definition, (\ref{eq:Bjxrho1})). 
By a simple computation, one finds 
\BEQ \esper[\tilde{\eta}^j(x)\tilde{\eta}^j(x')]\approx (1+d^j(\Del,\Del'))^d e^{-c_0 d^j(\Del,\Del')} 
\esper[\eta^j_{t_0}(x-x_{\Del})\eta^j_{t_0}(x'-x_{\Del'})]\gtrsim \esper[\eta^j_{t_0}(x)\eta^j_{t_0}(x')] \EEQ
if $x\in\Del,\ x'\in\Del'$ and 
\BEQ \Del,\Del'\in \tilde{B}^j(x,\rho)\cup\partial \tilde{B}^j(x,\rho):=\{\Del\in\D^j_{\mu} \ |\ 
\Del\cap \tilde{B}(x,2^{j/2}\rho)\not=\emptyset\}.\EEQ
Applying Proposition \ref{prop:Slepian} with $\phi(\eta^j_{t_0})=\psi\left( 2^{j(1+d_{\phi})}
\fint_{B(x,2^{j/2}\rho)\cap\D'_{\mu}}dy |\eta^j_{t_0}(y)| \right)$ where $\psi$ is any convex, increasing function
\footnote{Observe that $\psi_1\circ\psi_2$ is convex if $\psi_1:\R_+\to\R$ is convex and increasing and 
$\psi_2:\R^n\to\R_+$ is convex, since $\nabla^2(\psi_1\circ\psi_2)=\psi''_1\circ\psi_2 \ \cdot\ \nabla\psi_2\otimes
\nabla \psi_2 + \psi'_1\circ\psi_2 \ \cdot\ \nabla^2 \psi_2$.} on $\R_+$,
\BEQ \esper\left[\psi\left( 2^{j(1+d_{\phi})}
\fint_{B(x,2^{j/2}\rho)\cap\D'_{\mu}}dy |\eta^j_{t_0}(y)| \right)\right] \le 
\esper\left[\psi\left( 2^{j(1+d_{\phi})}
\fint_{B(x,2^{j/2}\rho)\cap\D'_{\mu}}dy |\tilde{\eta}^j(y)| \right)\right]. \label{eq:6.21} \EEQ

As follows from the discussion in the previous paragraph,
\BEQ \fint_{B(x,2^{j/2}\rho)\cap \D^j_{\mu}} dy |\tilde{\eta}^j(y)|\le \max_{B^j} \frac{\sum_{\Del\in B^j} \tilde{Y}_{\Del}}{\sharp B^j},\EEQ
where
\BEQ \tilde{Y}_{\Del}:=\sup_{\Del}|\tilde{\eta}^j|\EEQ
and the $B^j$ are a finite number (depending only on $d$) of subsets of cubes such that 
$\tilde{B}^j(x,\rho)\subset B^j\subset \tilde{B}^j(x,\rho)\cup\partial\tilde{B}^j(x,\rho).$

\item[(iii)]
By construction, see (\ref{eq:new-random-field}),
\BEQ \tilde{Y}_{\Del}\le \sum_{\Del'\in\tilde{B}^j(x,\rho)} e^{-c_0 d^j(\Del,\Del')} \sup|\eta^j_{\Del'}|. \label{eq:6.27}\EEQ
We have seen in (\ref{eq:suplele}) that $2^{j(1+d_{\phi})} \sup|\eta^j_{\Del'}|\lele C(|Z_{\Del'}|+1)$ if
$Z_{\Del'}\sim {\cal N}(0,1)$. Since the fields $(\eta^j_{\Del})_{\Del}$ are independent, we may by the above
cited Strassen theorem define a coupling of the field $\eta^j\big|_{\D^j_{\mu}}$ with i.i.d. standard Gaussian
variables $(Z_{\Del})_{\Del\in\D^j_{\mu}}$ in such a way that 
\BEQ 2^{j(1+d_{\phi})} \sup|\eta^j_{\Del}|\le C(|Z_{\Del}|+1). \label{eq:coupling} \EEQ
Hence 
\BEQ 2^{j(1+d_{\phi})} \sum_{\Del\in B^j} \tilde{Y}_{\Del}\le C\sum_{\Delta\in \tilde{B}^j(x,\rho)}(|Z_{\Del}|+1) \EEQ
-- note that the  bound in the right hand side does not depend on the choice of $B^j$ --  and
\BEQ \esper\left[ \psi\left( 2^{j(1+d_{\phi})}
\fint_{B(x,2^{j/2}\rho)\cap\D'_{\mu}}dy |\tilde{\eta}^j(y)| \right)\right] \le  \esper 
\left[  \psi\left( \frac{C}{n} \sum_{\Del\in\tilde{B}^j(x,\rho)} (|Z_{\Del}|+1) \right)
 \right].
        \label{eq:6.29} \EEQ

We rewrite the expectation as an integral by integration by parts,

\BEQ \esper 
\left[  \psi\left( \frac{C}{n} \sum_{\Del\in B^j} (|Z_{\Del}|+1) \right)
 \right] = \int_0^{+\infty} dA\, \psi'(A)\,  \proba\left[\frac{C}{n} \sum_{\Del\in B^j} (|Z_{\Del}|+1) >A\right]
 +\psi(0).
\EEQ

 Finally, $\sum_{\Del\in{B}^j} |Z_{\Del}|$ is
a sum of $n$ independent copies of $|Z|$, where $Z\sim {\cal N}(0,1)$, to
which
we may apply standard large deviation arguments in a trivial setting, 
\BEA  \proba[\sum_{\Del\in{B}^j} |Z_{\Del}|>nA] & \le & \min\left( 1, \min_{t\ge 0} e^{-tnA} \esper[e^{t\sum_{\Del} |Z_{\Del}|}] \right) \le 
\min\left( 1, 2^n \min_{t\ge 0} e^{-tnA+nt^2/2} \right) \nonumber\\
&=& \min\left( 1, 2^n e^{-nA^2/2} \right)\le C e^{-n(A-C)_+^2/2}.  \label{eq:standard-LD} \EEA
Thus we may choose $\psi(A)=e^{cn(A-C)_+^2} {\bf 1}_{A>C}+{\bf 1}_{A\le C}$ so that 
 $\esper\left[\psi\left( \frac{C}{n}\sum_{\Del\in B^j} (|Z_{\Del}|+1) \right)\right]\lesssim 1$.
Collecting  (\ref{eq:mu}), (\ref{eq:6.21}) and (\ref{eq:6.29}), one obtains
by Markov's inequality 
\BEQ \proba[M^{j(1+d_{\phi})} \fint_{B(x,2^{j/2}\rho)} dy |\eta^j_{t_0}(y)|>A]\lesssim \frac{1}{\psi(A)}\lesssim 
e^{-cn(A-C)^2},\qquad 
A\ge C.\EEQ 

For each fixed $n\ge 1$, the set $\{B^j(x,\rho),\rho\ge 0 \ |\ \sharp B^j(x,\rho)=n\}\cup\{
B^j(x,\rho)\cup\partial B^j(x,\rho),\rho\ge 0 \ |\ \sharp B^j(x,\rho)\cup\partial B^j(x,\rho)=n\}$
consists of $0,1$ or $2$ elements.
Thus, using (\ref{eq:mu}),
\BEQ \proba[(\eta^j)^*(t_0,x)>A 2^{-j(1+d_{\phi})}] \lesssim \min\left(1, \sum_{n\ge 1}  e^{-cn(A-C)_+^2} \right)
\lesssim
 e^{-c(A-C')_+^2}.\EEQ

Finally, we use a scaled version of (\ref{eq:fsharp-bounded}), 
\BEQ \sup_{B(x,2^{j/2})}(\eta^j_{t_0})^{\sharp}\lesssim \sup_{B(x,22^{j/2})} |\eta^j_{t_0}|+(\eta^j_{t_0})^{\sharp}(x),\EEQ
from which we conclude that
\BEQ \proba[\sup_{B(x,2^{j/2})} (\eta^j_{t_0})^*>A 2^{-j(1+d_{\phi})}]\le e^{-c(A-C)_+^2}. \EEQ 
In particular,
\BEQ \proba[\exists x\in\R^d\ |\ (\eta^j_{t_0})^*(x)=+\infty]\le \sum_{\Del\in\D^j} \proba[\sup_{\Del} (\eta^j_{t_0})^*=
+\infty]=0.\EEQ

\end{itemize}

\hfill \eop


\subsection{Large deviations for the exponential of the noise}


We now turn to large deviation estimates for $(e^{\lambda 2^j |\eta^j_{t_0}|})^*(x)$ and prove
the following result.

\begin{Theorem} \label{th:LD}
Let $j\in\N$, $\lambda>0$ and $t_0\in\R_+$. Then the function $x\mapsto (e^{\lambda 2^j |\eta^j_{t_0}|})^*(x)$
is a.s. everywhere defined (i.e. finite). Furthermore, the following large deviation estimates holds for every
$x\in\R^d$,
\BEQ \proba[ \sup_{B(x,2^{j/2})} \ln (e^{\lambda 2^j |\eta^j_{t_0}|})^*(x)>\eps A]\lesssim 
A^{-c\ln(A)}, \qquad A\ge 1\EEQ
where $\eps=\lambda 2^{-jd_{\phi}}$ and $c>0$ is some constant.
\end{Theorem}

The proof is essentially similar to that of Lemma \ref{lem:LD-main}, except that it is based on large deviation
estimates for log-normal variables. We cite a result by Nagaev,  show how to apply it in our context, and prove
a few technical lemmas before
turning to the proof of Theorem \ref{th:LD}.


\subsubsection{Log-normal large deviations}


\begin{Proposition} (see \cite{Nag2}, Corollary 1.8)
Let $X$ be a real-valued random variable such that $\esper[X]=0$ and $\esper[|X|^t]<\infty$ for some $t\ge 2$, and
$X_1,\ldots,X_n$ n i.i.d. copies of $X$, $S_n:=X_1+\ldots+X_n$. Then
\BEQ \proba[S_n> A]\lesssim n\esper[X^t {\bf 1}_{X>0}] A^{-t}+e^{-2(t+2)^{-2} e^{-t} A^2/n\esper[X^2]}. \label{eq:Nag}\EEQ
\end{Proposition}
Note that this general bound mixes the two regimes (\ref{eq:LD1}) and (\ref{eq:LD2}).

\begin{Corollary} \label{cor:Nagaev0}
Let $(Z_i)_{i=1,\ldots,n}$, $n\ge 1$ be i.i.d. standard Gaussian variables, and
let 
\BEQ S_n:=\sum_{i=1}^n \left( e^{\eps|Z_i|}-\esper[e^{\eps|Z_i|}] \right),\qquad \tilde{S}_n:=
\sum_{i=1}^n \left( e^{\eps|Z_i|}-1\right)\EEQ
where $0<\eps\ll 1$. Let finally $A\gg n\eps$ and $B\gg \ln(n)$. Then there exists a constant $c>0$ such that
\begin{itemize}
\item[(i)]
\BEQ \proba[S_n>A]\lesssim (A/\eps)^{-c\ln(A/\eps)} \label{eq:cor:Nag01}\EEQ
or equivalently
\item[(ii)]
\BEQ \proba[\ln\left(\frac{S_n}{\eps}\right)>B]\lesssim e^{-cB^2}; \label{eq:cor:Nag02}\EEQ
\item[(iii)]
\BEQ \proba[\ln S_n>A]\lesssim (A/\eps)^{-c\ln(A/\eps)}. \label{eq:cor:Nag03}\EEQ 
\end{itemize}
Furthermore, the same estimates (\ref{eq:cor:Nag01}), (\ref{eq:cor:Nag02}), (\ref{eq:cor:Nag03}) still hold
if one replaces $S_n$ by $\tilde{S}_n$. 
\end{Corollary}

{\bf Proof.}

Note that (ii) is equivalent to (i), and (iii) follows directly from (ii) since $e^A\gg n\eps$ and
$\frac{e^A}{\eps}\gg \frac{A}{\eps}$ if $A\gg n\eps$. Also, since $\esper[e^{\eps |Z_i|}]=1+O(\eps)$,
$\tilde{S}_n-S_n=O(n\eps)\ll A$, so the same estimates hold indifferently for $S_n$ or $\tilde{S}_n$ (up to
the choice of $c$).

Therefore we need only prove (i) for $S_n$. 
We apply the above Proposition with $X=e^{\eps|Z|}-\esper[e^{ \eps|Z|}]$,
where $Z$ is any of the variables $Z_i$. 
 One finds  
$\esper X^2=\esper[e^{2\eps|Z|}]-(\esper[e^{\eps|Z|}])^2\approx \eps^2$ and
$\esper[X^t {\bf 1}_{X>0}]<\esper[e^{ t \eps |Z|}]\le 2e^{t^2 \eps^2/2}.$ The bound (\ref{eq:Nag})
is close to optimal if one chooses $t=\half \ln(A/\eps)\gg 1$; we then find (using $\ln^2(A/\eps)\ll (A/\eps)^{\kappa}$
for all $\kappa>0$)
\BEQ \proba[S_n> A]\lesssim ne^{\frac{\eps^2}{8} \ln^2(A/\eps)} A^{-\half \ln(A/\eps)}+
e^{-\frac{c}{n} (A/\eps)^{\frac{3}{2}-\kappa}}. \label{eq:Nag02}\EEQ

\medskip

\begin{itemize}
\item[(i)]
 Assume first that $A\gtrsim e^{1/\eps}\gg 1$. 
The second term in the right-hand side of (\ref{eq:Nag02}) is then the smaller one since (for $\kappa<\half$)
\BEQ e^{-\frac{c}{n} (A/\eps)^{\frac{3}{2}-\kappa}}\le e^{-\frac{c}{n} (A/\eps)  (A/\eps)^{\half-\kappa}}
 \le  e^{-\half  \ln(A/\eps)^2}= (A/\eps)^{-\half \ln(A/\eps)} \le A^{-\half\ln(A/\eps)}.  \label{eq:6.25} \EEQ
As for the first term, it is bounded by $A^{-c\ln(A/\eps)}$ since (using $A\gg 1$)
\BEA e^{\frac{\eps^2}{8} \ln^2(A/\eps)} &\le & e^{\frac{\eps}{8} \ln^2(A/\eps)}\le e^{\frac{1}{8}
(\eps \ln(A)+1) \ln(A/\eps)} \nonumber\\
&\le & e^{\frac{1}{3} \ln(A)\ln(A/\eps)}=A^{\frac{1}{3}\ln(A/\eps)} \label{eq:6.36} \EEA
and 
\BEQ n\ll A/\eps=e^{\ln(A/\eps)} \lesssim A^{\frac{1}{12}\ln(A/\eps)}. \EEQ
 All together one has obtained
\BEQ \proba[S_n>A]\lesssim A^{-c\ln(A/\eps)}\lesssim (A/\eps)^{-c'\ln(A/\eps)},\qquad A\gtrsim e^{1/\eps}.
\label{eq:A>>}\EEQ

\item[(ii)] We now assume that $A\lesssim  e^{1/\eps}$, implying that
$t=\half \ln(A/\eps)\lesssim \frac{1}{\eps}$. Then $A^{t}$ is not
necessarily small, so we must first improve our bound on $\esper[X^t {\bf 1}_{X>0}]$:
\BEA \esper[X^t {\bf 1}_{X>0}] &\le & 2\esper[(e^{\eps Z}-1)^t {\bf 1}_{Z>0}] \nonumber\\
&\lesssim & \int_0^{1/\eps} dz (e^{\eps z}-1)^t e^{-z^2/2} + \int_{1/\eps}^{+\infty} dz (e^{\eps z}-1)^t e^{-z^2/2} 
\nonumber\\
&\lesssim & \int_0^{+\infty} dz (\eps z)^t e^{-z^2/2} + \esper[e^{t\eps Z}] e^{-1/2\eps^2} \nonumber\\
&=& 2^{-\half} (\eps\sqrt{2})^t \Gamma(\frac{t}{2}+1) + e^{\half(\eps^2 t^2-\frac{1}{\eps^2})}. \label{eq:6.39} \EEA
Since $t\lesssim \frac{1}{\eps}$, we find 
$e^{\half(\eps^2 t^2-\frac{1}{\eps^2})}\lesssim e^{-c/\eps^2}\ll e^{-\frac{1}{\eps}|\ln\eps|}\lesssim e^{-t|\ln\eps|}
=\eps^t.$ Hence
\BEQ n\esper[X^t {\bf 1}_{X>0}] A^{-t}\lesssim n t^t (A/\eps)^{-c\ln(A/\eps)}\lesssim (A/\eps)^{-c'\ln(A/\eps)}.\EEQ
As for the second term, clearly $e^{-\frac{c}{n} (A/\eps)^{5/4}}\le (A/\eps)^{-c'\ln(A/\eps)}$ (see (\ref{eq:6.25})).  All together,
\BEQ \proba[S_n>A]\lesssim (A/\eps)^{-c\ln(A/\eps)},\qquad A\lesssim e^{1/\eps}.  \label{eq:A>} \EEQ

\end{itemize}
\hfill \eop

{\bf Remark.} The above results are actually valid as soon as $A\gg n^{\kappa}\eps$ with $\kappa>\half$, as the
reader may easily check (choose $t=c\ln(A/\eps)$ with $c$ small enough and see how (\ref{eq:Nag02}) and (\ref{eq:6.25})
are modified). The condition $A\gg n^{\kappa}\eps$ may certainly be further improved with some extra effort.

\medskip

Corollary \ref{cor:Nagaev0} has the following generalization.

\begin{Corollary}[block large deviation estimates] \label{cor:Nagaev}
Let $Z:=\sum_{i'=1}^{n'} |\tilde{Z}_{i'}|$, where $(\tilde{Z}_{i'})_{i'=1,\ldots,n'}$, $n'\gg 1$ are i.i.d. standard
Gaussian variables; $n\in\N^*$ a multiple of $n'$, $Z_i,\ i=1,\ldots,n/n'$ i.i.d. copies of $Z$;  
\BEQ S_n:=\sum_{i=1}^{n/n'} \left( e^{\eps Z_i}-\esper[e^{\eps Z_i}] \right),\qquad \tilde{S}_n:=
\sum_{i=1}^{n/n'} \left( e^{\eps Z_i}-1\right)\EEQ
where $0<\eps\ll 1$ and $\eps n'\ll 1$. Let finally $A\gg n\eps$ and $B\gg \ln(n)$. Then there exists a constant $c>0$ such that
\begin{itemize}
\item[(i)]
\BEQ \proba[S_n>A]\lesssim (A/n'\eps)^{-c\ln(A/n'\eps)} \label{eq:cor:Nag1}\EEQ
or equivalently
\item[(ii)]
\BEQ \proba[\ln\left(\frac{S_n}{n'\eps}\right)>B]\lesssim e^{-cB^2}; \label{eq:cor:Nag2}\EEQ
\item[(iii)]
\BEQ \proba[\ln S_n>A]\lesssim (A/n'\eps)^{-c\ln(A/n'\eps)}. \label{eq:cor:Nag3}\EEQ 
\end{itemize}
Furthermore, the same estimates (\ref{eq:cor:Nag1}), (\ref{eq:cor:Nag2}), (\ref{eq:cor:Nag3}) still hold
if one replaces $S_n$ by $\tilde{S}_n$. 
\end{Corollary}

\medskip

{\bf Proof.} The result is exactly the one stated in Corollary \ref{cor:Nagaev0} if $n'=1$. We want to prove
the same kind of result for blocks of size $n'$. Standard large deviation arguments apply to $Z$, yielding
(see (\ref{eq:standard-LD})) $\proba[Z>A]\le ce^{-\frac{1}{2n'} (A-cn')_+^2}$, hence (letting as before
$X:=e^{\eps Z}-\esper[e^{\eps Z}]$), $\esper[e^{\eps Z}]=1+O(n'\eps)$,
$\esper[X^2]\lesssim \eps^2 \Var(Z)=O(n'\eps^2)$, and
\BEA \esper[X^t {\bf 1}_{X>0}]\le \esper[e^{t\eps Z}] & \lesssim & t\eps \int_0^{+\infty} e^{t\eps z} e^{-\frac{1}{2n'}(z-cn')_+^2} dz\nonumber\\
&\le & \int_0^{2cn'}  t\eps e^{t\eps z} dz +\int_{-\infty}^{+\infty} e^{t\eps z-\frac{1}{8n'}z^2} dz \nonumber\\
&\lesssim & e^{2Cn't\eps}+e^{Cn'(t\eps)^2/2}.\EEA
We set $t:=\half \ln(A/n'\sqrt{\eps})$ and distinguish two regimes according to whether $A\gtrless e^{1/(\eps\sqrt{n'})}$,
corresponding to $t\gtrless \frac{1}{\eps\sqrt{n'}}$. Thus
\BEQ e^{2Cn't\eps}=e^{Cn'\eps \ln(A/\eps\sqrt{n'})}\le c^{\ln(A/\eps\sqrt{n'})}\ll A^{\frac{1}{3}\ln(A/\eps\sqrt{n'})} \EEQ
instead of (\ref{eq:6.36}), and 
\BEA \esper[X^t {\bf 1}_{X>0}] &\le & \esper[(e^{\eps Z}-1)^t] \nonumber\\
&\lesssim & \int_0^{1/\eps} dz (e^{\eps z}-1)^t e^{-\frac{1}{2n'}(z-cn')^2_+} + \int_{1/\eps}^{+\infty}
dz (e^{\eps z}-1)^t e^{-\frac{1}{2n'}(z-cn')^2_+} \nonumber\\
&\lesssim & \int_0^{2cn'} dz (\eps z)^t + \int_0^{+\infty} dz (\eps z)^t e^{-z^2/8n'} +
\int_{1/\eps}^{+\infty} dz\, e^{t\eps z} \,  e^{-\frac{z^2}{8n'}-\frac{c}{n'\eps^2}} \nonumber\\
&\lesssim & (2\eps Cn')^{t+1}+ (Cn')^{\half(t+1)}\eps^t \Gamma(\frac{t}{2}+1) + e^{c(n'\eps^2 t^2-\frac{1}{n'\eps^2})}
\EEA
instead of (\ref{eq:6.39}).  Hence all estimates contained in the proof of Corollary \ref{cor:Nagaev0} hold if one
replaces $A/\eps$ by $A/n'\eps$ or $A/\sqrt{n'}\eps$.

\hfill\eop


\subsubsection{Mayer expansion}


We also  need in the course of the proof of Theorem \ref{th:LD} a technical result which we choose to
state separately for the sake of clarity. In the sequel, $B^j$ is one of the
$\mu$-dependent subsets of cubes with $\tilde{B}^j(x,\rho)\subset B^j\subset \tilde{B}^j(x,\rho)\cup \partial\tilde{B}^j(x,\rho)$ introduced in section 6.2.

\begin{Lemma}["Mayer expansion"]  \label{lem:Mayer0}
Let $(z_{\Del})_{\Del\in \tilde{B}_j(x,\rho)}\in\R_+$ and $c>0$. Define
\BEQ y_{\Del}:=\sum_{\Del'\in \tilde{B}_j(x,\rho)} e^{-c_0 d^j(\Del,\Del')} z_{\Del'} \EEQ
(see eq. (\ref{eq:6.27})).
\begin{itemize}
\item[(i)] ("Mayer expansion")
Let
\BEQ S_{\del}((z_{\Del})):=\sum_{\Del\in \tilde{B}_j(x,\rho)}  \left( e^{e^{-c_0\del}\eps z_{\Del}} -1\right), \qquad
\del\ge 0 \label{eq:6.60} \EEQ
and 
\BEQ S_0((y_{\Del})):=\sum_{\Del\in B_j} \left(e^{\eps y_{\Del}}-1\right).\EEQ
Then
\BEQ 0\le  S_0((y_{\Del}))\le \sum_{m\ge 1} \frac{1}{(m-1)!} \sum_{\del_1<\ldots<\del_m} \prod_{p=1}^m S_{\del_p}
((z_{\Del})) \label{eq:Mayer01}\EEQ
where $\del_i,\ i=1,2,\ldots$ range among the set $\{d^j(\Del,\Del'),\Del,\Del'\in B^j\}$.

More generally, if $\del_{max}\in\R_+$, then
\BEA &&  S_0((y_{\Del}))\le \sum_{m\ge 1} \frac{1}{(m-1)!} \sum_{\del_{max}<\del_1<\ldots<\del_m} \prod_{p=1}^m S_{\del_p}
((z_{\Del})) \nonumber\\
&& \qquad +\sum_{m\ge 1} \sum_{m'=1}^{m-1} \frac{1}{(m-1-m')!} \sum_{ \del_1<\ldots<\del_{m'}\le\del_{max}< \del_{m'+1}<\ldots<\del_m}  S_0
((z_{\Del}))^{1+O(e^{-c_0})} \prod_{p=m'+1}^{m} S_{\del_p}((z_{\Del})) \nonumber\\
&& \qquad +\sum_{m\ge 1}  \sum_{\del_1<\ldots<\del_m\le\del_{max}}  S_0
((z_{\Del}))^{O(e^{-c\del_1})}. \label{eq:Mayer0-generalized} \EEA

\item[(ii)]
Let
\BEQ T((z_{\Del})):=\sum_{\Del\in \tilde{B}_j(x,\rho)}   e^{\eps z_{\Del}}\EEQ
and similarly
\BEQ T((y_{\Del})):=\sum_{\Del\in B_j}   e^{\eps y_{\Del}}.\EEQ
Then 
\BEQ T((y_{\Del}))\le \left(T((z_{\Del}))\right)^{1+O(e^{-c})}. \label{eq:Mayer2} \EEQ
\end{itemize}
\end{Lemma}

\medskip

Note that (by invariance by translation) $\del_i\in \{d^j(\Del,\Del'),\Del'\in\D^j\}=\{0<d_1<d_2<\ldots\}$
where $\Del$ is some arbitrary fixed cube, and $d_i\approx_{i\to\infty} i^{1/d}$. The indexation is easier if we
 choose a distorted distance instead of the Euclidean distance, i.e., if we define e.g. $|x-y|=\sup_{i=1,\ldots,d} |({\cal R}(x-y))_i|$
where ${\cal R}$ is a generic rotation, so that the set $\{\Del'\in \D^j\ |\ d^j(\Del,\Del')
=\del\}$ contains at most one element for $\Del$, $\del$ fixed, which we denote by $\Del(\del)$. The nickname
"Mayer expansion" refers to a common expansion of the free energy in equilibrium statistical physics where $e^{\beta{\cal H}}$,
${\cal H}$ being the local energy density, is expanded
into $\left(e^{\beta{\cal H}}-1\right)+1$, which is exactly what we do in (i).

\medskip
{\bf Proof.}

\begin{itemize}
\item[(i)] 
Let $a_{\Del}(\del):=e^{e^{-c_0\del}\eps z_{\Del}}-1$ $(\Del\in\tilde{B}^j(x,\rho))$    and $a_{\Del}:=a_{\Del}(0)=e^{\eps z_{\Del}}-1$.
Rewriting $S_0((y_{\Del}))$ in terms of the $z$-variables and expanding  the product of terms of the form $a_{\Del(\del)}(\del)+1$ yields
\BEA  &&S_0((y_{\Del})) =\sum_{\Del\in B^j}  \left[ \prod_{\del} (a_{\Del(\del)}(\del)+1) \right]-1= \sum_{\Del\in B^j} \sum_{m=1}^n \sum_{\del_1<\del_2<\ldots 
< \del_m} \prod_{p=1}^m a_{\Del(\del_p)}(\del_p) \nonumber\\
 && \le \sum_{m=1}^n \frac{1}{(m-1)!} \sum_{\del_1<\del_2<\ldots<\del_m} \sum_{\Del\in B^j}
a_{\Del(\del_1)}(\del_1) \nonumber\\ && \qquad \qquad   \left( \sum_{\Del_2\in \{\Del(\del_2),\ldots,\Del(\del_m)\}} a_{\Del_2}(\del_2) \left( 
\sum_{\Del_3\in \{\Del(\del_2),\ldots,\Del(\del_m)\}\setminus \Del_2 } a_{\Del_3}(\del_3) \left(
\cdots \right)\right)
\right)  \label{eq:Mayer-entretemps} \\
&\le &  \sum_{m=1}^n \frac{1}{(m-1)!} \sum_{\del_1<\del_2<\ldots<\del_m} \left( \sum_{\Del\in B^j}
a_{\Del}(\del_1) \right) \left( \sum_{\Del_2\in B^j} a_{\Del_2}(\del_2) \left( \cdots \right)
\right) \nonumber\\
&=& \sum_{m\ge 1} \frac{1}{(m-1)!} \sum_{\del_1<\ldots<\del_m} \prod_{p=1}^m S_{\del_p}((z_{\Del})).
\label{eq:6.52} \EEA

For the proof of (\ref{eq:Mayer0-generalized}), we fix $\del_{max}\ge 0$, start from (\ref{eq:Mayer-entretemps}) and
pick its $p$-th factor,
$A_p=\sum_{\Del_p\in\{\Del(\del_2),\ldots,\Del(\del_m)\}\setminus\{\Del_2,\ldots,\Del_{p-1}\} }
a_{\Del_p}(\del_p)$. If $\del_p>\del_{max}$ we bound $A_p$ by $S_{\del_p}((z_{\Del}))$ as before. Otherwise we use the
identity $(x-1)^{\kappa}\le x^{\kappa}-1$ ($x\ge 1,\kappa\ge 1$) and
H\"older's inequality to get  $a_{\Del_p}(\del_p)\le \left(e^{\eps z_{\Del_p}}-1\right)^{e^{-c_0\del_p}}=a_{\Del_p}^{e^{-c_0\Del_p}}$ and
\BEQ A_p\le (m-p+1) \left(\sum_{\Del\in B^j} a_{\Del}\right)^{e^{-c\del_p}}=(m-p+1)S_0((z_{\Del}))^{e^{-c_0\del_p}}.\EEQ
Finally, $\sum_p e^{-c_0\del_p}\le 1+e^{-c_0 d_1}+e^{-c_0 d_2}+\ldots=1+O(e^{-c_0}).$

\item[(ii)]

One finds (all sums or supremums in the next expressions range over subsets of $\tilde{B}^j(x,\rho)$, unless
otherwise stated)
\BEA \sum_{\Del\in B^j} e^{\eps y_{\Del}} & \le &\sum_{\Del\in B^j}
\prod_{\Del'} e^{e^{-c_0 d^j(\Del,\Del')} \eps z_{\Del'}}  \nonumber\\
&\le & \frac{1}{(n-1)!} \sum_{\Del_1} e^{\eps z_{\Del_1}} \left( \sum_{\Del_2\not=\Del_1} e^{e^{-c_0 d_1}\eps z_{\Del_2}}
\left( \sum_{\Del_3\not=\Del_1,\Del_2} e^{e^{-c_0 d_2} \eps z_{\Del_3}} \left(\cdots\right)\right)\right) \nonumber\\
&\le & \frac{1}{(n-1)!}  \left( \sum_{\Del_1} e^{\eps z_{\Del_1}}\right) \sup_{\Del_1} \left( \sum_{\Del_2\not=
\Del_1} e^{e^{-c_0 d_1}\eps z_{\Del_2}}\right) \cdots \nonumber\\
&\le & \frac{1}{(n-1)!}  \left( \sum_{\Del_1} e^{\eps z_{\Del_1}}\right) \cdot \left[ (n-1)
\left(\sum_{\Del_2} e^{\eps z_{\Del_1}}\right)^{e^{-c_0 d_1}} \right] \cdots
\EEA
(H\"older's inequality was used in the last line). 
The product of the prefactors in the last expression, $(n-1)(n-2)\cdots$ is exactly compensated by the factorial $\frac{1}{(n-1)!}$,
and there remains 
\BEQ \sum_{\Del\in B^j} e^{\eps y_{\Del}} \le \left(\sum_{\Del\in \tilde{B}^j(x,\rho)}
e^{\eps z_{\Del}} \right)^{1+O(e^{-c_0})}. \EEQ

\end{itemize}
\hfill \eop

\medskip

Again, this lemma has a block generalization. Roughly speaking, we  want to group together all cubes $\Del'$ at
distance $\del\approx 3^k$ of a given cube $\Del$ and sum over $k$, instead of summing over the $\del_i$'s which (as a
detailed computation proves)
increase too slowly with $i$ to give a converging series. Actually, we bother to do so only for
$\del>\del_{max}$, in a region where the exponential decay governs essentially the estimates; the value of 
$\del_{max}$ is fixed later in the text. In order to avoid blocks with "holes" and overlaps between
blocks, we introduce the following definitions. Let $\D^{j,k}$, $k\ge \log_3 \del_{max}$ be the set of blocks ${\bf \Del}^k=[3^k 2^{j/2}k_1,3^k 2^{j/2}(k_1+1)]\times\ldots
\times[3^k 2^{j/2}k_d,3^k 2^{j/2}(k_d+1)]$ of size $3^k$ included in $\tilde{B}^j(x,\rho)$. The $3^d-1$ blocks of size $3^k$, $[x_{\Del,1}+\eps_1 3^k 2^{j/2}(k_1-\half),
x_{\Del,1}+\eps_1 3^k 2^{j/2}(k_1+\half)]\times\ldots  [x_{\Del,d}+\eps_d 3^k 2^{j/2}(k_d-\half) x_{\Del,d}+\eps_d 3^k 2^{j/2}(k_d+\half) ]$, where $\vec{\eps}=(\eps_1,\ldots,\eps_d)\in\{-1,0,1\}^d\setminus\{0,\ldots,0\}$, are all situated at
a scaled distance  $\ge\del=3^k$ of $\Del$. We denote them by $\vec{\Del}(\vec{\del})$, where  $\vec{\del}:=(\del,\vec{\eps})$
is a composite index including both  the distance $\del$ and a discrete index $\vec{\eps}$ ranging in a fixed finite set.
Then, for smaller distances $\del<3^k$, we set $\vec{\Del}(\del)=\Del(\del)$ as in the previous lemma, and 
$\vec{\del}=\del$ simply.
All together the blocks $(\vec{\Del}(\vec{\del}))_{\vec{\del}}$, $\vec{\del}=(\del,\vec{\eps})$ ($\del\ge 3^k$) or
$\del$ ($\del<3^k$) define for every fixed cube $\Del$ a partition of $\R^d$.
We choose in the sequel some arbitrary total ordering $<$ of the  indices $\vec{\del}$ such that 
$(\vec{\del}=(\del,\vec{\eps})$  or $\del, \ \vec{\del}'=(\del',\vec{\eps}')$ or $\del',\ \del<\del')\Rightarrow 
\vec{\del}<\vec{\del'}.$

\begin{Lemma}[block "Mayer expansion"]  \label{lem:Mayer}
Let $(z_{\Del})_{\Del\in B_j}\in\R_+$ and $c>0$. Define as in the previous lemma
\BEQ y_{\Del}:=\sum_{\Del'\in \tilde{B}^j(x,\rho)} e^{-c_0 d^j(\Del,\Del')} z_{\Del'},\quad S_0((y_{\Del})):=
\sum_{\Del\in B^j} \left(e^{\eps y_{\Del}}-1 \right),\quad S_0((z_{\Del})):=
\sum_{\Del\in \tilde{B}^j(x,\rho)} \left(e^{\eps z_{\Del}}-1 \right) \EEQ and
let, for $k\in\N$,
\BEQ \vec{S}_{3^k}((z_{\Del})):=\sum_{\Del^k\in\D^{j,k}}  \left( e^{e^{-c_0 3^k}\eps \sum_{\Del\in\Del^k\cap B^j} z_{\Del}} -1\right), \qquad
k \ge 0, \EEQ
a block version of (\ref{eq:6.60}) distinguished by the boldface letter. Choose some value of  $\del_{max}$ and order the 
indices $\vec{\del}$ as indicated above.
Then

\BEA &&  S_0((y_{\Del}))\le \sum_{m\ge 1} \frac{1}{(m-1)!} \sum_{\del_{max}<\vec{\del}_1<\ldots<\vec{\del}_m} \prod_{p=1}^m \vec{S}_{\del_p}
((z_{\Del})) \nonumber\\
&& \qquad +\sum_{m\ge 1} \sum_{k=1}^{m-1} \sum_{ \del_1<\ldots<\del_{m'}\le\del_{max}< \vec{\del}_{m'+1}<\ldots<
\vec{\del}_m}  S_0
((z_{\Del}))^{1+O(e^{-c_0})} \prod_{p=1}^k \vec{S}_{\del_p}((z_{\Del})) \nonumber\\
&& \qquad +\sum_{m\ge 1}  \sum_{\del_1<\ldots<\del_m\le \del_{max}}  S_0
((z_{\Del}))^{1+O(e^{-c_0})}. \label{eq:Mayer-generalized} \EEA

\end{Lemma}

{\bf Proof.} If $\vec{\Del}$ is a block of size $3^k$ for some $k\ge 0 $, we let $a_{\vec{\Del}}(\del):=
e^{e^{-c\del}\eps \sum_{\Del\in \vec{\Del}} z_{\Del}}-1$. Thus
\BEA  S_0((y_{\Del})) & \le &\sum_{\Del\in B^j} \left[ \prod_{\vec{\del}} a_{\vec{\Del}(\vec{\del})}(\del)+1\right]-1 \nonumber\\
&=& \sum_{\Del\in B^j} \sum_{m\ge 1} \sum_{\vec{\del}_1>\vec{\del}_2>\ldots>\vec{\del}_m} \prod_{p=1}^m
a_{\vec{\Del}(\vec{\del}_p)}(\del_p).\EEA
We then expand as in (\ref{eq:Mayer-entretemps}) and (\ref{eq:6.52}), and forget the unnecessary factorials
in the denominator (which would require a short discussion in any case since there is no symmetry factor for terms
belonging to blocks with different sizes). \hfill \eop

\hfill \eop



\subsubsection{Proof of Theorem \ref{th:LD}}


We shall  use several times the following elementary lemma.

\begin{Lemma} \label{lem:elementary}
Let $X_1,\ldots,X_n$, $n\ge 1$ be real-valued random variables, and $\lambda_1,\ldots,\lambda_n\in\R_+$ such
that $\lambda_1+\ldots+\lambda_n=1$. Then
\BEQ \proba[\lambda_1 X_1+\ldots+\lambda_n X_n>A]\le \proba[\sup_{p=1,\ldots,n} X_p>A]\le \sum_{p=1}^n \proba[X_p>A].\EEQ
\end{Lemma}

{\bf Proof of Theorem \ref{th:LD}.} 
The general scheme of the proof is the same as that of  Lemma \ref{lem:LD-main}. 
Applying Proposition \ref{prop:Slepian} with $\phi(\eta^j_{t_0})=\psi\left( \frac{1}{\eps}
\fint_{B(x,2^{j/2}\rho)\cap\D^j_{\mu}}dy \left(e^{\lambda 2^j |\eta^j_{t_0}(y)|}-1\right) \right)$ where $\psi:\R_+\to\R$ is any convex, increasing function yields instead of (\ref{eq:6.21})
\BEQ \esper\left[\psi\left( \frac{1}{\eps}
\fint_{B(x,2^{j/2}\rho)\cap\D^j_{\mu}}dy \left( e^{\lambda 2^j |\eta^j_{t_0}(y)|}-1\right) \right)\right] \le 
\esper\left[\psi\left( \frac{1}{\eps}
\fint_{B(x,2^{j/2}\rho)\cap\D^j_{\mu}}dy  \left(e^{\lambda 2^j |\tilde{\eta}^j(y)|}-1\right) \right)\right].
\nonumber\\
  \EEQ

Then we bound the last integral by sums of local suprema $\tilde{Y}_{\Del}=\sup_{\Del}|\tilde{\eta}^j(y)|$,
\BEQ \fint_{B(x,2^{j/2}\rho)\cap \D^j_{\mu}} dy  e^{\lambda 2^j |\tilde{\eta}^j(y)|} \le  \max_{B^j}  \frac{\sum_{\Del\in B^j} e^{\lambda 2^j \tilde{Y}_{\Del}} }{\sharp B_j} 
 \label{eq:6.25bis}
\EEQ
where $B^j$ is a union of cubes ranging over a finite set as in subsections 6.1 and 6.2, and $\sharp B_j\le
\sharp \tilde{B}^j(x,\rho)=n$. 
Eq. (\ref{eq:6.27}), (\ref{eq:coupling}) imply 
\BEQ \lambda 2^j \tilde{Y}_{\Del}\le C\eps \sum_{\Del'\in\tilde{B}^j(x,\rho)} e^{-c_0 d^j(\Del,\Del')} (|Z_{\Del'}|+1)
\le C'\eps + C'\eps \sum_{\Del'\in \tilde{B}^j(x,\rho)} e^{-c_0 d^j(\Del,\Del')} |Z_{\Del'}|. \label{eq:6.79} \EEQ

Finally, we use the formula
\BEA  \esper\left[\psi\left( \frac{1}{n\eps} \sum_{\Del\in B^j} \left(e^{\lambda 2^j \tilde{Y}_{\Del}}-1\right)\right)\right]
& =& \int_0^{+\infty} dA \,  \psi'(A) \, \proba[\sum_{\Del\in B^j} e^{\lambda 2^j \tilde{Y}_{\Del}}>n(1+\eps A)] +\psi(0) \EEA
Clearly, $e^{-C'\eps}(1+\eps A)\ge (1-C'\eps)(1+\eps A)\ge 1+c\eps A$ if $A\gg 1$. Hence, assuming
$\supp\psi'\in [C,+\infty]$ with $C$ large enough, one finds
\BEQ \esper\left[\psi\left( \frac{1}{n\eps} \sum_{\Del\in B^j} \left(e^{\lambda 2^j \tilde{Y}_{\Del}}-1\right)\right)\right]
\le  \int_0^{+\infty} dA \, \psi'(A) \,  \proba[\sum_{\Del\in B^j} e^{c\eps Y_{\Del}}>n(1+\eps A)]+\psi(0) \EEQ
with
\BEQ Y_{\Del}:=\sum_{\Del'\in \tilde{B}^j(x,\rho)} e^{-c_0 d^j(\Del,\Del')} |Z_{\Del'}|\EEQ
(see (\ref{eq:coupling}) and (\ref{eq:6.79})).

Below we prove that 
\BEQ \proba[ \sum_{\Del\in B^j} e^{\eps Y_{\Del}}
 >n(1+\eps A)] \le e^{-c\ln^2(nA)}, \qquad A\gg 1. \label{eq:below}\EEQ
 Hence $\proba[ \sum_{\Del\in B^j} e^{\eps Y_{\Del}}
 >n(1+\eps A)] \lesssim \psi^{-2}(A)$ with
 \BEQ \psi(A):=e^{c\ln^2(n(A-C))} {\bf 1}_{A>C(1+\frac{1}{n})}+e^{c\ln^2(C)} {\bf 1}_{A<C(1+\frac{1}{n})}.\EEQ
 Therefore,
\BEQ\esper\left[\psi\left( \frac{1}{n\eps}\sum_{\Del\in B^j} \left( e^{\lambda 2^j \tilde{Y}_{\Del}} - 1 \right) \right)\right] 
\lesssim [-\frac{1}{\psi(A)}]_{0}^{+\infty}+\psi(0)=O(1)\EEQ
and, by Markov's inequality,
\BEQ \proba[\frac{1}{n} \sum_{\Del\in B^j} e^{\lambda 2^j \tilde{Y}_{\Del}}>1+\eps A] \lesssim \frac{1}{\psi(A)} \lesssim e^{-c\ln^2(nA)}, \qquad A\gg 1 \EEQ
so
\BEQ \proba\left[ \ln\left(\frac{1}{n} \sum_{\Del\in B^j} e^{\lambda 2^j \tilde{Y}_{\Del}} \right)>\eps A\right]=
\proba\left[ \frac{1}{n} \sum_{\Del\in B^j} e^{\lambda 2^j \tilde{Y}_{\Del}}>e^{\eps A} \right]\lesssim e^{-c\ln^2(nA)}.\EEQ
(Note that,  if $\eps A\gg 1$, one obtains in fact a Gaussian queue distribution,
\BEQ \proba\left[ \ln\left(\frac{1}{n} \sum_{\Del\in B^j} e^{\lambda 2^j \tilde{Y}_{\Del}} \right)>\eps A\right]\lesssim
e^{-c(\eps A+\ln \frac{n}{\eps})^2}\le e^{-c(\eps A)^2}, \EEQ
with a very bad coefficient $\eps$ in front of $A$ however.)
Thus
\BEQ \proba\left[\ln (e^{\lambda 2^j |\eta^j_{t_0}|})^*(x)>\eps A\right] \lesssim \sum_{n\ge 1}  e^{-
c\ln^2(nA)}
\lesssim A^{-c\ln(A)}.\EEQ

\vskip 2cm

It remains to prove the key estimate (\ref{eq:below}). Recall $Y_{\Del}=\sum_{\Del'\in \tilde{B}^j(x,\rho)} e^{c_0 d^j(\Del,\Del')} |Z_{\Del'}|$,
where $(Z_{\Del})_{\Del}$ are i.i.d. standard Gaussian variables.
It turns out that there are four different large deviation regimes, according to the value of
$S_0((|Z_{\Del}|)_{\Del})=\sum_{\Del\in B^j} \left(e^{\eps |Z_{\Del}|}-1\right)$,
written as $\tilde{S}_n$ in Corollary \ref{cor:Nagaev},  or
$S_0((Y_{\Del})_{\Del})=\sum_{\Del\in B^j} \left(e^{\eps Y_{\Del}}-1\right)$
(see Lemma \ref{lem:Mayer}). By assumption 
\BEQ \sum_{\Del\in B^j} e^{\eps Y_{\Del}} >   n(1+\eps A) \EEQ
 with $A\gg 1$; in other terms,
 \BEQ S_0((Y_{\Del})_{\Del})> n\eps A.\EEQ

\begin{itemize}
\item[(i)]
(Gaussian regime) Assume $S_0((|Z_{\Del}|)_{\Del})=\tilde{S}_n\le 2$. Then $\eps |Z_{\Del}|=O(1)$ for all $\Del$, hence $\eps Y_{\Del}=O(1)$
too,  so $ S_0((Y_{\Del})_{\Del})
\lesssim S_0((|Z_{\Del}|)_{\Del})\approx \eps \sum_{\Del\in\partial B^j} |Z_{\Del}|$. Therefore
\BEQ   \proba[\sum_{\Del\in B^j} e^{\eps Y_{\Del}} > n(1+\eps A)] \le 
\proba[S_0((Y_{\Del})_{\Del})>n\eps A]\le \proba[\sum_{\Del\in\tilde{B}^j(x,\rho)} |Z_{\Del}|>cnA]\lesssim
e^{-c'n(A-C)^2}. \EEQ
Clearly this last quantity is bounded by $e^{-c''\ln^2(nA)}$ for $A$ large enough.

\item[(ii)] (very large deviation regime) Assume $S_0(Y_{\Del})_{\Del})\gg n^{1+O(e^{-c_0})}$, or equivalently
$A\eps\gg n^{O(e^{-c_0})}$. Then  the "Mayer expansion"
(see Lemma \ref{lem:Mayer0} (i)) is not needed; we use Lemma \ref{lem:Mayer0} (ii) and find successively
$S_0((Y_{\Del})_{\Del})\approx T((Y_{\Del})_{\Del})\le \left( T((|Z_{\Del}|)_{\Del}) \right)^{1+O(e^{-c_0})}$, so $T((|Z_{\Del}|)_{\Del}))\gg n$, hence again $T((|Z_{\Del}|)_{\Del})) \approx
S_0((|Z_{\Del}|)_{\Del}))$. All together we have found $S_0((Y_{\Del})_{\Del}) \lesssim
\left( S_0((|Z_{\Del}|)_{\Del}) \right)^{1+O(e^{-c_0})}$.

So we may apply Corollary \ref{cor:Nagaev0}, to the result that
\BEQ \proba[\sum_{\Del\in B^j} e^{\eps Y_{\Del}} > n (1+\eps A)] \le
\proba[\tilde{S}_n>(n(1+\eps A))^{1/(1+O(e^{-c_0}))}]\lesssim e^{-c''\ln^2(nA)}. \EEQ

\item[(iii)] Assume $n\eps\lesssim 1$ and $A\eps\lesssim n^{O(e^{-c_0})}$. Then we use the generalized block "Mayer expansion" (\ref{eq:Mayer-generalized}) with $\del_{max}=0$:
\BEQ \proba[ S_0((Y_{\Del})_{\Del})> n\eps A] \le \proba\left[\sum_{m\ge 1}  \sum_{ \vec{\del}_1>\ldots>\vec{\del}_m} \prod_{p=1}^m \vec{S}_{\del_p}
((|Z_{\Del}|)_{\Del}) >n\eps A\right].\EEQ
Since $\sum_{m\ge 1}\sum_{\vec{\del}_1>\ldots>\vec{\del}_m} e^{-\frac{c_0}{2}(\del_1+\ldots+\del_m)}
\approx \sum_{m\ge 1} \frac{1}{m!} \sum_{\vec{\del_1},\ldots,\vec{\del}_m} e^{-\frac{c_0}{2}(\del_1+\ldots+\del_m)}
\approx 1$, we get by Lemma \ref{lem:elementary}
\BEA  \proba[ S_0((Y_{\Del})_{\Del})> n\eps A]  &\le &  \sum_{m\ge 1}  \sum_{\vec{\del}_1<\ldots<
\vec{\del}_m}
\proba\left[\prod_{p=1}^m \frac{\vec{S}_{\del_p}((|Z_{\Del}|)_{\Del})}{n\eps e^{-c_0\del_p}}> \prod_{p=1}^m
\left(e^{\frac{c_0}{2}\bar{\del}} (n\eps)^{\frac{1}{m}-1} A^{1/m} \right) \right] \nonumber\\
&\le & \sum_{m\ge 1}\sum_{\vec{\del}_1<\ldots<\vec{\del}_m} \sum_{p=1}^m 
\proba\left[\frac{\vec{S}_{\del_p}((|Z_{\Del}|)_{\Del})}{n\eps e^{-c_0\del_p}}> 
e^{\frac{c_0}{2}\bar{\del}} (n\eps)^{\frac{1}{m}-1}  A^{1/m}\right]
 \nonumber\\ \EEA
 where $\bar{\del}:=\frac{1}{m}(\del_1+\ldots+\del_m)$. Note that the expression $\vec{S}_{\del_p}$ is a sum
 over blocks of size $n'\approx \del_p^d\ll (\eps e^{-c_0\del_p})^{-1}$, and $e^{\frac{c_0}{2}\bar{\del}} (n\eps)^{\frac{1}{m}-1}  A^{1/m} \gg 1$, hence we are in the large deviation regime studied in Corollary \ref{cor:Nagaev}. Also,
 \BEQ n e^{\frac{c_0}{2}\bar{\del}} (n\eps)^{\frac{1}{m}-1}  A^{1/m}\ge ne^{\frac{c_0}{2}\bar{\del}}(n\eps)^{-\half}
 \ge e^{\frac{c_0}{2}\bar{\del}} A^{\half} n^{\half-O(e^{-c_0})}  \label{eq:6.94} \EEQ
 if $m\ge 2$. For $m=1$ the estimates of Corollary \ref{cor:Nagaev} give directly a log-normal queue, so we
 sum over $m\ge 2$. By construction $\del_m\ge 3^{m/(3^d-1)}$, so $\bar{\del}>\frac{\del_m}{m}\gtrsim \sqrt{\del_m}$
 and $n'\approx \del_p^d\le \del_m^d\ll e^{\frac{c_0}{2}\bar{\del}}$. Hence 
 \BEQ \proba\left[\frac{\vec{S}_{\del_p}((|Z_{\Del}|)_{\Del})}{\eps e^{-c_0\del_p}}> 
n e^{\frac{c_0}{2}\bar{\del}} (n\eps)^{\frac{1}{m}-1}  A^{1/m}\right] \lesssim e^{-c\ln^2 \left(\frac{n}{n'} e^{c_0\bar{\del}}
A \right)} \lesssim e^{-c\ln^2 \left(n e^{\frac{c_0}{2}\bar{\del}}
A \right)} \lesssim e^{-c''\bar{\del}^2} (nA)^{-c''\ln(nA)}.\EEQ

Let $V_m(r):=\sharp\{(\del_1,\ldots,\del_m)\ |\ \bar{\del}<r\}=\sharp\{(\del_1,\ldots,\del_m)\ |\ \sum_{p=1}^m \del_p<mr\}$: clearly $V_m(r)\lesssim \sharp\{(i_1,\ldots,i_m)\in \N^m\ |\ i_1+\ldots+i_m=mr\}$, hence  
\BEQ  V_m(r) \lesssim  \int_0^{+\infty} dx_1\ldots dx_m {\bf 1}_{\sum_p x_p<mr} =\frac{(mr)^m}{m!}\lesssim (Cr)^m.\EEQ

Thus  (with an extra factor $\frac{1}{m!}$ due to the ordering $\vec{\del}_1>\ldots>\vec{\del}_m$)
\BEA  \sum_{m\ge 2} \sum_{\vec{\del}_1>\ldots>\vec{\del}_m} \proba\left[ \frac{\vec{S}_{\del_p}((|Z_{\Del}|)_{\Del}) }{n\eps e^{-c_0\del_p}}> 
e^{\frac{c_0}{2}\bar{\del}} (n\eps)^{\frac{1}{m}-1}  A^{1/m}\right] & \lesssim&  \sum_m \frac{1}{m!}  (nA)^{-c\ln(nA)}
\sum_{r=0}^{+\infty} (V_m(r)-V_m(r-1))e^{-cr^2} \nonumber\\ &\le&   (nA)^{-c\ln(nA)} 
\sum_{r=0}^{+\infty} e^{-cr^2} \sum_m \frac{V_m(r)}{m!} \nonumber\\
&\lesssim & (nA)^{-c\ln(nA)}  \EEA

\item[(iv)] Finally, assume $n\eps\gtrsim 1$ and $A\eps\lesssim n^{0(e^{-c_0})}$, and apply the generalized "Mayer expansion"
with $\del_{max}=\frac{8}{c_0}\ln(n\eps)$ defined in such a way that $e^{\frac{c_0}{8}\del}\ge n\eps$
for $\del\ge\del_{max}$.
Since 
\BEQ  \sum_{m\ge 1}\sum_{m'=0}^m  \sum_{\del_1<\ldots<\del_{m'}< \del_{max}\le \vec{\del}_{m'+1}<\ldots<\vec{\del}_m} e^{-\frac{c_0}{2}(\del_{m'+1}+\ldots+\del_m)} 
 \lesssim \sharp{\cal P}(\{1,\ldots,\del_{max}\})=2^{\del_{max}}=(n\eps)^{2\ln(2)/c_0}, \label{eq:part} \EEQ
 the generalized "Mayer expansion", together with Lemma \ref{lem:elementary}, yield
\BEA &&  \proba[ S_0((Y_{\Del})_{\Del})> n\eps A]  \le
  \sum_{m\ge 1} \sum_{m'=0}^m    \sum_{\del_1<\ldots<\del_{m'}\le \del_{max}< \vec{\del}_{m'+1}<\ldots<\vec{\del}_m}
  \nonumber\\ && \qquad \qquad  
\proba\left[ S_0((|Z_{\Del}|)_{\Del})^{1+O(e^{-c_0})} \prod_{p=m'+1}^{m}
 \frac{\vec{S}_{\del_p}((|Z_{\Del}|)_{\Del})}{n\eps e^{-c_0\del_p}} > (n\eps)^{1-\frac{2\ln(2)}{c_0}}  A \prod_{p=m'+1}^m \frac{e^{\frac{c_0}{2}\bar{\del}}}{n\eps} \right] \nonumber\\
&&\le \sum_{m\ge 1} \left\{ \sum_{\del_1<\ldots<\del_m\le\del_{max}} \proba\left[\frac{S_0((|Z_{\Del}|)_{\Del})^{1+O(e^{-c_0})}}{\eps}>\eps^{-2\ln(2)/c_0} n^{1-\frac{2\ln(2)}{c_0}} A \right]  \right. \label{eq:6.100}\\
&&\left.
+\sum_{m'=0}^{m-1} \sum_{\del_1<\ldots<\del_{m'}\le\del_{max}< \vec{\del}_{m'+1}<\ldots<\vec{\del}_m}
\left(\proba\left[\frac{S_0((|Z_{\Del}|)_{\Del})^{1+O(e^{-c_0})}}{\eps}>\eps^{-2\ln(2)/c_0} n^{1-\frac{2\ln(2)}{c_0}} e^{\frac{c_0}{4}\bar{\del}} A \right]   \right.\right.\nonumber\\ && \left.\left. \qquad \qquad \qquad\qquad +
\sum_{p=m'+1}^m  \proba\left[  \frac{\vec{S}_{\del_p}((|Z_{\Del}|)_{\Del})}{\eps e^{-c_0\del_p}}  >  \frac{e^{\frac{c_0}{4}\bar{\del}}}{\eps} \right] \right)\right\} \label{eq:6.101}
 \EEA
 where $\bar{\del}:=\frac{1}{m-m'} (\del_{m'+1}+\ldots+\del_m)$.

For $c_0$ large enough (recall $c_0$ has been multiplied by $m_0$, thus it suffices to choose $m_0$ large enough)
\BEQ n^{1-\frac{2\ln(2)}{c_0}} \eps^{-2\ln(2)/c_0} e^{\frac{c_0}{4}\bar{\del}} A\gtrsim e^{\frac{c_0}{4}\bar{\del}}
An^{\kappa} \gtrsim n^{\kappa} \EEQ
with $\kappa>\half$, and
\BEQ \frac{e^{\frac{c_0}{4} \bar{\del}}}{\eps}\ge \max\left(n, \frac{e^{\frac{c_0}{6} \bar{\del}} (n\eps)^{2/3}}{\eps}
\right) \gtrsim \max\left(n, e^{\frac{c_0}{6}\bar{\del}} A^{\frac{1}{3}} n^{\frac{2}{3}-O(e^{-c_0})} \right).\EEQ
Thus we are in the large deviation regime (see remark after Corollary \ref{cor:Nagaev0}, and Corollary \ref{cor:Nagaev}),
and the lower bounds are as in (\ref{eq:6.94}), yielding a bound $O((nA)^{-c\ln(nA)})$ for the sum (\ref{eq:6.101})
over $m'$ and $\del_{m'+1},\ldots,\del_m$. As for the first sum  over $\del_1<\ldots<\del_m\le\del_{max}$ in
(\ref{eq:6.100}), or $\del_1<\ldots<\del_{m'}$ in (\ref{eq:6.101}), it produces as in (\ref{eq:part}) a supplementary
multiplicative factor  of order $(n\eps)^{2\ln(2)/c_0}$, of no  incidence on the result since
$(n\eps)^{2\ln(2)/c_0} (nA)^{-c\ln(nA)}\le n^{2\ln(2)/c_0} (nA)^{-c\ln(nA)}
\lesssim (nA)^{-c'\ln(nA)}.$

\end{itemize}

\hfill\eop


\subsection{Proof of Theorem \ref{th:LD1}}


We may now finally prove Theorem \ref{th:LD1}. 
Let 
\BEQ F(\del t,[0,t]):= 2^{-j}\del t \sum_{p=0}^{\lfloor t/\del t\rfloor} (e^{-2^{-j} \del t})^p \ln \left[
\left(e^{\lambda 2^{j} \fint_{t-p\del t}^{t-(p-1)\del t} |\eta^j(s)| ds} \right)^*(x)\right],\EEQ
and 
\BEQ c_{\del t}:=\left[ 2^{-j}\del t \sum_{p=0}^{\lfloor t/\del t\rfloor} (e^{-2^{-j} \del t})^p \right]^{-1}
\gtrsim 1.\EEQ

Clearly,
\BEQ F(\del t,[0,t])\le 2^{-j}\del t \sum_{p=0}^{\lfloor t/\del t\rfloor} (e^{-2^{-j} \del t})^p \ln \left[
\left(e^{\lambda 2^{j} ||\tilde{\eta}^j||_{\infty, [t-p\del t,t-(p-1)\del t]} }  \right)^*(x)\right]\EEQ
where $||\tilde{\eta}^j||_{\infty, [t-p\del t,t-(p-1)\del t]}(x):=\sup_{s\in [t-p\del t,t-(p-1)\del t]}
|\eta^j(s,x)|.$ Thus
\BEA   c_{\del t} F(\del t) & \le & \sup_p (e^{-2^{-j} \del t})^p \ln \left[
\left(e^{\lambda 2^{j} ||\tilde{\eta}^j||_{\infty, [t-p\del t,t-(p-1)\del t]} }  \right)^*(x)\right] \nonumber\\
& \le & \sup_{q\in\N} e^{-q} \ln \left[
\left(e^{\lambda 2^{j}  ||\tilde{\eta}^j||_{\infty, [t-q 2^j,t-(q-1) 2^j]} }  \right)^*(x)\right]
.\EEA

   The estimates we developed for $|\eta^j_{t_0}|$ in Theorem \ref{th:LD}  extend to 
$||\tilde{\eta}^j||_{\infty, [t-q2^j,t-(q-1)2^j]}(x)$ by using the BTIS inequality once again.
Hence 
\BEA \proba[F(\del t,[0,t])>\eps A] & \lesssim& \sum_q \proba\left [
\ln \left[
\left(e^{\lambda 2^{j} ||\tilde{\eta}^j||_{\infty, [t-q2^j,t-(q-1)2^j]} }  \right)^*(x)\right] >c\eps
e^q A\right] \nonumber\\
&\lesssim & \sum_{q=0}^{+\infty} \left( e^q A \right)^{-c' \ln \left( e^q A \right)}
\lesssim A^{-c'' \ln A}\EEA
by Theorem \ref{th:LD}.

Now, by H\"older's inequality,
\BEA  F(\del t,[0,t]) & \le &  \half 2^{-j}\del t \sum_{p=0}^{\lfloor t/\del t\rfloor} (e^{-2^j \del t})^p  \left\{\ln \left[
\left(e^{\lambda 2^j \fint_{t-p\del t}^{t-(p-1/2)\del t} |\eta^j(s)| ds} \right)^*(x)\right]+ \ln \left[
\left(e^{\lambda 2^j \fint_{t-(p-1/2)\del t}^{t-(p-1)\del t} |\eta^j(s)| ds} \right)^*(x)\right] \right\} \nonumber\\
& \le &  e^{2^j \del t/2} F(\del t/2,[0,t])\EEA
so $|||\eta^j|||_{\lambda,j}(\R_+,x)\le \limsup_{\del t\to 0,t\to +\infty} \frac{1}{\lambda} F(\del t,[0,t])$
by (\ref{eq:etaj1}), and by monotone convergence
we get $\proba[|||\eta^j|||_{\lambda,j}(\R_+,x)> A 2^{-jd_{\phi}}]\le \limsup_{\del t\to 0,t\to +\infty} \proba[F(\del t,[0,t])>
\eps A]\lesssim A^{-c\ln A}.$ 

\medskip

The estimates for $||| 2^{j/2} \frac{\eta^j(.,.+\del x)-\eta^j(.,.)}{|\del x|} |||_{\lambda,j}(\R_+,x)$ are proved in the same way since $2^{j/2}\frac{\eta^j(.,.+\del x)-\eta^j(.,.)}{|\del x|}=O((2^{-j})^{1+d_{\phi}})$ scales like 
$\eta^j_t(x)$ (see \ref{eq:5.20bis}).

\hfill \eop


\bigskip
{\bf Acknowledgements.} We wish to express our thanks to S. Benachour and G. Barles for useful discussions.


\end{document}